# On Culler-Shalen seminorms and Dehn filling

By S. Boyer and X. Zhang*

### Introduction

If $\Gamma$ is a finitely generated discrete group and $G$ a complex algebraic Lie group, the $G$-character variety of $\Gamma$ is an affine algebraic variety whose points correspond to characters of representations of $\Gamma$ with values in $G$. Marc Culler and Peter Shalen developed the theory of $SL_2(\mathbf{C})$-character varieties of finitely generated groups and applied their results to study the topology of 3-dimensional manifolds in the papers [6], [7], [8].

Consider the exterior $M$ of a hyperbolic knot lying in a closed, connected, orientable 3-manifold. The Mostow rigidity theorem implies that the holonomy representation $\bar{\rho} : \pi_1(M) \to \text{Isom}_+(\mathbf{H}^3) \cong \text{PSL}_2(\mathbf{C})$ is unique up to conjugation and taking complex conjugates. The orientability of $M$ can be used to show $\bar{\rho}$ lifts to a representation $\rho \in \text{Hom}(\pi_1(M), SL_2(\mathbf{C}))$ whose character determines an essentially unique point of $\chi_\rho$ of $X(\pi_1(M))$, the $SL_2(\mathbf{C})$-character variety of $\pi_1(M)$. Culler and Shalen [8] proved that the component $X_0$ of $X(\pi_1(M))$ which contains $\chi_\rho$ is a curve. One of their major contributions was to show how $X_0$ determines a norm on $H_1(\partial M; \mathbf{R})$ which encodes many topological properties of $M$. In particular it provides information on the Dehn fillings of $M$. Their construction may be applied to arbitrary curves in the $SL_2(\mathbf{C})$-character variety of a connected, compact, orientable, irreducible 3-manifold whose boundary is a torus, though in this generality one can only guarantee that it will define a seminorm. The first half of this paper is devoted to the development of the general theory of Culler-Shalen seminorms defined for curves of $PSL_2(\mathbf{C})$-characters. By working over $PSL_2(\mathbf{C})$ we obtain a theory that is more generally applicable than its $SL_2(\mathbf{C})$ counterpart, while being only mildly more difficult to set up.

In the second half of this paper we apply the theory of Culler-Shalen seminorms to study the Dehn filling operation. In particular we examine the relationship between fillings which yield manifolds having a positive dimen-

*Boyer was partially supported by grants: NSERC OGP 0009446 and FCAR EQ 3518. Zhang was supported by a postdoctoral fellowship from the Centre de Recherches Mathématiques.

sional $PSL_2(\mathbf{C})$-character variety with those that yield manifolds having a finite or cyclic fundamental group. In one interesting application of this work we show that manifolds resulting from a nonintegral surgery on a knot in the 3-sphere tend to have a zero-dimensional $PSL_2(\mathbf{C})$-character variety (Corollary 6.7). As a consequence we obtain an infinite family of closed, orientable, hyperbolic Haken manifolds which have zero-dimensional $PSL_2(\mathbf{C})$-character varieties (Theorem 1.8). It is of interest to compare this result with a theorem of Culler and Shalen [7, §2] which states that if a closed irreducible 3-manifold has a positive dimensional $SL_2(\mathbf{C})$-character variety, then it is Haken.

One of the advantages gained through the introduction of seminorms is that they provide a unified context in which to study the Dehn filling operation. In order to explain this, let $M$ be the exterior of a knot in a closed, connected, orientable 3-manifold. The Dehn fillings of $M$ are parametrized by the $\pm$ pairs of primitive homology classes in $H_1(\partial M)$ (see §1). When $M$ is the exterior of a hyperbolic knot, inequalities were found in [6] and [3] for the Culler-Shalen norm of the classes in $H_1(\partial M)$ associated to Dehn fillings with a cyclic or finite fundamental group. Further it was shown how these inequalities combine with the geometry of the Culler-Shalen norm to determine sharp bounds on the number of such fillings. Suppose now that $M$ is the exterior of an arbitrary knot. In Theorem 6.2 of this paper we show that analogous inequalities hold for the seminorm determined by any curve of $PSL_2(\mathbf{C})$-characters of $\pi_1(M)$ which contains the character of an irreducible representation. Moreover we show how these inequalities combine with the geometry of nontrivial, indefinite seminorms to allow us to compare fillings with a finite or cyclic fundamental group with those that are, for instance, reducible or Seifert-fibered. This method depends on our being able to construct nontrivial, indefinite Culler-Shalen seminorms. We accomplish this in many interesting situations through an analysis of the essential surfaces in $M$ associated to the ideal points of a curve of $PSL_2(\mathbf{C})$-characters (see e.g. Proposition 5.6). Once this is done, these special seminorms are exploited in a way that is related to the methods used in [1] and [9] to study generalized triangle groups. For instance, our proof of Theorem 1.2 follows the same outline as the original proof of Corollary 1.3 [17], which depended on showing that generalized triangle groups are nontrivial. Our proof of this corollary, through Theorem 1.2, reduces to essentially the same fact, though reinterpreted in the language of seminorms.

The paper is organized as follows. Section 1 provides a discussion of the exceptional Dehn filling problem and our applications of Culler-Shalen seminorm theory to it. In Section 2 we compile some technical lemmas which will be used in the later sections of the paper. The theory of $PSL_2(\mathbf{C})$-character varieties and Culler-Shalen seminorms is then developed in Sections 3, 4, 5 and

6. The remainder of the paper is devoted to providing proofs of the results discussed in Section 1.

The authors gratefully acknowledge the many technical and stylistic suggestions made by the referee. They helped us to improve all facets of the text.

## 1. The main applications

In this section we shall describe our applications of the seminorm method to the topology of 3-manifolds, most particularly to the Dehn filling operation.

Let $M$ be a connected, compact, orientable, irreducible 3-manifold such that $\partial M$ is a torus. A *slope* on $\partial M$ is a $\partial M$-isotopy class of unoriented, essential, simple closed curves. A slope $r$ determines a primitive homology class $H_1(\partial M)$, well-defined up to sign, obtained by orienting a representative curve for $r$ and considering the homology class of the associated 1-cycle. Conversely, any primitive element of $H_1(\partial M)$ can be represented by a nonseparating, oriented, simple closed curve. This curve is well-defined up to isotopy, and so corresponds to some slope $r$ on $\partial M$. It will be convenient for us to use the symbol $\alpha(r)$ to represent either of the two homology classes in $H_1(\partial M)$ associated to a slope $r$. In spite of the ambiguity in the choice of sign, most results in which we use this notation will be sign-independent. In the few instances where this is not the case, we shall specify orientations of the slopes under consideration.

The *distance* between two slopes $r_1$ and $r_2$, denoted $\Delta(r_1, r_2)$, is the minimal geometric intersection number amongst curves representing them. It may be calculated from the identity $\Delta(r_1, r_2) = |\alpha(r_1) \cdot \alpha(r_2)|$, the absolute value of the algebraic intersection number between $\alpha(r_1)$ and $\alpha(r_2)$. In particular, if we fix a basis $\{\gamma_1, \gamma_2\}$ for $H_1(\partial M)$ and we write $\alpha(r_i) = p_i\gamma_1 + q_i\gamma_2$, then

$$\Delta(r_1, r_2) = |p_1 q_2 - p_2 q_1|.$$

Fix a slope $r$ on $\partial M$ and let $M(r)$ denote the manifold obtained by attaching a solid torus to $M$ in such a way that the meridional slope on the boundary of the solid torus is identified with $r$. We say that $M(r)$ is the *Dehn filling $M$* along $\partial M$ with slope $r$. A fundamental result of A. Wallace [46] and W. B. R. Lickorish [27] states that each closed, orientable 3-manifold results from filling the exterior of some link in the 3-sphere. Thus a natural approach to 3-manifold topology is to analyze to what extent various aspects of the topology of a manifold $M$, as above, are inherited by the manifolds $M(r)$. For instance, one could try to describe when a closed, essential surface in $M$ becomes inessential in some $M(r)$, or when an irreducible $M$ could produce a reducible $M(r)$. An excellent survey of this topic may be found in [13].

Another example of some importance arises as follows. A 3-manifold $W$ is called *atoroidal* if every incompressible torus in $W$ is parallel into $\partial W$. It is called *simple* if it is irreducible and atoroidal. Thurston [43, Th. 2.6] has shown that the interior of a compact, connected, orientable 3-manifold which is simple and not Seifert-fibered admits a complete hyperbolic metric of finite volume. Moreover, he has shown that if $M$ (as above) is simple and non-Seifert, then all but finitely many of the closed manifolds $M(r)$ admit hyperbolic structures [42]. Call $r$ a *nonhyperbolic* filling slope if $M(r)$ is not a hyperbolic manifold. The nonhyperbolic filling slopes on $\partial M$ include the slopes on $\partial M$ whose associated fillings are either

- manifolds with finite or cyclic fundamental groups, or
- manifolds which are reducible, or
- manifolds which are Seifert-fibered spaces, or
- manifolds which admit an incompressible torus.

Thurston's geometrization conjecture [43] predicts that the remaining slopes yield fillings which are hyperbolic manifolds.

A basic problem then is to describe the set of slopes on a torally bounded manifold $M$ which are exceptional in the sense that they produce nongeneric fillings. An appropriate description should include an upper bound on the number of such slopes as well as a qualitative measure of their relative positions determined by a bound on their mutual distances. We shall refer to a slope $r$ on $\partial M$ as a *cyclic filling slope* if $M(r)$ has a cyclic fundamental group. Similarly we shall refer to a slope $r$ as either a *finite* filling slope, *reducible* filling slope, or *Seifert* filling slope if the filled manifold $M(r)$ is of the specified type.

Precise estimates are known on the number of and on the distance between exceptional filling slopes of certain kinds. Thus if $M$ is not a simple Seifert-fibered manifold, then there are at most three cyclic filling slopes on $\partial M$ and the distance between two such slopes is at most 1 [6]; if $M$ is not Seifert-fibered and is not a union along a torus of a cable space and $I(K)$, the twisted $I$-bundle over the Klein bottle, then there are at most six slopes which are either finite filling slopes or cyclic filling slopes, and the distance between two such slopes is at most five [3]; the distance between two reducible filling slopes on $\partial M$ is at most one (hence there are at most three such fillings) [18]; if $M$ is simple, then the distance between toral filling slopes is at most eight [14]. In another direction, we remark that a manifold which has a cyclic or finite fundamental group, or is reducible, or Seifert-fibered, or contains a non-boundary parallel, incompressible torus, does not admit a negatively curved Riemannian metric. S. Bleiler and C. Hodgson [2] have combined certain cusp volume estimates due to Colin Adams with results of Gromov and Thurston to show that there are at most 24 slopes on the boundary of a simple, non-Seifert manifold $M$ which do not produce fillings which admit Riemannian metrics of strictly negative sectional curvature, and the distance between two such slopes is at most 21. Given the results stated above and empirical evidence, it is clear that the estimates obtained in [2] are not optimal. Indeed it is conjectured that the distance between the nonhyperbolic filling slopes is bounded above by eight, and if we exclude four specific manifolds $M$, the conjectured bound is five [15, Conj. 3, 4]. Our first result shows that this latter bound holds in many instances.

Call a slope on $\partial M$ a *big Seifert* filling slope if the associated filled manifold is Seifert-fibered but does not admit a Seifert structure whose base orbifold is a 2-sphere with three or fewer cone points.

THEOREM 1.1. *Let $M$ be a compact, connected, orientable, simple, non-Seifert 3-manifold with $\partial M$ a torus. Let $r_1$ and $r_2$ be two slopes on $\partial M$ whose associated fillings are either reducible manifolds, big Seifert manifolds, or manifolds with a finite fundamental group. If neither $M(r_1)$ nor $M(r_2)$ is $\mathbf{R}P^3 \# \mathbf{R}P^3$ or a union of two copies of $I(K)$, then $\Delta(r_1, r_2) \leq 5$.*

Theorem 1.1 is a consequence of the more refined results we shall explain below, together with the work of other authors. We remark that though the theorem should be able to be generalized to include the cases where one of $M(r_1)$ or $M(r_2)$ is either $\mathbf{R}P^3 \# \mathbf{R}P^3$ or a union of two copies of $I(K)$, our methods break down in these cases. Indeed, the techniques we use to prove Theorem 1.1 may also be applied to simple Seifert manifolds to prove similar results. The twisted $I$-bundle over the Klein bottle, $I(K)$, is such a space and it has a slope $r_0$ such that $I(K)(r_0) \cong \mathbf{R}P^3 \# \mathbf{R}P^3$, as well as other finite filling slopes whose distances from $r_0$ are arbitrarily large.

We say that $M$ is *cabled*, or more precisely, $M$ is a *cable* on a manifold $M_1$, if $M = C \cup_T M_1$ where $C$ is a cable space [12, §2], $\partial M \subset \partial C$ and $T = \partial C \cap \partial M_1$ is an incompressible torus in $M_1$.

A finite group which is the fundamental group of a 3-manifold must belong to one of the following types [28]:
  - C-type: cyclic groups,
  - D-type: dihedral-type groups,
  - T-type: tetrahedral-type groups,
  - O-type: octahedral-type groups,
  - I-type: icosahedral-type groups,
  - Q-type: quaternionic-type groups.

THEOREM 1.2. *Let $M$ be a compact, connected, orientable, irreducible 3-manifold with $\partial M$ a torus. Assume that $M$ is neither a simple Seifert-fibered manifold nor a cable on $I(K)$. Fix slopes $r_1$ and $r_2$ on $\partial M$ and suppose that $M(r_1)$ is a reducible manifold.*

(1) *If $M(r_2)$ has a cyclic fundamental group, then $\Delta(r_1, r_2) \leq 1$.*

(2) If $r_2$ is a finite filling slope, then $\Delta(r_1, r_2) \leq 5$ unless $M(r_1) = \mathbf{R}P^3 \# \mathbf{R}P^3$ and $\pi_1(M(r_2))$ is a D-type group or a Q-type group.

(3) If $M$ contains an essential torus and $r_2$ is a finite filling slope, then $\Delta(r_1, r_2) \leq 1$ unless

   (i) $M$ is a cable on $I(K)$, the twisted I-bundle over the Klein bottle, or

   (ii) $M$ is a cable on a simple, non-Seifert manifold $M_1$ for which there are slopes $r'_1, r'_2$ on $\partial M_1$ such that $M_1(r'_1) \cong \mathbf{R}P^3 \# \mathbf{R}P^3$, $M_1(r'_2)$ has a D-type or Q-type fundamental group, and $\Delta(r'_1, r'_2) \geq 8$.

Theorem 1.2 (1) is sharp in the sense that the distance 1 can be realized by slopes on manifolds satisfying the hypotheses. If we require $M$ to be a simple, non-Seifert manifold, Example 7.8 shows that part (1) of the theorem is still sharp, though in this case it seems likely that the inequality in part (2) can be improved and the restrictions on $M(r_1)$ and $M(r_2)$ removed. The exclusion of cables on $I(K)$ in part (3) of the theorem is necessary, as the method of Example 9.11 shows, but it seems likely that these manifolds provide the only exceptions to the conclusion of this part of the theorem.

An *integral* slope on the boundary of the exterior of a knot $K \subset S^3$ is any slope of distance 1 from the meridian slope of the knot. In terms of the standard meridian, longitude coordinates for the first homology of the boundary torus, the integral slopes are those representing homology classes with longitude coordinate $\pm 1$. Note that the meridian slope of $K$ is a cyclic filling slope.

COROLLARY 1.3 ([17]). *If $M$ is the exterior of a knot in $S^3$ and if $M(r)$ is a reducible manifold, then $r$ is an integral slope.*

*Proof.* If $M$ is a simple Seifert-fibered manifold then the knot is a torus knot [32, Th. 2] and the corollary follows from [32]. Otherwise by Theorem 1.2 (1), the distance between $r$ and the meridian slope is 1, and thus $r$ is an integral slope. □

We note that it follows from Corollary 1.3 and the reducible surgery theorem [18], that for a knot $K$ in $S^3$, there are at most two reducible *surgeries* on $K$ (that is, fillings of the exterior of $K$) and if there are two, the surgery slopes correspond to successive integers. The *cabling conjecture*, which is still an open problem, states that the only knots in the 3-sphere which admit reducible surgeries are torus knots and *cabled knots* (that is, knots whose exteriors are cabled). Combining Theorem 1.2, the reducible surgery theorem [18], and the cyclic surgery theorem [6] we obtain the following corollary.

COROLLARY 1.4.  *Let $M$ be a compact, connected, orientable, irreducible 3-manifold with $\partial M$ a torus. Suppose further that $M$ is not a simple Seifert-fibered space. If for $i = 1, 2$ $M(r_i)$ is either a reducible manifold or a manifold with a cyclic fundamental group, then $\Delta(r_1, r_2) \leq 1$. Consequently, there are a total of at most three Dehn fillings on $M$ which are either cyclic or reducible.*

Next we consider the case of fillings yielding spaces with the fundamental group of a Seifert-fibered manifold. As all but at most one filling of a Seifert manifold is Seifert, we shall restrict this discussion to manifolds $M$ which are not Seifert.

Suppose that $r_1$ is a slope on $\partial M$ such that $M(r_1)$ has the fundamental group of a Seifert-fibered manifold $W$. How are the finite filling slopes and cyclic filling slopes on $\partial M$ related to $r_1$? Suppose first of all that $W$ admits a Seifert fibration with base orbifold the 2-sphere having no more than three exceptional fibers, and if there are three such fibers, assume that their indices form a platonic triple, i.e. $M(r_1)$ has a fundamental group which is cyclic or finite. It was shown in [3, Ths. 1.1 and 1.2] that as long as $M$ is not a simple Seifert manifold or a cable on $I(K)$,

$$\Delta(r_1, r_2) \leq \begin{cases} 2 & \text{if } M(r_2) \text{ has a cyclic fundamental group,} \\ 5 & \text{if } M(r_2) \text{ has a finite fundamental group.} \end{cases}$$

These inequalities were also shown to be sharp. Our next two theorems deal with most of the remaining cases.

THEOREM 1.5.  *Let $M$ be a compact, connected, orientable, simple, non-Seifert 3-manifold with $\partial M$ a torus. Suppose that $r_1$ is a slope on $\partial M$ such that $M(r_1)$ has the same fundamental group as a Seifert-fibered space which admits no Seifert fibration having base orbifold the 2-sphere with exactly three cone points. Let $r_2$ be a finite or cyclic filling slope on $\partial M$. Then*

(1)  $\Delta(r_1, r_2) \leq 1$ *if $M(r_2)$ has a cyclic fundamental group;*

(2)  $\Delta(r_1, r_2) \leq 5$ *if $M(r_2)$ has a finite fundamental group, unless $M(r_1)$ is either $\mathbf{R}P^3 \# \mathbf{R}P^3$ or a union of two copies of $I(K)$, and $\pi_1(M(r_2))$ is a D-type group nor a Q-type group.*

Part (1) of Theorem 1.5 is sharp (Example 7.8), but it seems likely that part (2) is not. Firstly the upper bound 5 can most likely be improved, and secondly it is probably unnecessary to exclude the possibility that $M(r_1)$ is either $\mathbf{R}P^3 \# \mathbf{R}P^3$ or a union of two copies of $I(K)$ or that $\pi_1(M(r_2))$ is either a D-type group or a Q-type group.

THEOREM 1.6.  *Let $M$ be a non-Seifert, compact, connected, orientable, irreducible 3-manifold with $\partial M$ a torus and suppose that $M$ contains an essential torus. Consider two slopes $r_1$ and $r_2$ on $\partial M$ such that $M(r_1)$ has the* *fundamental group of a Seifert-fibered space and $M(r_2)$ has a finite or cyclic fundamental group.*

(1) *If $\Delta(r_1, r_2) > 1$ then $M$ is a cable on a manifold $M_1$ which is either simple or Seifert-fibered. Furthermore, $M_1$ admits a finite or cyclic filling according to whether $r_2$ is a finite or a cyclic filling slope.*

(2) *Suppose that $\Delta(r_1, r_2) > 1$ where $r_2$ is a cyclic filling slope and $M(r_1)$ has the fundamental group of a Seifert-fibered space which admits no Seifert fibration having base orbifold the 2-sphere with exactly three cone points. Then $M$ is a cable on a Seifert manifold admitting a cyclic filling.*

It appears likely that under the hypotheses of Theorem 1.6, the inequality $\Delta(r_1, r_2) > 1$ implies that $M$ is the union of a cable space and a Seifert-fibered manifold admitting a finite or a cyclic filling slope. Our proof of Theorem 1.6 shows that the only other possibility is for $M$ to be a cable on a simple, non-Seifert manifold $M_1$ for which there are slopes $r'_1, r'_2$ on $\partial M_1$ such that $M(r_1) \cong M_1(r'_1), M(r_2) \cong M_1(r'_2)$, and $\Delta(r'_1, r'_2) \geq 8$. Notice then that Theorem 1.5 constrains the topology of $M(r_1)$ and $M(r_2)$. In Examples 9.11 and 9.12 we show that there are certain manifolds $M$, each a cable on a Seifert-fibered space, which have big Seifert filling slopes $r_1$ and finite filling slopes $r_2$ for which $\Delta(r_1, r_2)$ is arbitrarily large.

Next we consider the case of Seifert surgery on a knot $K$ in the 3-sphere. It is conjectured that if $K$ is neither a torus knot nor a cable on a torus knot, then only surgeries along integral slopes can yield a Seifert-fibered space. The following corollary verifies this conjecture in many instances. We remind the reader that a *simple knot* is a knot whose exterior is simple, and a compact, connected, orientable 3-manifold is called *Haken* if it is irreducible and contains a properly embedded, 2-sided, incompressible surface.

COROLLARY 1.7. *Let $M$ be the exterior of a nontrivial knot $K$ in $S^3$ and $r$ a slope on $\partial M$ such that $M(r)$ is a Seifert-fibered space.*

(1) *If $M(r)$ is Haken, then $r$ is an integral slope.*

(2) *If $K$ is a satellite knot which is not a cable on a simple knot, then $r$ is an integral slope.*

*Proof.* Satellite knots contain a nonboundary parallel, incompressible torus in their exteriors [12, §2] and conversely, the solid torus theorem [37, 4.C.1] implies that knots whose exteriors contain such a torus are satellite knots. Thus part (2) of the corollary applies only to knots whose exteriors are not simple.

Suppose first of all that $M$ admits the structure of a Seifert-fibered space. Then $M$ is a torus knot exterior [32, Th. 2] and so has a Seifert structure whose base orbifold is the 2-disk with exactly two cone points [32]. In particular this implies that $M$ is a simple manifold, and so we need only show that part (1) of the corollary holds. Now it is shown in [32] that a filling of $M$ is either a connected sum of two nontrivial lens spaces or admits the structure of a Seifert-fibered space whose base orbifold is the 2-sphere with at most three cone points. Hence if $M(r)$ is Haken then $H_1(M(r))$ is infinite [24, VI.13], and therefore $r$ is the longitudinal slope of $K$. As this slope is integral, part (1) of the corollary holds.

Next assume that $M$ is simple but not Seifert-fibered. If $M(r)$ admits no Seifert fibration whose base orbifold is the 2-sphere with exactly three cone points, then we may apply Theorem 1.5 (1), taking $r_1 = r$ and $r_2 = \mu$, the meridional slope of $K$, to deduce that $\Delta(r, \mu) = 1$. Thus $r$ is an integral slope. On the other hand, if $M(r)$ is Haken and admits some Seifert fibration whose base orbifold is the 2-sphere with exactly three cone points, then, as observed in the previous paragraph, $r$ must correspond to the longitudinal slope of $K$, and so is integral. We are therefore done in this case too.

We now assume that $M$ admits a nonboundary parallel, incompressible torus. That is, we take $K$ to be a satellite knot. Note that $M$ does not admit a Seifert fibration as otherwise it would be a torus knot exterior and therefore simple. Suppose that $r$ is not an integral slope, so that $\Delta(r, \mu) > 1$ where $\mu$ is the meridional slope of $K$. We may apply Theorem 1.6 (1), with $r_1 = r$ and $r_2 = \mu$ to see that $M$ is a cable on a nontrivial knot exterior $M_1 \subset S^3$ [37, 4.C.1] which is either simple or Seifert-fibered. If $M_1$ is Seifert-fibered, then it is a torus knot exterior, and so in both cases $K$ is a cable on a simple knot. Thus Corollary 1.7 (2) holds. To deduce part (1) we apply Theorem 1.6 (2) to see that if $M(r)$ is Haken then either $M(r)$ admits a Seifert structure whose base orbifold is the 2-sphere with exactly three cone points, or $M_1$ is a torus knot exterior. The former case cannot arise because then $r$ would be the longitudinal slope of $K$ and so $\Delta(r, \mu) = 1$. We shall finish the proof by showing that the latter case does not arise either.

Let $C$ be the cable space $M \setminus \text{int}(M_1)$ and denote by $\phi$ the slope on $\partial M$ represented by a fiber of the Seifert structure on $C$. It follows from [12, Cors. 7.3 and 7.4] that $\partial M_1$ is an essential torus in $M(r)$. Up to isotopy, $\partial M_1$ is either horizontal (i.e. is transverse to the fibers of $M(r)$) or vertical (i.e. consists of a disjoint union of fibers) [24, VI.34]. Now it cannot be horizontal, for then [24, VI.34] implies that $H_1(M(r))$ would surject onto either $\mathbf{Z}$ or $\mathbf{Z}/2 \oplus \mathbf{Z}/2$. Neither are possible because $H_1(M(r))$ is a finite cyclic group (as $\Delta(r, \mu) \geq 2$). Thus we may assume that $\partial M_1$ is vertical in $M(r)$. It then follows from the uniqueness of the Seifert structure on $M_1$ [24, VI.18] that the structure on $M(r)$ extends that on $M_1$. Hence there is an essential annulus, properly embedded in $M(r) \setminus \text{int}(M_1)$ and vertical in $M(r)$, whose boundary consists of two fibers of the Seifert structure on $M_1$. The slope of such a fibre is an integral slope of the knot $K_1$ whose exterior is $M_1$ [32]. Isotope the annulus so that it intersects $C$ in a properly embedded, essential planar surface $P$. By [16, Lemma 3.1], the surfaces of this sort which intersect $\partial M_1$ in an integral slope of $K_1$ also intersect $\partial M$ in a nonempty collection of curves of integral slope. But this is impossible as $P \cap \partial M$ consists of curves of slope $r$ and we have supposed that $\Delta(r, \mu) \geq 2$. This contradiction completes the proof of the final case of Corollary 1.7. □

It follows from [6] and [17] (see Corollary 1.3 also) that for nontorus knots both cyclic filling slopes and reducible filling slopes are integral. Thus for knots in $S^3$ which are neither torus knots nor cables on a torus knot, the problem of determining whether a Seifert filling slope is necessarily integral reduces to the case where the Seifert manifold is non-Haken and has base orbifold the 2-sphere with exactly three cone points.

Part (2) of Corollary 1.7 has been derived independently by Miyazaki and Motegi [29].

*Proof of Theorem* 1.1. Suppose that $M$ is as described in the hypotheses of the theorem. A fundamental result of Gordon and Luecke states that the distance between two reducible filling slopes is no more than 1 [18]. To see that the distance between a reducible filling slope and a big Seifert filling slope is bounded above by 5, we observe first of all that a big Seifert manifold is either reducible or admits an essential torus. This is seen by examining its base orbifold. Thus the distance estimate follows from [18] and [16, Prop. 6.1]. It follows from [18] and [14] that the distance between two big Seifert filling slopes is at most 5. The remaining cases of the theorem which have to be verified follow from Theorem 1.2 and Theorem 1.5 of this paper, and Theorem 1.1 of [3]. □

The final result we wish to discuss is another application of the seminorm method. It is shown in [7] that if the $SL_2(\mathbf{C})$-character variety of a closed, irreducible 3-manifold $W$ has positive dimension, then $W$ contains a closed incompressible surface: i.e., $W$ is a Haken manifold. Culler and Shalen prove this result by showing that the ideal points of a fixed curve of $SL_2(\mathbf{C})$-characters give rise to nontrivial actions of $\pi_1(W)$ on simplicial trees. Actually the same argument holds true for curves of $PSL_2(\mathbf{C})$-characters (see Theorem 4.3), and so $W$ will be Haken if it has a positive dimensional $PSL_2(\mathbf{C})$-character variety. The converse of this is not true in general; i.e., a closed Haken manifold $W$ may not have a positive dimensional $SL_2(\mathbf{C})$- or $PSL_2(\mathbf{C})$-character variety. Examples of such Haken manifolds have been described by K. Motegi in [33], obtained by gluing two torus knot exteriors according to a well-chosen gluing homeomorphism of their boundaries. Clearly none of Motegi's examples is hyperbolic. As a by-product of the techniques developed in this paper, we prove the following result.

THEOREM 1.8. *An infinite family of hyperbolic Haken manifolds may be constructed, each member of which has a 0-dimensional* $\mathrm{PSL}_2(\mathbf{C})$-*character variety.*

## 2. Preliminaries

We work in the smooth category. All manifolds are understood to be orientable unless otherwise specified. By an *essential* surface in a compact 3-manifold, we mean a properly embedded, incompressible surface such that no component of the surface is $\partial$-parallel and no 2-sphere component of the surface bounds a 3-ball. A 3-manifold is called *irreducible* if it does not contain an essential 2-sphere, and *reducible* otherwise.

Suppose that $M$ is a compact 3-manifold such that $\partial M$ contains a torus component $T$. We shall let $M(T; r)$ denote the manifold obtained by Dehn filling $M$ along $T$ with slope $r$, but shall simplify this notation to $M(r)$ when $\partial M = T$. A slope $r$ on $T$ is called a *boundary* slope if there is an essential surface $F$ in $M$ such that $\partial F \cap T$ is a nonempty set of parallel simple closed curves on $T$ of slope $r$. A boundary slope $r$ on $\partial M$ is called a *strict boundary slope* if there is an essential surface $F$ in $M$ which is not the fiber in any representation of $M$ as a fiber bundle over the circle, and such that $\partial F \cap T$ is a nonempty set of parallel simple closed curves on $T$ of slope $r$.

We shall refer to the following result often in this paper.

LEMMA 2.1 ([6, Th. 2.0.3 and Add. 2.0.4]). *Let $M$ be a compact, connected and irreducible 3-manifold such that $\partial M$ is a torus. Assume that the first Betti number of $M$ is 1 and that $r$ is a boundary slope on $\partial M$. Then either*

(1) $M(r)$ *is a Haken manifold, or*

(2) $M(r)$ *is a connected sum of two nontrivial lens spaces, or*

(3) *$M$ contains a closed, essential surface which remains essential in $M(r')$ whenever $\Delta(r, r') > 1$, or*

(4) *$M$ fibers over $S^1$ with fiber a planar surface having boundary slope $r$. Further, if $r$ is a strict boundary slope and this case arises, then case (3) also occurs.*

COROLLARY 2.2. *Let $M$ be a compact, connected and irreducible 3-manifold such that $\partial M$ is a torus. Assume that the first Betti number of $M$ is 1. If $M(r)$ is a reducible manifold, then one of the following three possibilities occurs*:

(1) $M(r) = L(p,s) \# L(q,t)$ *is a connected sum of two lens spaces with $1 < p, q < \infty$; or*

(2) $M$ *contains a closed, essential surface $S$ such that $S$ is compressible in $M(r)$ and is incompressible in $M(r')$ whenever $\Delta(r,r') > 1$; or*

(3) $M(r) = S^1 \times S^2$.

*Proof.* Consider the collection of all properly embedded surfaces in $M$ resulting from the intersection of $M$ and an essential 2-sphere in $M(r)$. As $M$ is irreducible, no such surface has an empty boundary, and so a standard argument shows that any surface from this collection having a minimal number of boundary components is essential in $M$. In particular $r$ is a boundary slope. We may now apply Lemma 2.1 to deduce that either the conclusion of part (1) or the conclusion of part (3) of this corollary occurs, or there is an essential surface $S$ in $M$ which remains incompressible in $M(r')$ whenever $\Delta(r,r') > 1$. The proof that $S$ may be chosen to compress in $M(r)$ proceeds as follows. Assume that $M(r)$ is not of the form described in (1) or (3). Choose a separating, essential surface $F$ in $M$ which has a nonempty boundary of slope $r$ and which, subject to these conditions, has a minimal number of boundary components. If $F$ is nonplanar, we can use [6, Add. 2.2.2] and the remarks that precede it to find the desired surface $S$. When $F$ is planar, we use the argument in the last paragraph of [6, page 285]. □

LEMMA 2.3. *Let $M$ be a compact, connected, irreducible and $\partial$-irreducible 3-manifold such that $\partial M$ is not connected and contains a torus component $T$. Suppose further that $M$ is not homeomorphic to $T \times [0,1]$. Consider two slopes $r_1$ and $r_2$ on $\partial M$ corresponding to curves on $T$.*

(1) *If $M(T;r_1)$ is $\partial$-reducible and $M(T;r_2)$ is a reducible manifold, then $\Delta(r_1, r_2) \leq 1$.*

(2) *If $M(T;r_i)$ is $\partial$-reducible for $i = 1, 2$, then either $\Delta(r_1, r_2) \leq 1$ or there is an essential annulus $A$ in $M$ with one boundary component of $A$ in a component $S \neq T$ of $\partial M$ and the other in $T$. Moreover if $r_0$ is the slope on $T$ determined by $A$, then $S$ is compressible in $M(T;r)$ if and only if $\Delta(r_0, r) \leq 1$.*

*Proof.* (1). According to the main theorem of [38], our hypotheses imply that there is another essential torus $T'$ in the interior of $M$ and a cable space $C \subset M$ for which $\partial C = T \cup T'$. Further, $r_2$ is the slope of the cabling annulus. We may decompose $M$ as $C \cup_{T'} M_0$ where $T'$ is necessarily incompressible in $M_0$. Consideration of a compressing disk for $\partial(M(T; r_1))$ shows that $T'$ is compressible in $C(T; r_1)$. Hence by [12, Lemma 7.2], $\Delta(r_1, r_2) = 1$.

(2). For $i = 1$ and 2, fix a component $S_i \neq T$ of $\partial M$ which compresses in $M(T; r_i)$. It is shown in [6, Th. 2.4.5] that $\Delta(r_1, r_2) \leq 1$ as long as $S_1 \neq S_2$ and so we may assume that $S_1 = S_2$. If we suppose that $\Delta(r_1, r_2) > 1$, then, as Wu has shown, there exists a properly embedded annulus $A$ in $M$ having one boundary component in $S_1$ and the other in $T$ ([47, Th. 1]). If $r_0$ denotes the slope on $T$ determined by $A$, then it follows from [6, Th. 2.4.3] that $S_1$ compresses in $M(T; r)$ if and only if $\Delta(r, r_0) \leq 1$. Moreover, in this situation, $S_1$ is the only component of $\partial M \setminus T$ which can compress in some $M(T; r)$ ([6, Th. 2.4.5]). □

We say that $M$ is a *generalised* 1-*iterated torus knot exterior* if $M$ is a nontrivial cable on a Seifert-fibered space $M_1$ where $M_1$ admits a Seifert structure whose base orbifold is the 2-disc with exactly two cone points.

LEMMA 2.4. *Let $M$ be a compact, connected, irreducible and $\partial$-irreducible 3-manifold whose boundary is a torus and suppose that $M$ admits an infinite cyclic filling slope $r_1$ as well as a finite filling slope $r_2$. If $M$ is not a cable on $I(K)$ then*

$$\Delta(r_1, r_2) \leq \begin{cases} 1 & \text{if } M \text{ contains an essential torus} \\ 2 & \text{if } M \text{ is simple but not Seifert-fibered.} \end{cases}$$

*Proof.* If $M$ is not a generalised 1-iterated torus knot exterior, the inequalities follow from [3, Th. 1.1 (2) and Th. 1.2 (1)]. Suppose then that $M$ is a generalized 1-iterated torus knot exterior which is not a cable on $I(K)$. Then $M = M_1 \cup C$ where $C$ is a cable space of type $(m, n)$, $n \geq 2$, and $M_1$ admits a Seifert structure whose base orbifold is a 2-disk with two cone points of orders $p$ and $q$ satisfying $2 \leq p \leq q$ and $3 \leq q$. If $r_2$ is a cyclic filling slope, then $\Delta(r_1, r_2) \leq 1$ by the cyclic surgery theorem [6] (we note that though this theorem is stated for non-Seifert-fibered manifolds, it holds for nonsimple Seifert manifolds which are not cables on $I(K)$ as well). Assume then that $r_2$ is a finite, noncyclic filling slope. We shall show that $p = q = 3$ and then argue that $\Delta(r_1, r_2) = 1$.

It is proved in [12, Lemma 7.2] that if $r_0$ and $r_0'$ denote the slopes on $\partial M$ and $\partial M_1 = T$ corresponding to the cabling annulus in $C$, then $M(r_0) \cong M_1(r_0') \# L$ where $L$ is a nontrivial lens space. Note then that $r_1 \neq r_0$ and as $T$ compresses in $M(r_1)$, another application of [12, Lemma 7.2] implies that there is a slope $r_1'$ on $\partial M_1$ such that $\Delta(r_0', r_1') = n$ and $M(r_1) \cong M_1(r_1')$. Now if we fill $M_1$ along the fiber slope $\phi$ of its given Seifert fibration, we obtain the connected sum of two nontrivial lens spaces. Thus $r_1' \neq \phi$ and so $M_1(r_1')$ admits a Seifert structure whose base orbifold $\mathcal{B}$ is the 2-sphere with at least two cone points of order $p, q \geq 2$. Observe also that the core of the filling solid torus is a fiber of multiplicity $q_1 = \Delta(r_1', \phi) \geq 1$ in $M_1(r_1')$. As $r_1'$ is an infinite cyclic filling slope, it follows that $M_1(r_1') \cong S^1 \times S^2$ and that $\mathcal{B}$ has exactly two cone points, so $q_1 = 1$. An analysis of the Seifert structures on $S^1 \times S^2$ now forces us to conclude that $3 \leq p = q$.

The constraints on $p$ and $q$ we have derived imply that no filling of $M_1$ is simply connected. Thus from the previously noted decomposition $M(r_0) \cong M_1(r_0') \# L$ we can deduce that $r_2 \neq r_0$. With the argument as above, the compressibility of $T$ in $M(r_2)$ guarantees that there is a slope $r_2'$ on $\partial M_1$ such that $\Delta(r_0', r_2') = n$ and $M(r_2) \cong M_1(r_2')$. Now $r_2'$ cannot be the fiber slope $\phi$ on $\partial M_1$ since $M_1(r_2')$ has a finite fundamental group, and therefore $M(r_2')$ admits a Seifert structure whose base orbifold is that 2-sphere with two cone points of orders $p$ and a third whose order is $q_2 = \Delta(r_2', \phi) \geq 1$. But then the fact that $\pi_1(M_1(r_2'))$ is a finite noncyclic group requires that $p = 3$ and $q_2 = 2$.

Since $\Delta(r_0', r_1') = \Delta(r_0', r_2') = n$, we can orient representative curves for $r_0', r_1'$ and $r_2'$ so that the associated homology classes satisfy

(2.5) $$\alpha(r_2') = \alpha(r_1') + a\alpha(r_0')$$

for some $a \in \mathbf{Z}$. Then $n^2 \Delta(r_1, r_2) = \Delta(r_1', r_2') = |\alpha(r_1') \cdot \alpha(r_2')| = |a||\alpha(r_0') \cdot \alpha(r_1')| = |a|\Delta(r_0', r_1') = |a|n$, i.e. $|a| = n\Delta(r_1, r_2)$. Equation 2.5 also implies that $2 = q_2 = \Delta(r_2', \phi) = |\alpha(r_2') \cdot \alpha(\phi)| = |\alpha(r_1') \cdot \alpha(\phi) + a\alpha(r_0') \cdot \alpha(\phi)| = |\Delta(r_1', \phi) \pm a\Delta(r_0', \phi)| = |q_1 \pm a\Delta(r_0', \phi)| = |1 \pm a\Delta(r_0', \phi)|$. Hence $n\Delta(r_1, r_2)\Delta(r_0', \phi) = |a|\Delta(r_0', \phi)$ is either 1 or 3. But $n \geq 2$ and so we conclude, in particular, that $\Delta(r_1, r_2) = 1$. This completes the proof of Lemma 2.4. $\square$

LEMMA 2.5. *Let $M$ be a simple Seifert-fibered manifold whose boundary is an incompressible torus and suppose that $r_1$ is a non-Seifert filling slope on $\partial M$ while $r_2$ is a finite filling slope on $\partial M$. Then $\Delta(r_1, r_2) \leq 5$.*

*Proof.* $M$ admits a Seifert structure having base orbifold a 2-disc with exactly two cone points of orders $p$ and $q$ say. As $r_1$ is a non-Seifert filling slope, $M \neq I(K)$ (i.e. $(p, q) \neq (2, 2)$) and $r_1$ is the slope on $\partial M$ corresponding to the fiber. Hence $r_2$, being a finite filling slope, cannot coincide with $r_1$, and so $M(r_2)$ admits a Seifert structure having base orbifold the 2-sphere with exactly three cone points of orders $p, q$ and $\Delta(r_1, r_2)$. The group $\pi_1(M(r_2))$ surjects onto the $(p, q, \Delta(r_1, r_2))$-triangle group, which is infinite as long as $\Delta(r_1, r_2) > 5$. Thus the lemma follows. $\square$

We now briefly review some of the facts about incompressible surfaces in Seifert-fibered spaces needed in this paper.

Let $W$ be a closed Seifert-fibered space which contains an incompressible (orientable) surface $S$. Up to isotopy, $S$ is either horizontal (i.e. is transverse to the fibers of $W$) or vertical (i.e. consists of a disjoint union of fibers) [24, VI.34]. When $S$ is horizontal and does not separate $W$, it is a fiber in a realization of $W$ as the total space of a locally trivial surface bundle over the circle. Further, the algebraic intersection of a fiber of the Seifert fibration with a 2-cycle carried by $S$ is nonzero. When $S$ is horizontal and does separate $W$, then $W = M_1 \cup_S M_2$ where $M_1 \cap M_2 = \partial M_1 = \partial M_2 = S$ and each $M_i$ is a twisted $I$-bundle over a closed, nonorientable surface. Hence we have:

LEMMA 2.6. *If $S$ is a horizontal incompressible surface in a closed Seifert-fibered space $W$, then $\pi_1(S)$ is a normal subgroup of $\pi_1(W)$ and $\pi_1(W)/\pi_1(S) \cong \mathbf{Z}$ or $\mathbf{Z}/2 * \mathbf{Z}/2$.*

LEMMA 2.7. *Let $S_1, S_2$ be two closed, essential surfaces in a closed, irreducible Seifert-fibered manifold $W$ such that $S_1$ is isotopic to a horizontal surface.*

(1) *If $S_2$ is isotopic to a vertical surface, then $S_1 \cap S_2 \neq \emptyset$.*

(2) *If $S_1 \cap S_2 = \emptyset$, then $S_2$ is isotopic to $S_1$.*

*Proof.* As we noted above, either $S_1$ splits $W$ into two pieces, each a twisted $I$-bundle over an unoriented surface, or $S_1$ is a fiber of a realization of $W$ as the total space of a locally trivial surface bundle over the circle. To prove part (1), note that by possibly passing to a 2-fold cover we can assume that the latter occurs. Now as $S_2$ is isotopic to a vertical surface, it contains a closed 1-cycle whose algebraic intersection with the 2-cycle carried by $S_1$ is nonzero. This intersection number is a homology invariant, so in particular we must have $S_1 \cap S_2 \neq \emptyset$.

To prove part (2), assume that $S_1 \cap S_2 = \emptyset$. If $W$ is the total space of a locally trivial $S_1$-bundle over the circle, then our hypothesis implies that $S_2 \subset S_1 \times I \subset W$. Hence [45, Prop. 3.1] shows that $S_1$ and $S_2$ are isotopic. On the other hand suppose that $W = M_1 \cup_{S_1} M_2$ where $M_1, M_2$ are twisted $I$-bundles over an unorientable surface. From our hypothesis we may assume that $S_2 \subset M_1$ without loss of generality. Now $M_1$ is doubly covered by $\tilde{M}_1 \cong S_1 \times I$. Let $\tilde{S}_2$ be a component of the inverse image of $S_2$ in $\tilde{M}_1$. Another application of [45, Prop. 3.1] implies that $\tilde{S}_2$ is isotopic to $S_1 \times \{0\} \subset S_1 \times I$. In particular, composing this isotopy with the covering projection we obtain a homotopy of $\tilde{S}_2 \times I \to M_1$ which when restricted to $\tilde{S}_2 \times \{0\}$ is a cover to $S_2$ and when restricted to $\tilde{S}_2 \times \{1\}$ has image in $S_1$. Then [45, Lemma 5.3] now shows that $S_2$ is isotopic to $S_1$. □

## 3. Varieties of $\mathrm{PSL}_2(\mathbf{C})$-characters and central $\mathbf{Z}/2$-extensions

**of groups**

In this section we describe algebraic sets of $\text{PSL}_2(\mathbf{C})$ representations and characters, and relate them to their $\text{SL}_2(\mathbf{C})$ counterparts. The algebraic sets we consider will all be complex and affine. We shall reserve the term *variety* to denote an irreducible algebraic set.

Fix a finitely generated group $\Gamma$ and let $R(\Gamma)$ denote the algebraic set consisting of the representations of $\Gamma$ in $\text{SL}_2(\mathbf{C})$. There is an associated algebraic set $X(\Gamma)$ which consists of the characters of the representations in $R(\Gamma)$ [7, §1]. When $\Gamma$ is the fundamental group of a path-connected space $B$ based at a point $b_0$, we shall suppress the base point from our notation and use $R(B)$ and $X(B)$ to denote $R(\pi_1(B; b_0))$ and $X(\pi_1(B; b_0))$ respectively. The natural surjective map $R(\Gamma) \xrightarrow{t} X(\Gamma)$ which sends a representation $\rho$ to its character $\chi_\rho$ is a regular function.

There is a similar theory over $\text{PSL}_2(\mathbf{C})$. Consider $\bar{R}(\Gamma)$, the set of representations of $\Gamma$ in $\text{PSL}_2(\mathbf{C})$. The adjoint representation $\text{Ad} : \text{PSL}_2(\mathbf{C}) \to \text{Aut}(\text{sl}_2(\mathbf{C}))$ is injective and has image a closed, algebraic subgroup of $\text{SL}_3(\mathbf{C})$. Thus $\bar{R}(\Gamma)$ may be identified with a complex algebraic set. Now composition with the natural quotient map $\Phi : \text{SL}_2(\mathbf{C}) \to \text{PSL}_2(\mathbf{C})$ determines a regular map $\Phi_* : R(\Gamma) \to \bar{R}(\Gamma)$, and a simple calculation shows that each of its fibers either is empty or is an orbit of the free $H^1(\Gamma; \mathbf{Z}/2) = \text{Hom}(\Gamma, \{\pm I\})$-action on $R(\Gamma)$ defined by multiplication: $(\varepsilon \cdot \rho)(\gamma) = \varepsilon(\gamma)\rho(\gamma)$. The condition determining whether or not the fiber above $\bar{\rho} \in \bar{R}(\Gamma)$ is nonempty, i.e., whether or not $\bar{\rho}$ lifts to a representation $\rho \in R(\Gamma)$, is well-understood. The homomorphism $\bar{\rho}$ induces a continuous function $B(\bar{\rho}) : B(\Gamma) \to B(\text{PSL}_2(\mathbf{C}))$ between classifying spaces, and hence a principal $\text{PSL}_2(\mathbf{C})$-bundle $\xi_{\bar{\rho}}$ over $B(\Gamma)$. Then $\bar{\rho}$ lifts to $\text{SL}_2(\mathbf{C})$ if and only if the Stiefel-Whitney class $w_2(\xi_{\bar{\rho}})$ vanishes. Now the isomorphism class of $\xi_{\bar{\rho}}$ depends only on the topological component of $\bar{\rho}$ in $\bar{R}(\Gamma)$ (see [11, §4.5] for instance), and so the image of $R(\Gamma) \to \bar{R}(\Gamma)$ is a union of topological components of $\bar{R}(\Gamma)$ ([5, Th. 4.1]). Indeed, $R(\Gamma) \to \Phi_*(R(\Gamma))$ is a regular covering space with group $H^1(\Gamma; \mathbf{Z}/2)$.

The natural action of $\text{PSL}_2(\mathbf{C})$ on $\bar{R}(\Gamma)$ has an algebro-geometric quotient $\bar{X}(\Gamma)$ [35, Th. 3.3.5] and there is a surjective "quotient" map $\bar{t} : \bar{R}(\Gamma) \longrightarrow \bar{X}(\Gamma)$ which is constant on conjugacy classes of representations. The affine algebraic set $\bar{X}(\Gamma)$ is determined by the condition that its coordinate ring be isomorphic to the ring of invariants of the natural action of $\text{PSL}_2(\mathbf{C})$ on $\mathbf{C}[\bar{R}(\Gamma)]$. Thus for each $\gamma \in \Gamma$, the function $\bar{X}(\Gamma) \to \mathbf{C}$ given by $\bar{t}(\bar{\rho}) \mapsto (\text{trace}(\bar{\rho}(\gamma)))^2$ is regular.

*Definition.* A representation $\bar{\rho} \in \bar{R}(\Gamma)$ is called *irreducible* if it is not conjugate to a representation whose image lies in

$$\lambda\{\pm \begin{pmatrix} a & b \\ 0 & a^{-1} \end{pmatrix} \mid a, b \in \mathbf{C}, \quad a \neq 0 \rho\}.$$

Otherwise it is called *reducible*.

Our notion of irreducibility coincides with the notion of stability found, for instance, in [26, p. 53]. It follows that if $\bar{\rho}$ is irreducible, then $\bar{t}^{-1}(\bar{t}(\bar{\rho}))$ is the orbit of $\bar{\rho}$ under conjugation [35, Cor. 3.5.2]. When the analogous construction is performed on the $\mathrm{SL}_2(\mathbf{C})$-action on $R(\Gamma)$, the resulting quotient may be canonically identified with $X(\Gamma)$, and so by analogy, $\bar{X}(\Gamma)$ is called the set of $\mathrm{PSL}_2(\mathbf{C})$-*characters* of $\Gamma$ and $\bar{t}(\bar{\rho})$ is denoted by $\chi_{\bar{\rho}}$.

The function $\Phi_* : R(\Gamma) \to \bar{R}(\Gamma)$ induces a regular mapping $\Phi_\# : X(\Gamma) \to \bar{X}(\Gamma)$ and as above, the fibers of this map are either empty or the orbits of the $H^1(\Gamma; \mathbf{Z}/2) = \mathrm{Hom}(\Gamma, \{\pm I\})$-action on $X(\Gamma)$ defined by $\varepsilon \cdot \chi_\rho = \chi_{\varepsilon \cdot \rho}$. In general this action is not free, though the isotropy groups are readily determined ([31, pp. 513, 514]). The image of $\Phi_\#$ in $\bar{X}(\Gamma)$ is a union of topological components. Note that as $X(\Gamma) \to \bar{X}(\Gamma)$ is finite-to-one, $\dim(X(\Gamma)) \leq \dim(\bar{X}(\Gamma))$, and the inequality can be strict for some groups $\Gamma$ (see Example 3.2 below).

Suppose now that $\Gamma$ is generated by $x_1, x_2, \ldots, x_n$. The character of an $\mathrm{SL}_2(\mathbf{C})$-representation $\rho$ of $\Gamma$ is determined by its values on the $m = \frac{n(n^2+5)}{6}$ (not necessarily distinct) elements $y_1, y_2, \ldots, y_m$ of $\Gamma$ contained in the set $\{x_i \mid 1 \leq i \leq n\} \cup \{x_i x_j \mid 1 \leq i < j \leq n\} \cup \{x_i x_j x_k \mid 1 \leq i < j < k \leq n\}$ [44]. A similar result holds for $\mathrm{PSL}_2(\mathbf{C})$-characters. To describe it, let $F_n$ denote the free group on the symbols $\xi_1, \xi_2, \ldots, \xi_n$.

LEMMA 3.1. *Suppose that $\Gamma$ is generated by $x_1, x_2, \ldots, x_n$ and that $\bar{\rho}, \bar{\rho}' \in \bar{R}(\Gamma)$. Choose matrices $A_1, A_2, \ldots, A_n, B_1, B_2, \ldots, B_n \in \mathrm{SL}_2(\mathbf{C})$ satisfying $\bar{\rho}(x_i) = \pm A_i$ and $\bar{\rho}'(x_i) = \pm B_i$ for each $i$. Define $\rho, \rho' \in R(F_n)$ by requiring that $\rho(\xi_i) = A_i$ and $\rho'(\xi_i) = B_i$ for each $i \in \{1, 2, \ldots, n\}$. Let $y_1, y_2, \ldots, y_m$ be the $m = \frac{n(n^2+5)}{6}$ elements of $F_n$ associated to the generators $\xi_1, \ldots, \xi_n$, as described above. Then $\chi_{\bar{\rho}} = \chi_{\bar{\rho}'}$ if and only if there is a homomorphism $\epsilon \in \mathrm{Hom}(F_n, \{\pm 1\})$ for which $\mathrm{trace}(\rho'(y_j)) = \epsilon(y_j) \mathrm{trace}(\rho(y_j))$ for each $j \in \{1, 2, \ldots, m\}$.*

*Proof.* There is a natural surjection $F_n \to \Gamma$ which sends $\xi_i$ to $x_i$ and which induces inclusions $\bar{R}(\Gamma) \to \bar{R}(F_n)$ and $\bar{X}(\Gamma) \to \bar{X}(F_n)$. Thus it suffices to prove the lemma in the case where $\Gamma = F_n$.

Consider the commutative diagram of surjective maps

$$\begin{array}{ccc} R(F_n) & \xrightarrow{\Phi_*} & \bar{R}(F_n) \\ t \downarrow & & \downarrow \bar{t} \\ X(F_n) & \xrightarrow{\Phi_\#} & \bar{X}(F_n). \end{array}$$

From the diagram we see that $\chi_{\bar{\rho}} = \chi_{\bar{\rho}'}$ if and only if $\Phi_{\#}(\chi_\rho) = \Phi_{\#}(\chi_{\rho'})$. Now the latter identity holds if and only if there is a homomorphism $\epsilon \in \text{Hom}(F_n, \{\pm 1\})$ such that $\chi_{\rho'} = \epsilon \cdot \chi_\rho = \chi_{\epsilon\rho}$. But as we observed above, $\chi_{\rho'} = \chi_{\epsilon \cdot \rho}$ if and only if $\chi_{\rho'}(y_j) = \chi_{\epsilon \cdot \rho}(y_j)$ for each $j \in \{1, 2, \ldots, m\}$, that is, if and only if the conclusion of the lemma holds. □

*Example* 3.2. Suppose that $\Gamma = \mathbf{Z}/p * \mathbf{Z}/q$ where $1 < p, q < \infty$ and let $x$ generate $\mathbf{Z}/p$ and $y$ generate $\mathbf{Z}/q$. We shall show that $\bar{X}(\mathbf{Z}/p * \mathbf{Z}/q)$ consists of $([\frac{p}{2}]+1)([\frac{q}{2}]+1)$ algebraic components, all mutually disjoint, of which $[\frac{p}{2}][\frac{q}{2}]$ are isomorphic to complex lines, the rest being isolated points.

There is a homeomorphism

$$\bar{R}(\mathbf{Z}/p * \mathbf{Z}/q) \longrightarrow \{(a, b) \in \text{PSL}_2(\mathbf{C})^2 | a^p = b^q = \pm I\}$$

which associates the pair $(\bar{\rho}(x), \bar{\rho}(y))$ to a representation $\bar{\rho} \in \bar{R}(\mathbf{Z}/p * \mathbf{Z}/q)$. It is then a simple exercise to verify that both $\bar{R}(\mathbf{Z}/p * \mathbf{Z}/q)$ and $\bar{X}(\mathbf{Z}/p * \mathbf{Z}/q)$ have $([\frac{p}{2}]+1)([\frac{q}{2}]+1)$ topological components, the component of $\bar{\rho} \in \bar{R}(\mathbf{Z}/p * \mathbf{Z}/q)$ or its character $\chi_{\bar{\rho}} \in \bar{X}(\mathbf{Z}/p * \mathbf{Z}/q)$ being determined by the pair $(|\chi_\rho(x)|, |\chi_\rho(y)|) = (2\cos(\pi j/p), 2\cos(\pi k/q))$ for some $j, k$ satisfying $0 \leq j \leq [\frac{p}{2}]$ and $0 \leq k \leq [\frac{q}{2}]$.

Each topological component of $\bar{X}(\mathbf{Z}/p * \mathbf{Z}/q)$ is also an algebraic component, and is either a curve or an isolated point. More precisely, let $C(j, k)$ be a component of $\bar{X}(\mathbf{Z}/p * \mathbf{Z}/q)$ corresponding to the pair $(j, k)$. If $j = 0$ or $k = 0$, $C(j, k)$ is a point. Suppose then that $1 \leq j \leq [\frac{p}{2}]$ and $1 \leq k \leq [\frac{q}{2}]$. Set

$$\lambda = e^{\pi i j/p}, \quad \mu = e^{\pi i k/q}, \quad \tau = \mu + \mu^{-1},$$

and for each $a \in \mathbf{C}$ define $\bar{\rho}_a \in \bar{R}(\mathbf{Z}/p * \mathbf{Z}/q)$ by

$$\bar{\rho}_a(x) = \pm \begin{pmatrix} \lambda & 0 \\ 0 & \lambda^{-1} \end{pmatrix}, \quad \bar{\rho}_a(y) = \pm \begin{pmatrix} a & 1 \\ a(\tau - a) - 1 & \tau - a \end{pmatrix}.$$

If $\chi_{\bar{\rho}} \in C(j, k)$ is the character of an irreducible representation, then $\bar{\rho}$ is conjugate to $\bar{\rho}_a$ for some $a \in \mathbf{C}$. This can be seen by first conjugating $\bar{\rho}$ to make $\bar{\rho}(x)$ diagonal, and then conjugating $\bar{\rho}$ by diagonals to make the upper left hand entry of $\bar{\rho}(y)$ equal to 1. On the other hand, if $\chi_{\bar{\rho}} \in C(j, k)$ is the character of an reducible representation, one can show that $\chi_{\bar{\rho}} = \chi_{\bar{\rho}_a}$ where $a = \mu$ or $a = \mu^{-1}$. Thus there is a surjective regular map

$$\Psi : \mathbf{C} \to C(j, k), \ a \mapsto \chi_{\bar{\rho}_a}.$$

In fact $\Psi$ is a bijection unless $j = p/2$ or $k = q/2$, in which case $\Psi(a) = \Psi(b)$ if and only if $b = a$ or $b = \tau - a$. To prove this we appeal to Lemma 3.1 to see that $\Psi(a) = \Psi(b)$ if and only if there are $\epsilon_x, \epsilon_y \in \{\pm 1\}$ such that the following equations hold.

$$(\lambda + \lambda^{-1}) = \epsilon_x(\lambda + \lambda^{-1}),$$

$$\tau = \epsilon_y \tau,$$
$$(\lambda - \lambda^{-1})b + \lambda^{-1}\tau = \epsilon_x \epsilon_y((\lambda - \lambda^{-1})a + \lambda^{-1}\tau).$$

Now $j = p/2$ if and only if $(\lambda + \lambda^{-1}) = 0$, while $k = q/2$ if and only if $\tau = 0$. Hence if $j \neq p/2$ and $k \neq q/2$, then the first two of the equations above show that $\epsilon_x = \epsilon_y = 1$. If we plug these values into the third equation we deduce that $a = b$. Hence $\Psi$ is a bijection. On the other hand if $j = p/2$ or $k = q/2$, these equations readily imply that $\Psi(a) = \Psi(b)$ if and only if $b = a$ or $b = \tau - a$.

Next we show that when $1 \leq j \leq [\frac{p}{2}]$ and $1 \leq k \leq [\frac{q}{2}]$, the curve $C(j,k)$ is in fact isomorphic to a complex line. The function $g_1 : C(j,k) \to \mathbf{C}$ given by $g_1(\chi_{\bar\rho_a}) = (\mathrm{trace}(\bar\rho_a(xy))^2$ is regular. Since the derivative of $g_1 \circ \Psi$ vanishes only at $a = -\lambda^{-1}\tau/(\lambda - \lambda^{-1})$, $C(j,k)$ is smooth away from $\Psi(-\lambda^{-1}\tau/(\lambda - \lambda^{-1}))$. On the other hand, a similar argument using $g_2 : C(j,k) \to \mathbf{C}$, $g_2(\chi_{\bar\rho_a}) = (\mathrm{trace}(\bar\rho_a(x^{-1}y))^2$, shows that $C(j,k)$ is smooth away from $\Psi(\lambda\tau/(\lambda - \lambda^{-1}))$. If $\tau \neq 0$ and $\lambda \neq i$ (i.e. $j \neq p/2$ and $k \neq q/2$) then $-\lambda^{-1}\tau/(\lambda - \lambda^{-1}) \neq \lambda\tau/(\lambda - \lambda^{-1})$, and therefore $C(j,k)$ is smooth and $\Psi$ is an isomorphism. Again, on the other hand, if $j = p/2$ or $k = q/2$ it is easy to verify that $g_1$ is in fact a regular, birational isomorphism, and therefore an isomorphism by [41, p. 105, Cor. 2].

We have shown then that $\dim \bar{X}(\mathbf{Z}/p * \mathbf{Z}/q) = 1$ and there are $[\frac{p}{2}][\frac{q}{2}] \geq 1$ components of $\bar{X}(\mathbf{Z}/p * \mathbf{Z}/q)$ which are isomorphic to complex lines. Note that the curves are precisely the components of $\bar{X}(\mathbf{Z}/p * \mathbf{Z}/q)$ which contain the character of an irreducible representation. Further, the characters on $C(j,k)$ corresponding to reducible representations are those given by the values $a = \mu, \mu^{-1}$ of the parameter. There are exactly two such characters if both $j \neq p/2$ and $k \neq q/2$, and one otherwise.

Finally we note that a similar analysis shows that $X(\mathbf{Z}/p * \mathbf{Z}/q)$ consists of $([\frac{p}{2}] + 1)([\frac{q}{2}] + 1)$ mutually disjoint algebraic components precisely $([\frac{p+1}{2}] - 1)([\frac{q+1}{2}] - 1)$ of which are isomorphic to complex lines, the rest being isolated points. Note then that if $p = 2$, $0 = \dim X(\mathbf{Z}/p * \mathbf{Z}/q) < 1 = \dim \bar{X}(\mathbf{Z}/p * \mathbf{Z}/q)$. □

Consider once again an arbitrary finitely generated group $\Gamma$. Given a representation $\bar\rho : \Gamma \to \mathrm{PSL}_2(\mathbf{C})$, there is a central extension $\phi : \hat\Gamma \to \Gamma$ of $\Gamma$ by $\mathbf{Z}/2$ and a representation $\hat\rho : \hat\Gamma \to \mathrm{SL}_2(\mathbf{C})$ for which the following commutative diagram has exact lines.

$$\begin{array}{ccccccccc}
1 & \to & \mathbf{Z}/2 & \to & \hat\Gamma & \xrightarrow{\phi} & \Gamma & \to & 1 \\
 & & \| & & \hat\rho \downarrow & & \bar\rho \downarrow & & \\
1 & \to & \mathbf{Z}/2 & \to & \mathrm{SL}_2(\mathbf{C}) & \xrightarrow{\Phi} & \mathrm{PSL}_2(\mathbf{C}) & \to & 1.
\end{array}$$

The triple $(\hat\Gamma, \phi, \hat\rho)$ is (essentially) uniquely determined by these conditions and we shall refer to it as the central $\mathbf{Z}/2$-*extension of* $\Gamma$ *lifting* $\bar\rho$.

The central extensions of $\Gamma$ by $\mathbf{Z}/2$ are classified by $H^2(\Gamma; \mathbf{Z}/2)$ ([23, Th. 10.3]) and so in particular, $\bar{\rho}$ determines an element $w_2(\bar{\rho}) \in H^2(\Gamma; \mathbf{Z}/2)$ which is zero if and only if $\hat{\Gamma} \cong \Gamma \times \mathbf{Z}/2$, i.e. if and only if $\bar{\rho}$ lifts to an element of $R(\Gamma)$. In fact, the class $w_2(\bar{\rho})$ is precisely the element $w_2(\xi_{\bar{\rho}})$ discussed earlier in this section. Note then that if $H^2(\Gamma; \mathbf{Z}/2) \cong 0$, $\hat{\Gamma} \cong \Gamma \times \mathbf{Z}/2$ for every $\bar{\rho}$, and so each $\bar{\rho} \in \bar{R}(\Gamma)$ lifts to a representation $\rho \in R(\Gamma)$. As an example, consider a compact, connected, irreducible, orientable 3-manifold $M$ whose boundary is a torus and which satisfies $H_1(M; \mathbf{Z}/2) = \mathbf{Z}/2$. Then $M$ is aspherical and Lefschetz duality shows that $H^2(\pi_1(M); \mathbf{Z}/2) \cong 0$.

The next lemma will be used in our construction of actions of $\Gamma$ on trees associated to ideal points of curves of $\mathrm{PSL}_2(\mathbf{C})$-characters of $\Gamma$.

LEMMA 3.3. *Fix $\bar{\rho} \in \bar{R}(\Gamma)$ and suppose that $(\hat{\Gamma}, \phi, \hat{\rho})$ is the central $\mathbf{Z}/2$-extension of $\Gamma$ lifting $\bar{\rho}$. There is a bijective correspondence between the varieties $R_0 \subset \bar{R}(\Gamma)$ in which $\bar{\rho}$ is a smooth point, and the varieties $S_0 \subset R(\hat{\Gamma})$ in which $\hat{\rho}$ is a smooth point. In fact, if $R_0 \subset \bar{R}(\Gamma)$ corresponds to $S_0 \subset R(\hat{\Gamma})$, then the homomorphisms $\Phi : \mathrm{SL}_2(\mathbf{C}) \to \mathrm{PSL}_2(\mathbf{C})$ and $\phi : \hat{\Gamma} \to \Gamma$ induce a regular function $\phi_* : S_0 \to R_0$ which topologically is a finite, regular cover.*

*Proof.* Let $C(\bar{\rho})$ be the connected component of $\bar{\rho}$ in $\bar{R}(\Gamma)$. Now composition with $\phi : \hat{\Gamma} \to \Gamma$ determines an inclusion $\bar{R}(\Gamma) \subset \bar{R}(\hat{\Gamma})$ which is regular and so in particular, $C(\bar{\rho})$ includes into the topological component $C(\bar{\rho} \circ \phi)$ of $\bar{\rho} \circ \phi$ in $\bar{R}(\hat{\Gamma})$. Denote by $\varepsilon$ the nontrivial element of the central $\mathbf{Z}/2 \subset \hat{\Gamma}$. The evaluation map $C(\bar{\rho} \circ \phi) \to \mathrm{PSL}_2(\mathbf{C})$ sending $\rho'$ to $\rho'(\varepsilon)$ is continuous and has image contained in the set of elements of $\mathrm{PSL}_2(\mathbf{C})$ whose order is either one or two. This set has precisely two components, one consisting of the elements of order two in $\mathrm{PSL}_2(\mathbf{C})$ and the other consisting of only the identity element. As $\bar{\rho} \circ \phi(\varepsilon) = 1$, $\varepsilon$ must be sent to the identity by each representation in $C(\bar{\rho} \circ \phi)$, from which it follows that the image of $C(\bar{\rho})$ in $\bar{R}(\hat{\Gamma})$ is exactly $C(\bar{\rho} \circ \phi)$, i.e. $C(\bar{\rho}) \cong C(\bar{\rho} \circ \phi)$.

Next recall the natural function $\Phi_* : R(\hat{\Gamma}) \to \Phi_*(R(\hat{\Gamma}))$. We have already discussed the fact that either $\Phi_*^{-1}(C(\bar{\rho} \circ \phi)) \subset R(\hat{\Gamma})$ is empty or $\Phi_*| : \Phi_*^{-1}(C(\bar{\rho} \circ \phi)) \to C(\bar{\rho} \circ \phi)$ is a regular cover with group $H^1(\hat{\Gamma}; \mathbf{Z}/2)$. The identity $\bar{\rho} \circ \phi = \Phi \circ \hat{\rho}$ shows that $\bar{\rho} \circ \phi$ lifts to $\hat{\rho} \in R(\hat{\Gamma})$, and so the latter case must arise. It follows then that the natural map $C(\hat{\rho}) \to C(\bar{\rho} \circ \phi) \equiv C(\bar{\rho})$ is a finite-to-one regular cover. The commutative diagram in Figure 3.1 below summarizes the situation.

Now consider a variety $R_0 \subset \bar{R}(\Gamma)$ which contains $\bar{\rho}$ as a smooth point. The closed set $R_0$ is connected ([34, Cor. 4.16]) and so is a subset of $C(\bar{\rho})$. Therefore its image in $C(\bar{\rho} \circ \phi) \cong C(\bar{\rho})$, $\hat{R}_0$ say, is closed in $\bar{R}(\hat{\Gamma})$ and so by [41, Th. 6, p. 50], is a variety. From the previous paragraph we know that $\Phi_*^{-1}(\hat{R}_0)$ is the total space of a regular $H^1(\hat{\Gamma}; \mathbf{Z}/2)$-cover of $\hat{R}_0$ and thus our hypotheses imply that $\hat{\rho}$ is a simple point in $\Phi_*^{-1}(\hat{R}_0)$. But then $\hat{\rho}$ is contained in a unique algebraic component $S_0$ of $\Phi_*^{-1}(\hat{R}_0)$ and is a smooth point there. Thus $R_0$ determines an appropriate variety $S_0 \subset R(\hat{\Gamma})$.

The inverse operation is described as follows. Suppose that $S_0 \subset R(\hat{\Gamma})$ is a variety containing $\hat{\rho}$ as a smooth point. Since $S_0$ is a closed subset of $C(\hat{\rho})$ and $\Phi_*$ is a closed function, we may apply [41, Th. 6, p. 50] to see that $\hat{R}_0 = \Phi_*(S_0)$ is a variety in $\bar{R}(\hat{\Gamma})$ containing $\phi \circ \hat{\rho}$ as a smooth point. For every $\rho' \in \hat{R}_0$ we have $\rho'(\varepsilon) = I$ and so $\hat{R}_0$ lies in the image of $\bar{R}(\Gamma)$. It follows that there is a variety $R_0$ containing $\bar{\rho}$ as a smooth point. This operation is clearly inverse to the previous one, thus completing the proof of the lemma. □

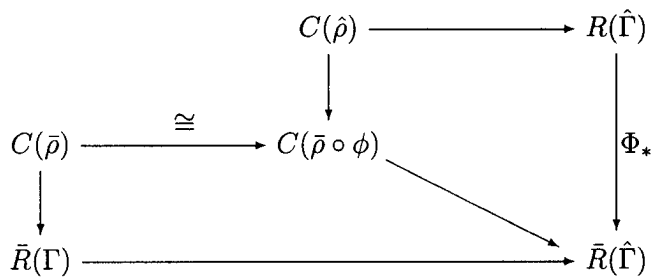

Figure 3.1.

## 4. $\mathrm{PSL}_2(\mathbf{C})$ Culler-Shalen theory

In this section we first review some of the ideas from the work of M. Culler and P. Shalen expounded in [7] and Chapter 1 of [6] with an eye to describing how they naturally extend to curves of $\mathrm{PSL}_2(\mathbf{C})$-characters. The reader is directed to [6] and [7] for more complete details.

One of the main results of [7] is that if $X(\Gamma)$ has positive dimension, then there is a nontrivial splitting of $\Gamma$ in terms of free products with amalgamation and HNN-extensions. An important implication of this result is that if $M$ is a compact 3-manifold with $X(M)$ having positive dimension, then there is a splitting of the manifold $M$ along essential surfaces [7, Prop. 2.3.1]. The construction of the splitting of $\Gamma$ begins by choice of a curve $X_0$ in $X(\Gamma)$ and a variety $R_0 \subset R(\Gamma)$ such that $t(R_0) = X_0$. The function field $F$ of $R_0$ may be considered a finite extension of the function field $K$ of $X_0$. An ideal point $x$ of the smooth, projective model of $X_0$ determines a discrete valuation on $K$ which may be extended to a discrete valuation on $F$. Following a construction of Tits and Serre, this valuation on $F$ determines an action of $\mathrm{SL}_2(F)$ on a simplicial tree $T$, an action which factors through $\mathrm{PSL}_2(F)$. We shall exploit this fact later on in this section.

Next the *tautological representation* $P : \Gamma \to \mathrm{SL}_2(F)$ is introduced. This is the homomorphism given by

$$P(\gamma) = \begin{pmatrix} a & b \\ c & d \end{pmatrix}$$

where the functions $a, b, c$ and $d$ are determined by the identity

$$\rho(\gamma) = \begin{pmatrix} a(\rho) & b(\rho) \\ c(\rho) & d(\rho) \end{pmatrix}$$

for all $\rho \in R_0$. The action of $\mathrm{SL}_2(F)$ on the tree $T$ constructed above pulls back via $P$ to an action of $\Gamma$ on $T$ which can be shown to be nontrivial. The nontrivial splitting of $\Gamma$ is then produced by applying Bass-Serre-Tits theory.

Consider now curves of $\mathrm{PSL}_2(\mathbf{C})$-characters of $\Gamma$. Morgan [30] has shown how to build actions of $\Gamma$ on trees associated to ideal points of a curve of $\mathrm{PSL}_2(\mathbf{C}) = \mathrm{SO}_{\mathbf{R}}^+(3,1)$-characters. We shall provide another such construction below which, though somewhat ad hoc, is convenient for our applications.

LEMMA 4.1. *For any curve $X_0 \subset \bar{X}(\Gamma)$, there is an algebraic component $R_0$ of $\bar{t}^{-1}(X_0)$ which satisfies $\bar{t}(R_0) = X_0$. If $X_0$ contains the character of an irreducible representation, then $R_0$ is 4-dimensional and is uniquely determined by $X_0$.*

*Proof.* Let $V_1, V_2, \ldots, V_n$ be the algebraic components of $\bar{t}^{-1}(X_0)$. The set $\bar{t}^{-1}(X_0)$ is closed under conjugation and so there is a regular map $g: V_i \times \mathrm{PSL}_2(\mathbf{C}) \to \bar{t}^{-1}(X_0)$ defined by $(\rho, g) \mapsto g\rho g^{-1}$. The image of this map contains $V_i$ and since $V_i \times \mathrm{PSL}_2(\mathbf{C})$ is a variety, $\overline{g(V_i \times \mathrm{PSL}_2(\mathbf{C}))} \subset \bar{t}^{-1}(X_0)$ is irreducible. By construction, this is only possible if the image of $g$ is $V_i$. Thus each $V_i$ is closed under conjugation.

According to [35, Th. 3.3.5(iv)], the image by $\bar{t}$ of a closed invariant subset of $\bar{R}(\Gamma)$ is a Zariski-closed subset of $\bar{X}(\Gamma)$. Thus for each $i$, $\bar{t}(V_i)$ is either a point or $X_0$. Now $\bar{t}$ is surjective and $X_0$ is an infinite set, so there is some index $i$ such that $\bar{t}(V_i) = X_0$. Taking $R_0$ to be such a $V_i$ proves the first part of the lemma.

Now assume that $X_0$ contains the character of an irreducible representation. After possibly reordering the components $V_1, V_2, \ldots, V_n$, there is some $k \in \{1, 2, \ldots, n\}$ such that $\bar{t}(V_i) = X_0$ when $1 \leq i \leq k$, while $\bar{t}(V_i)$ is a point when $k+1 \leq i \leq n$. Suppose that $k \geq 2$ so that in particular $\bar{t}(V_1) = \bar{t}(V_2) = X_0$. Now as $X_0$ contains the character of an irreducible representation, we may argue as in [7, Cor. 1.5.3] to see that $V_1$ and $V_2$ both have dimension 4. Further $X_0$ contains a Zariski open set $U$ consisting of characters of irreducible representations [7, Cor. 1.2.2]. As an irreducible representation is determined up to conjugation by its character, and both $V_1$ and $V_2$ are closed under conjugation, it follows that each one contains the 4-dimensional set $\bar{t}^{-1}(U)$. But then $V_1 = V_2$, which is clearly impossible. Hence $k = 1$ and the lemma follows. □

Fix a curve $X_0 \subset \bar{X}(\Gamma)$. Our first task is to construct a tautological representation of $\Gamma$ associated to $X_0$. This will be accomplished through the aid of Lemma 3.3. Let $R_0 \subset \bar{t}^{-1}(X_0)$ be a subvariety of $\bar{R}(\Gamma)$ guaranteed by Lemma 4.1. Choose a smooth point $\bar{\rho}_0 \in R_0$ and let $(\hat{\Gamma}, \phi, \hat{\rho}_0)$ be the central $\mathbf{Z}/2$-extension of $\Gamma$ lifting $\bar{\rho}_0$. According to Lemma 3.3, there are a variety $S_0 \subset R(\hat{\Gamma})$ containing $\hat{\rho}$ and a dominating regular map $\phi_* : S_0 \to R_0$. Denote the function fields $\mathbf{C}(X_0)$ and $\mathbf{C}(S_0)$ by $K$ and $F$ respectively and note that $F$ may be considered a finitely generated extension field of $K$ through the identification of $f \in K$ with $f \circ \bar{t} \circ \phi_* \in F$.

There is a tautological representation $\hat{P} : \hat{\Gamma} \to \mathrm{SL}_2(F)$ associated to $S_0$, previously discussed in this section, which, when composed with the projection $\mathrm{SL}_2(F) \to \mathrm{PSL}_2(F)$, factors through $\Gamma$ to give us the *tautological representation*

$$P : \Gamma \to \mathrm{PSL}_2(F).$$

Observe that

$$P(\gamma)(\hat{\rho}) = \rho(\gamma)$$

for each $\gamma \in \Gamma$ and for each $\rho \in R_0$ and $\hat{\rho} \in S_0$ such that $\phi_*(\hat{\rho}) = \rho$.

We now use $P$ to construct actions of $\Gamma$ on simplicial trees.

Let $\tilde{X}_0$ be the smooth, projective model of $X_0$. It is described in [6, §1.5] how if $X_0^\nu$ is the normalisation of $X_0$, then there is a birational equivalence $\nu : X_0^\nu \to X_0$ as well as a natural inclusion $X_0^\nu \subset \tilde{X}_0$ whose image consists of the *ordinary* points of $\tilde{X}_0$ and whose complement consists of the *ideal* points of $\tilde{X}_0$. Consider an ideal point $x$ of $\tilde{X}_0$. This point determines a discrete valuation $w : K^* \to \mathbf{Z}$ which, by [31, Lemma II.4.4], may be extended to a discrete valuation $v : F^* \to \mathbf{Z}$. Now the action of $\mathrm{SL}_2(F)$ on the associated tree $T$ ([40]) factors through $\mathrm{PSL}_2(F)$ and if $O_v$ denotes the valuation ring of $v$, then each vertex stabiliser of this action is conjugate in $PGL_2(F)$ to $\mathrm{PSL}_2(O_v)$. Observe that the action of $\mathrm{PSL}_2(F)$ on the tree $T$ constructed above pulls back via $P$ to an action of $\Gamma$ on $T$. We shall say that this is an action *associated* to $x$. With this construction, most of Chapter 1 of [6] goes through with only minor modifications, as we shall describe now.

Each element $\gamma \in \Gamma$ defines a regular function

$$f_\gamma : \bar{X}(\Gamma) \to \mathbf{C}$$

$$f_\gamma(\chi_{\bar{\rho}}) = \mathrm{trace}(\Phi^{-1}(\bar{\rho}(\gamma)))^2 - 4.$$

There is an associated holomorphic mapping $\widetilde{f_\gamma | X_0} : \tilde{X}_0 \to \mathbf{C}P^1$. As the curve $X_0$ will always be clear from the context, we shall abbreviate $\widetilde{f_\gamma | X_0}$ to

$$\tilde{f}_\gamma : \tilde{X}_0 \to \mathbf{C}P^1.$$

PROPOSITION 4.2 ([6, Prop. 1.2.6]).    *Let $\Gamma \times T \to T$ be an action of $\Gamma$ on a tree $T$ associated to an ideal point $x$ of a curve $X_0$ in $\bar{X}(\Gamma)$. This action* is nontrivial. Indeed, an element $\gamma \in \Gamma$ has a fixed point in $T$ if and only if $\tilde{f}_\gamma(x) \in \mathbf{C}$.

The following theorem is an easy consequence of Proposition 4.2.

THEOREM 4.3 ([7], [30]). *If $\bar{X}(\Gamma)$ has positive dimension, then $\Gamma$ has a nontrivial splitting. In particular if $\Gamma$ is the fundamental group of a compact 3-manifold $M$, then $M$ has a splitting along a nonempty essential surface.*

*Proof.* Fix a curve $\bar{X}_0$ in $\bar{X}(\Gamma)$ and an ideal point $x$ of $\tilde{X}_0$. Then Proposition 4.2 guarantees the existence of a nontrivial action of $\Gamma$ on a tree, and therefore a nontrivial splitting of $\Gamma$. □

We note that Theorem 4.3 is more generally applicable than the analogous theorem for curves of $\mathrm{SL}_2(\mathbf{C})$-characters [7, §2]. As an example, consider the manifold $M = L(2,1) \# L(p,q)$, $p \geq 2$. The splitting of its fundamental group arising from the connected sum is readily seen to be associated to an affine curve in $\bar{X}(M)$, but not one in $X(M)$ as the latter is zero-dimensional (see Example 3.2).

The next proposition will be used in the later sections to analyse, amongst other things, Dehn fillings of simple, non-Seifert manifolds which yield Seifert-fibered spaces. We continue to use the notation (e.g. $R_0, \phi, \hat{\Gamma}, S_0, F = \mathbf{C}(S_0)$) developed for our definition of a tautological representation $P : \Gamma \to \mathrm{PSL}_2(F)$ associated to a curve $X_0 \in \bar{X}(\Gamma)$.

PROPOSITION 4.4 ([6, Prop. 1.2.7]). *Let $\Gamma \times T \to T$ be an action of $\Gamma$ on a tree $T$ associated to an ideal point $x$ of a curve $X_0$ in $\bar{X}(\Gamma)$. If there is a normal subgroup $N$ of $\Gamma$ which fixes a point of $T$, then either $P(N) = \{\pm I\}$ or $P(\Gamma)$ contains a diagonalisable subgroup of index $2$.*

*Definition* 4.5. Let $q$ be a positive integer.

(1) We say that $X_0$ is *index $q$ virtually abelian* if it contains the character of an irreducible representation and there is an index $q$ subgroup $\Gamma_0$ of $\Gamma$ for which $\rho(\Gamma_0)$ is abelian for each $\rho \in R_0$.

(2) We say that $X_0$ is *index $q$ virtually reducible* if there is an index $q$ subgroup $\Gamma_0$ of $\Gamma$ for which $\rho|\Gamma_0$ is reducible for each $\rho \in \bar{R}(\Gamma)$ such that $\chi_\rho \in X_0$. If $X_0$ is not index $q$ virtually reducible, we say that it is *index $q$ virtually irreducible*.

Observe that if $X_0$ contains the character of an irreducible representation, then $R_0$ is uniquely determined by $X_0$ (Lemma 4.1). Thus there is no ambiguity in the second of the two definitions above.

LEMMA 4.6. *Let $X_0 \subset \bar{X}(\Gamma)$ be a curve.*

(1) *Suppose that $X_0$ contains the character of an irreducible representation. If there is an index $q$ subgroup $\Gamma_0$ of $\Gamma$ for which $P|\Gamma_0$ is diagonalisable, then $X_0$ is index $q$ virtually abelian.*

(2) *If $X_0$ is index $q$ virtually abelian, then it is index $q$ virtually reducible.*

(3) *If there is an index $q$ subgroup $\Gamma_0$ of $\Gamma$ for which $P|\Gamma_0$ is reducible, then $X_0$ is index $q$ virtually reducible.*

*Proof.* (1) Let $\Gamma_0 \subset \Gamma$ be a subgroup satisfying the hypotheses part (1) of the lemma. We must show that for each $\rho \in R_0$, $\rho(\Gamma_0)$ is abelian. But this is obvious because for each $\rho \in R_0$, there is a homomorphism from the abelian group $P(\Gamma_0)$ onto $\rho(\Gamma_0)$ which sends $P(\gamma) \in \text{PSL}_2(F)$ to $P(\gamma)(\hat{\rho}) = \rho(\gamma) \in \text{PSL}_2(\mathbf{C})$. Thus $\rho(\Gamma_0)$ is an abelian group.

(2) Suppose that $X_0$ contains the character of an irreducible representation and that $\Gamma_0$ is a subgroup of index $q$ of $\Gamma$ for which $\rho(\Gamma_0)$ is abelian for each $\rho \in R_0$. An abelian subgroup of $\text{PSL}_2(\mathbf{C})$ is either reducible (i.e. conjugate to a group of $\pm$-classes of upper triangular matrices), or is irreducible and isomorphic to $\mathbf{Z}/2 \oplus \mathbf{Z}/2$. As any two $\mathbf{Z}/2 \oplus \mathbf{Z}/2$ subgroups of $\text{PSL}_2(\mathbf{C})$ are conjugate, $\Gamma_0$ admits only finitely many $\text{PSL}_2(\mathbf{C})$-orbits of representations whose images are isomorphic to $\mathbf{Z}/2 \oplus \mathbf{Z}/2$. Furthermore, any such representation has a neighbourhood consisting entirely of irreducible representations. Hence our assumptions imply that if $\rho(\Gamma_0) \cong \mathbf{Z}/2 \oplus \mathbf{Z}/2$ for some $\rho \in R_0$, then there is a neighbourhood of $\rho$ in $R_0$ which lies in the orbit of $\rho$. But this is impossible as the orbit of $\rho$ is 3-dimensional while $R_0$ is 4-dimensional (Lemma 4.1). Thus $\rho|G_0$ is reducible for each $\rho \in R_0$, which is what we needed to prove.

(3) Let $\Gamma_0 \subset \Gamma$ be a subgroup satisfying the hypotheses of part (3) of the lemma. Now by hypothesis, $P(\Gamma_0)$ is conjugate into the image in $\text{PSL}_2(F)$ of the upper triangular matrices of $\text{SL}_2(F)$. Let $\hat{\Gamma}_0 = \phi^{-1}(\Gamma_0) \subset \hat{\Gamma}$. It follows from the definition of $P$ that $\hat{P}|\hat{\Gamma}_0$ is reducible.

Suppose that $\rho_1 \in \bar{R}(\Gamma)$ is a representation with $\chi_{\rho_1} \in X_0$, say $\chi_{\rho_1} = \chi_\rho$ where $\rho \in R_0$. Fix $\hat{\rho} \in S_0$ such that $\rho = \hat{\rho} \circ \phi$ and observe that for each $\hat{\gamma} \in [\hat{\Gamma}_0, \hat{\Gamma}_0]$, the commutator subgroup of $\hat{\Gamma}_0$, $\chi_{\hat{\rho}}(\hat{\gamma}) = \chi_{\hat{P}}(\hat{\gamma})(\hat{\rho}) = 2$. It follows then from [7, Lemma 1.2.1] that $\hat{\rho}$, and hence $\rho$, is reducible. Finally, the irreducibility of a representation in $\bar{R}(\Gamma)$ is determined by its character (see for example [7, Prop. 1.5.2]) and therefore $\rho_1$ is reducible. □

Consider a compact, connected, irreducible and $\partial$-irreducible 3-manifold $M$ with $\partial M$ a torus. Let $L$ denote the group $H_1(\partial M)$, which will be regarded as a lattice in the 2-dimensional real vector space $V = H_1(\partial M; \mathbf{R})$. We will let $e: L \to \pi_1(\partial M)$ denote the inverse of the Hurewicz isomorphism. This gives a canonical identification of $L$ with $\pi_1(\partial M)$. Further, $\pi_1(\partial M)$ will be identified with its image under the natural inclusion $\pi_1(\partial M) \to \pi_1(M)$. Note that even though this inclusion is only well-defined up to conjugation, if $\gamma \in \pi_1(\partial M)$, then the function $f_\gamma : \bar{X}(M) \to \mathbf{C}$ is nevertheless well-defined. Clearly $f_{\gamma^{-1}} = f_\gamma$. Hence a slope $r$ on $\partial M$ determines a well-defined function $f_{e(\alpha(r))}$, where $\alpha(r)$ is either one of the primitive elements of $H_1(\partial M)$ associated to $r$. For reasons of notational simplicity, if $\alpha \in H_1(\partial M)$ we shall use $f_\alpha$ to denote $f_{e(\alpha)}$, and if $r$ is a slope on $\partial M$, we shall let $f_r$ denote $f_{e(\alpha(r))}$.

Of particular interest to us is the behaviour of the zeros and poles of the functions $f_\gamma$, where $\gamma \in \pi_1(\partial M)$. Let $x$ be a point in the smooth, projective model $\tilde{X}_0$ of a curve $X_0$ in $\bar{X}(M)$. If $\tilde{f} \in \mathbf{C}(\tilde{X}_0)$, define

$$\Pi_x(\tilde{f}) \text{ to be the multiplicity of } x \text{ as a pole of } \tilde{f},$$

and define

$$Z_x(\tilde{f}) \text{ to be } \begin{cases} \text{the multiplicity of } x \text{ as a zero of } \tilde{f} & \text{if } f \not\equiv 0 \\ +\infty & \text{otherwise.} \end{cases}$$

PROPOSITION 4.7 ([6, Prop. 1.3.8, Prop. 1.3.9, and Lemma 1.4.1]). *Let $x$ be an ideal point of $\tilde{X}_0$, the smooth, projective model of a curve $X_0$ in $\bar{X}(M)$. Then there is a linear function $\phi_x : L \to \mathbf{Z}$ such that $\Pi_x(\tilde{f}_\alpha) = |\phi_x(\alpha)|$ for each $\alpha \in L$. Furthermore, if*

(1) *$\phi_x \equiv 0$, so that $\tilde{f}_\alpha(x) \in \mathbf{C}$ for each $\alpha \in L$, then there is a closed essential surface in $M$ associated to $x$;*

(2) *$\phi_x \not\equiv 0$, then there is a unique slope $r$ on $\partial M$ for which $\phi_x(\alpha(r)) = 0$, i.e. for which $\tilde{f}_r(x) \in \mathbf{C}$. Moreover, $r$ is necessarily a boundary slope, and if $X_0$ contains the character of an irreducible representation but is not index 2 virtually abelian, then $r$ is a strict boundary slope.*

PROPOSITION 4.8 ([6, Prop. 1.5.2]). *Let $x$ be an ordinary point in the smooth, projective model $\tilde{X}_0$ of a curve $X_0$ in $\bar{X}(M)$ which contains the character of an irreducible representation. Let $R_0$ be the subvariety of $\bar{R}(M)$ such that $\bar{t}(R_0) = X_0$. Suppose that $Z_x(\tilde{f}_\alpha) > Z_x(\tilde{f}_\beta)$ for some $\alpha, \beta \in L$. Then any representation in $R_0$ whose character is $\nu(x)$ sends $\alpha$ to $\pm I$. Further, there is a representation $\bar{\rho} \in \bar{t}^{-1}(X_0)$ such that $\bar{t}(\bar{\rho}) = \nu(x)$ and the image of $\bar{\rho}$ in $\mathrm{PSL}_2(\mathbf{C})$ is noncyclic.*

PROPOSITION 4.9 ([6, Prop. 1.6.1]). *Let $x$ be an ideal point in the smooth, projective model $\tilde{X}_0$ of a curve $X_0$ in $\bar{X}(M)$. Let $r$ be a slope on $\partial M$ and suppose that $Z_x(\tilde{f}_r) > Z_x(\tilde{f}_\beta)$ for some $\beta \in L$. If $r$ is not a boundary slope, then there is a closed essential surface in $M$ which remains essential in $M(r)$. If $X_0$ contains the character of an irreducible representation but is not index 2 virtually abelian, then we need only assume that $r$ is not a strict boundary slope.*

It is essential for our development and applications of the theory of $\mathrm{PSL}_2(\mathbf{C})$-character varieties that we understand the relationship between the behaviour of the functions $f_r$ at the ideal points of a curve $X_0$ and the topological properties of $M$. Proposition 4.9 provides such a relationship when $r$ is not a boundary slope. In Propositions 4.10 and 4.12 below we provide other such connections.

PROPOSITION 4.10. *Let $X_0$ be a curve of $\mathrm{PSL}_2(\mathbf{C})$-characters of $\pi_1(M)$ which contains the character of an irreducible representation. Let $R_0$ be the unique subvariety of $\bar{R}(M)$ for which $\bar{t}(R_0) = X_0$ and suppose that there is a slope $r$ on $\partial M$ such that $\rho(e(\alpha(r))) = \pm I$ for each $\rho \in R_0$. Suppose further that there is an ideal point $x$ of $\tilde{X}_0$ at which each $\tilde{f}_\alpha(x)$, $\alpha \in L$, is finite. Then $M$ contains a closed, essential surface $S$ such that if $S$ compresses in $M(r)$ and $M(r_1)$ for some slope $r_1$ on $\partial M$, then $\Delta(r, r_1) \leq 1$.*

*Proof.* The proof is modelled on that of [6, Prop. 1.6.1].

For a closed surface $S$, define $\chi_-(S) = \sum -\chi(S_i)$ where the sum is taken over all components $S_i$ of $S$ which are not 2-spheres, and define $\#(S)$ to be the number of components of $S$. The *complexity* of $S$ is the pair $(\chi_-(S), \#(S))$, which we order lexicographically.

According to Proposition 4.7 (1), there is a closed surface in $M$ associated to the action of $\pi_1(M)$ on a tree provided by the ideal point $x$. Amongst all such closed surfaces in $M$, let $S$ be one having minimal complexity. Then $S$ is essential in $M$ ([6, Lemma 1.6.4]), and therefore the irreducibility of $M$ implies that each component of $S$ has genus larger than zero.

Suppose that $r_1$ is a slope on $\partial M$ such that $\Delta(r, r_1) > 1$. We claim that there is a component of $S$ which remains essential in at least one of $M(r)$ or $M(r_1)$. Suppose to the contrary that every component of $S$ compresses in both of these manifolds. Denote by $W$ the component of $M \setminus S$ which contains $\partial M$ and observe that $W$ is an irreducible and $\partial$-irreducible manifold which is not homeomorphic to $\partial M \times [0,1]$. Furthermore, the Dehn fillings $W(\partial M; r)$ and $W(\partial M; r_1)$ of $W$ are $\partial$-reducible. As we have assumed that $\Delta(r, r_1) > 1$, it follows from Lemma 2.3 (2) that there is an essential annulus $A$ in $W$ with one component of $\partial A$ in a component of $\partial W$ different from $\partial M$ and with the other component of $\partial A$ in $\partial M$ having some slope $r_0$ such that $\Delta(r_0, r) = \Delta(r_0, r_1) = 1$. In particular we have $r_0 \neq r$.

Let $J$ denote the solid torus attached to $\partial M$ in forming $W(\partial M; r)$. By hypothesis there is a component $S_0$ of $S \cap \partial W$ which compresses in $W(\partial M; r)$. In particular there is a compressing disk $D$ in $W(\partial M; r)$ whose boundary lies in $S_0$ and which intersects $J$ in a finite nonempty collection of meridional disks. Let $Q_0$ be the surface obtained by first compressing $S_0$ in $W(\partial M; r)$ using the disk $D$, and then pushing the resulting surface into the interior of $W(\partial M; r)$ by a small isotopy which is fixed near $J$. Finally let $Q_0^- = Q_0 \cap W$ and let $Q^- = Q_0^- \cup (S \setminus S_0)$.

*Claim.* The surface $Q^- \subset M$ is associated to the action of $\pi_1(M)$ on the tree provided by the ideal point $x$.

*Proof of the Claim.* Set $\pi = \pi_1(M)$ and recall there are a central $\mathbf{Z}/2$-extension $\phi : \hat{\pi} \to \pi$, a subvariety $S_0 \subset R(\hat{\pi})$, and a regular morphism $\phi_* : S_0 \to R_0$ which are used to define the tautological representation $P : \pi \to \mathrm{PSL}_2(F)$, where $F$ is the function field of $S_0$. Let $\alpha(r) \in H_1(\partial M)$ be one of the primitive homology classes associated to the slope $r$. From the definition of $P$ and our hypotheses we have that for any $\hat{\rho} \in S_0$, $P(e(\alpha(r)))(\hat{\rho}) = \phi_*(\hat{\rho})(e(\alpha(r))) = \pm I$. Thus $P(e(\alpha(r))) = \pm I$. It follows that $e(\alpha(r))$ acts trivially on any tree associated to the ideal point $x$. Now the proof of the claim follows exactly as in the proof of [6, Claim 1.6.8], where we replace [6, Lemma 1.6.5] by our stronger condition that $e(\alpha(r))$ acts trivially on the whole tree. □

It is possible that $Q_0^-$ is not an essential surface in $W$, but after a finite sequence of compressions (which take place in $W$), as well as deletions of 2-spheres and $\partial$-parallel components, we may replace $Q_0^-$ by a (possibly empty) essential surface $Q_0^* \subset W$. According to [6, Prop. 1.3.6], the new surface $Q_0^* \cup (S \setminus S_0)$ is still associated to the action of $\pi_1(M)$ on the tree.

Observe first of all that no component of $Q_0^*$ can have nonempty boundary. For by construction, $W$ contains an annulus $A$ which intersects $\partial M$ in a slope $r_0 \neq r$. If $R$ is a component of $Q_0^*$ with $\partial R \neq \emptyset$, the slope on $\partial M$ associated to $\partial R$ is clearly $r$. But this contradicts the following lemma.

LEMMA 4.11.  *Let $M$ be an irreducible, $\partial$-irreducible 3-manifold whose boundary contains a torus $T$, and suppose that there is an essential annulus $A$ in $M$ having exactly one boundary component in $T$, and that this component has slope $r_0$. Then every essential surface in $M$ whose boundary is contained in $T$ has boundary slope $r_0$.*

*Proof.* Let $F$ be an essential surface in $M$ with boundary contained in $T$. Since $F$ is essential it can be isotoped so that it meets $A$ transversely and so that no component of $A \cap F$ is an inessential arc in $A$. But since $\partial F \subset T$ and since $A$ has exactly one boundary component in $T$, any arc component of $A \cap F$ would have both its endpoints in the same boundary component of $A$ and would therefore be inessential. It follows that no component of $A \cap F$ is an arc, and hence we must have $\partial F \cap \partial A = \emptyset$. It follows that $F$ has boundary slope $r_0$. □

We conclude then that $Q_0^*$ is closed. Now the complexity of the surface $Q_0^* \cup (S \setminus S_0)$ is less than or equal to that of $Q_0 \cup (S \setminus S_0)$ ([6, p. 264]), while the complexity of $Q_0 \cup (S \setminus S_0)$ is less than or equal to that of $S$ with a strict inequality as long as the genus of $S_0$ is at least 2 ([6, 1.6.7]). Thus our minimality assumption on the complexity of $S$ implies that $S_0$ must be a torus. But then $Q_0$ is a 2-sphere and so by construction, each component of $Q_0^*$ is also a 2-sphere. As $M$ is irreducible, it follows that $Q_0^*$ is empty. Hence $S \setminus S_0$ is associated to the action of $\pi_1(M)$ on the tree, which is impossible by the minimality assumption on our choice of $S$, and so $Q_0^*$ cannot be closed.

The proof of Proposition 4.10 is now complete. □

Combining the method of the proof of Proposition 4.10 and the method of the proof of [6, Prop. 1.6.1], we see that the following useful variant of [6, Prop. 1.6.1] also holds.

PROPOSITION 4.12. *Let $X_0$ be a curve of $\bar{X}(M)$ which contains the character of an irreducible representation and let $x \in \tilde{X}_0$ be an ideal point such that $\tilde{f}_\alpha(x) \in \mathbf{C}$ for all $\alpha \in L$. Suppose further that there is a slope $r$ and some nonzero element $\delta \in L$ such that $Z_x(\tilde{f}_r) > Z_x(\tilde{f}_\delta)$. Then $M$ contains a closed, essential surface $S$ such that if $S$ compresses in $M(r)$ and $M(r')$, then $\Delta(r, r') \leq 1$.*

## 5. Culler-Shalen seminorms

In this section we introduce our main tool, the seminorms which arise from the Culler-Shalen construction. We shall also examine a natural way in which indefinite seminorms arise, and relate them to distance and incompressible surfaces.

Let $M$ denote a compact, connected, irreducible and $\partial$-irreducible 3-manifold whose boundary is a torus. Fix a curve $X_0 \subset \bar{X}(M)$ and let $\tilde{X}_0$ denote its smooth, projective model. Recall from Section 4 that if $X_0^\nu$ is the normalisation of $X_0$, then there is a birational equivalence $\nu : X_0^\nu \to X_0$ as well as a natural inclusion $X_0^\nu \subset \tilde{X}_0$ whose image consists of the *ordinary* points of $\tilde{X}_0$ and whose complement consists of the finite set of *ideal* points of $X_0$.

We consider now the construction of Section 1.4 of [6] which is used to calculate the degree of $f_\alpha | X_0$ for $\alpha \in L$. As it is easier to deal with smooth objects, we calculate this quantity by realizing it as the degree of the holomorphic map $\tilde{f}_\alpha$. Let $\{x_1, x_2, \ldots, x_n\}$ denote the set of ideal points of $\tilde{X}_0$ and recall the functions $\phi_{x_1}, \phi_{x_2}, \ldots, \phi_{x_n}$ described in Proposition 4.7. For each $\alpha \in L$, the degree of $\tilde{f}_\alpha : \tilde{X}_0 \to \mathbf{C}P^1$ is given by $\sum_{i=1}^n |\phi_{x_i}(\alpha)|$. Extend each $\phi_{x_i}$ to a linear mapping $\Phi_{x_i} : V \to \mathbf{R}$ and consider the function

$$\| \cdot \|_{X_0} : H_1(\partial M; \mathbf{R}) \to [0, \infty)$$

$$v \mapsto \sum_{i=1}^{n} |\Phi_{x_i}(v)|.$$

This function is a seminorm on $H_1(\partial M; \mathbf{R})$ which we call the *Culler-Shalen seminorm* associated to the curve $X_0$.

A seminorm on a 2-dimensional real vector space is either a norm, a nonzero indefinite seminorm, or identically zero. In our case $\|\cdot\|_{X_0}$ is induced from an integer-valued seminorm defined on $L$ for which

$$\|\alpha\|_{X_0} = \text{degree}(f_\alpha | X_0) \in \mathbf{Z} \quad \text{for each } \alpha \in L.$$

Set

$$s(X_0) = \begin{cases} \min\{\|\alpha\|_{X_0} \mid \alpha \in L, \\ \qquad f_\alpha \text{ not constant on } X_0\} & \text{if some } f_\alpha \text{ is nonconstant on } X_0 \\ 0 & \text{otherwise} \end{cases}$$

and define the *fundamental ball* of $\|\cdot\|_{X_0}$ to be

$$B_{X_0} = \{v \in V \mid \|v\|_{X_0} \leq s(X_0)\}.$$

*Definition* 5.1. Let $x$ be an ideal point of a curve $X_0 \subset \bar{X}(M)$ and suppose that the function $\phi_x$ is not identically zero. The unique slope $r$ on $\partial M$ for which $\phi_x(\alpha(r)) = 0$ is said to *associated* to $x$.

PROPOSITION 5.2. *If $\|\cdot\|_{X_0}$ is a norm on $H_1(\partial M; \mathbf{R})$, then the following hold.*

(1) *The curve $X_0$ contains the character of an irreducible representation.*

(2) *Each of the functions $f_\alpha$, $0 \neq \alpha \in L$, is nonconstant on $X_0$.*

(3) *The fundamental ball $B_{X_0}$ is a compact, convex, finite-sided polygon which is balanced, that is, $B_{X_0} = -B_{X_0}$.*

(4) *Each vertex of $B_{X_0}$ is a rational multiple of a primitive element of $L$ which corresponds to a boundary slope on $\partial M$. Furthermore, the boundary slopes that so arise are precisely the set of slopes which are associated to some ideal point $x$ of $X_0$.*

(5) *If the boundary slope associated to a vertex of $B_{X_0}$ is not a strict boundary slope, then $X_0$ is index 2 virtually abelian.*

*Proof.* Recall that the character of a reducible representation is also the character of a diagonal representation. Thus if $X_0$ does not contain an irreducible character, it lies in the image of the natural map $\bar{X}(H_1(M)) \to \bar{X}(M)$. It follows from Lefschetz duality that the inclusion induced homomorphism $i_* : H_1(\partial M) \to H_1(M)$ has rank 1 and so there is a nonzero element $\alpha \in \ker(i_*)$.

Hence as $X_0 \subset \bar{X}(H_1(M)) \subset \bar{X}(M)$, we have $f_\alpha|X_0 \equiv 0$. But then $\|\alpha\|_{X_0} = 0$, which contradicts the fact that $\|\cdot\|_{X_0}$ is a norm. This proves part (1).

Part (2) of the proposition follows from the fact that if $0 \neq \alpha \in L$, then $0 \neq \|\alpha\|_{X_0} = \text{degree}(f_\alpha|X_0)$. As for the rest, one first notes that if a norm on $\mathbf{R}^2$ is given as the sum of the absolute values of a finite number of nonzero linear functionals, then a ball of positive radius of this norm is necessarily a compact, convex, finite-sided polygon which is balanced. Moreover, the vertices of such a ball lie on the lines given by the kernels of the linear functionals. Indeed on each such line there is a $\pm$-pair of vertices of the ball. In our case the linear functionals are given by the $\phi_x$ where $x$ ranges over the ideal points of $\tilde{X}_0$ for which $\phi_x \neq 0$. Thus part (3) holds and parts (4) and (5) can be seen to be consequences of Proposition 4.7 (2). □

As an example assume that $M$ is simple and non-Seifert, so that there is a complete hyperbolic metric of finite volume in the interior of $M$. If $X_0$ is chosen to be any irreducible component of $\bar{X}(M)$ which contains the character of a discrete, faithful representation $\pi_1(M) \to \text{PSL}_2(\mathbf{C})$, then $X_0$ has dimension 1 [8] and $\|\cdot\|_{X_0}$ is a norm [6], which we shall refer to as the *Culler-Shalen norm*.

*Remark* 5.3. Curves $X_0$ in $\bar{X}(M)$ for which $\|\cdot\|_{X_0}$ is a norm are often index 2 virtually irreducible. Indeed if we suppose that the boundary of any 2-fold cover of $M$ is connected, then we may argue as in the proof of Proposition 5.2 (1) to see that for any index 2 subgroup $\tilde{\pi}$ of $\pi_1(M)$, $X_0$ contains the character of a representation which restricts to an irreducible representation on $\tilde{\pi}$. One interesting set of examples where this condition on the 2-fold covers of $M$ holds occurs when $M$ is the exterior of a knot in the 3-sphere. Note that in this case Lemma 4.6 (2) and Proposition 5.2 (5) imply that the vertices of the fundamental polygon of any norm $\|\cdot\|_{X_0}$ are strict boundary slopes.

Next we consider nonzero indefinite seminorms.

PROPOSITION 5.4. *If $\|\cdot\|_{X_0}$ is a nonzero indefinite seminorm on $H_1(\partial M; \mathbf{R})$:*

(1) *There is a unique slope $r$ on $\partial M$ for which $f_r|X_0$ is constant.*

(2) *If $\alpha \in L$, then the function $f_\alpha$ is constant on $X_0$ if and only if $\alpha = n\alpha(r)$ for some $n \in \mathbf{Z}$.*

(3) *The fundamental ball $B_{X_0}$ is an infinite band centred on the line through $0$ and $\alpha(r)$. In particular, $\|\alpha(r')\|_{X_0} = \Delta(r', r)s(X_0)$ for any $\alpha \in L$.*

(4) *The slope $r$ is a boundary slope on $\partial M$. If $r$ is not a strict boundary slope, then $X_0$ is index $2$ virtually abelian.*

(5) If $x$ is an ideal point of $\tilde{X}_0$, then either $r$ is the only slope for which $\tilde{f}_r(x) \in \mathbf{C}$, or $\tilde{f}_\alpha(x) \in \mathbf{C}$ for each $\alpha \in L$.

*Proof.* The seminorm $\|\cdot\|_{X_0}$ is induced from an integer-valued seminorm defined on $L$. Thus its restriction to $L$ must also be an indefinite seminorm. In particular, there is a unique pair $\pm\alpha$ of primitive elements of $L$ such that $\|\pm\alpha\|_{X_0} = 0$. Taking $r$ to be the slope on $\partial M$ corresponding to $\pm\alpha$, we see that part (1) of the proposition holds.

To prove part (2), we note that for each nonzero $\alpha \in L$ there is a nonzero $n \in \mathbf{Z}$ and a primitive $\alpha_0 \in L$ such that $\alpha = n\alpha_0$. Let $r_0$ be the slope associated to $\alpha_0$ and observe that $\mathrm{degree}(f_\alpha|X_0) = \|\alpha\|_{X_0} = n\|\alpha_0\|_{X_0} = n\ \mathrm{degree}(f_{r_0}|X_0)$. According to part (1) of this proposition, the latter quantity is zero if and only if $r_0 = r$. Thus part (2) holds.

Part (3) of this proposition follows from the basic properties of a nonzero indefinite seminorm and the identity $\Delta(r', r) = |\alpha(r') \cdot \alpha(r)|$. To see that part (4) holds, observe that as $\|\cdot\|_{X_0}$ is not identically zero, there is an ideal point $x$ of $\tilde{X}_0$ for which $\phi_x \not\equiv 0$. By part (1), $\tilde{f}_r$ is constant on $\tilde{X}_0$ and so in particular $\tilde{f}_r(x) \in \mathbf{C}$. Then we may apply Proposition 4.7 (2) to deduce part (4).

Finally suppose that there is an ideal point $x$ of $\tilde{X}_0$ and a slope $r' \ne r$ such that $f_{r'}(x) \in \mathbf{C}$. Now as $\tilde{f}_r$ is constant on $\tilde{X}_0$, $\tilde{f}_r(x)$ is finite as well. Thus $\phi_x(\alpha(r)) = \phi_x(\alpha(r')) = 0$. But $\alpha(r)$ and $\alpha(r')$ span a rank 2 sublattice of $L$, and so $\phi_x \equiv 0$. This implies that part (5) of the proposition holds. □

Figure 5.1 depicts $B_{X_0}$ when $\|\cdot\|_{X_0}$ is a nonzero indefinite seminorm. In this figure $r$ is the unique slope on $\partial M$ for which $\|\alpha(r)\| = 0$.

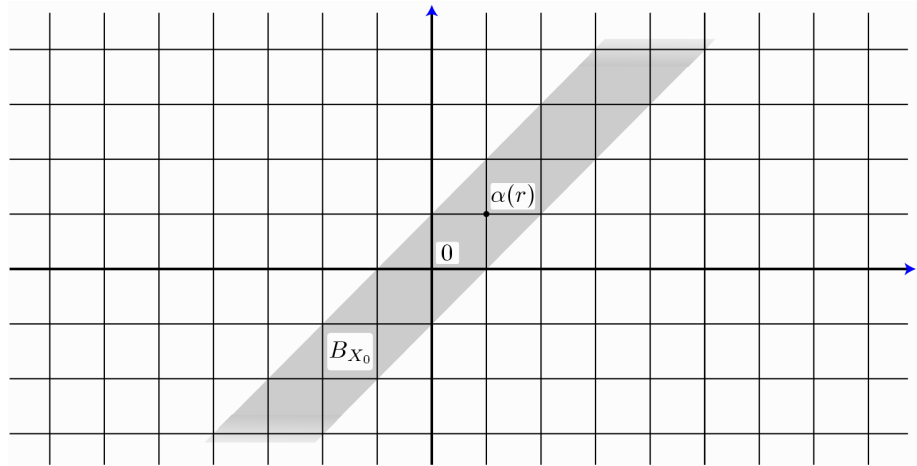

Figure 5.1: The fundamental ball of a nonzero indefinite seminorm.

PROPOSITION 5.5. *If $\|\cdot\|_{X_0}$ is identically zero, then each of the functions $f_\alpha, \alpha \in L$, is constant on $X_0$ and so there is an essential, closed surface in $M$.*

*Proof.* It is clear that each $f_\alpha$ is constant when restricted to $X_0$. In particular, for each $\alpha \in L$ and for any ideal point $x$ of $\tilde{X}_0$, $\tilde{f}_\alpha(x) \in \mathbf{C}$. Thus Proposition 4.7 (1) implies that there is an essential, closed surface in $M$. □

A natural way to construct curves $X_0$ for which $\|\cdot\|_{X_0}$ is indefinite occurs when there is a slope $r$ on $\partial M$ such that $\bar{X}(M(r))$ has positive dimension. For then the natural inclusion $\bar{X}(M(r)) \subset \bar{X}(M)$ allows us to think of any curve $X_0$ in $\bar{X}(M(r))$ as a curve in $\bar{X}(M)$. Clearly $f_r$ is identically zero on any such curve and so $\|\alpha(r)\|_{X_0} = 0$. As an example, suppose that $M(r) \cong L(p,s)\#L(q,t)$ where $1 < p, q < \infty$. According to Example 3.2, $\bar{X}(M(r))$ contains $[\frac{p}{2}][\frac{q}{2}] \geq 1$ disjoint curves. Each of them will produce an indefinite seminorm.

Our next result makes more precise the options for $\|\cdot\|_{X_0}$ when $X_0$ is a curve in $\bar{X}(M(r))$.

*Definition.* A curve $X_0 \subset \bar{X}(M)$ will be called an *r-curve* if $\|\cdot\|_{X_0}$ is nonzero and indefinite and $r$ is the unique slope for which $f_r|X_0$ is constant.

PROPOSITION 5.6. *Suppose that $X_0$ is a curve in $\bar{X}(M(r))$ which contains the character of an irreducible representation. Either*

(1) *$X_0$ is an r-curve and if $x \in \tilde{X}_0$ is an ideal point, then $r$ is the only slope $r'$ on $\partial M$ for which $\tilde{f}_{r'}(x) \in \mathbf{C}$, or*

(2) *$M$ contains a closed, essential surface $S$ associated to $x$ such that if $S$ compresses in $M(r)$ and $M(r_1)$, then $\Delta(r, r_1) \leq 1$.*

*Proof.* Consider a curve $X_0 \subset \bar{X}(M(r)) \subset \bar{X}(M)$ which contains the character of an irreducible representation and let $R_0$ be the unique subvariety of $\bar{R}(M(r)) \subset \bar{R}(M)$ such that $\bar{t}(R_0) = X_0$ (Lemma 4.1). Note that $\rho(e(\alpha(r))) = \pm I$ for each $\rho \in R_0$. If alternative (1) of Proposition 5.6 does not hold, then by Proposition 4.7 there is an ideal point $x$ of $X_0$ at which each $\tilde{f}_\alpha$ is finite, $\alpha \in L$. Hence we may now apply Proposition 4.10 to see that alternative (2) of Proposition 5.6 holds. □

A 3-manifold $W$ is called *small* if it contains no closed, essential surface. If $M$ is small, it is known that any irreducible component of $\bar{X}(M)$ is a curve [4, Prop. 2.4] and so the character varieties of small manifolds provide an excellent collection of spaces in which to consider the seminorms.

PROPOSITION 5.7. *Let $M$ be compact, connected, orientable, small 3-manifold whose boundary is a torus. Let $X_0$ be a curve in $\bar{X}(M)$.*

(1) If $x$ is an ideal point of $\tilde{X}_0$, then there is a unique slope $r$ on $\partial M$, possibly depending on $x$, such that $\tilde{f}_r(x) \in \mathbf{C}$. In particular, if $r' \neq r$ is any other slope, $f_{r'}$ is nonconstant. Thus $\|\cdot\|_{X_0}$ is nonzero.

(2) Suppose that $X_0$ is an $r$-curve and contains the character of an irreducible representation. If $M$ admits a cyclic filling slope $r_1 \neq r$, then $f_r \equiv 0$.

ADDENDUM 5.8. *Let $X_0$ be an $r$-curve of $\mathrm{PSL}_2(\mathbf{C})$-characters of the fundamental group of a small manifold $M$ and assume that it contains the character of an irreducible representation. If $r$ is a cyclic filling slope, then $M(r) \cong S^1 \times S^2$ and $X_0$ is index 2 virtually abelian.*

*Proof of Proposition* 5.7. Let $x$ be an ideal point of $\tilde{X}_0$. As $M$ is small, Proposition 4.7 shows that there is a unique slope $r$ on $\partial M$ such that $\tilde{f}_r(x) \in \mathbf{C}$. It follows that for any other slope $r'$, $\Pi_x(\tilde{f}_{r'}) > 0$, and so $f_{r'}$ is nonconstant on $X_0$. This proves part (1).

To prove part (2), assume that $X_0$ is an $r$-curve for some slope $r$ on $\partial M$, that $X_0$ contains the character of an irreducible representation, and finally that $r_1 \neq r$ is a cyclic filling slope on $\partial M$. We must find a zero of the constant function $\tilde{f}_r$ on $\tilde{X}_0$.

As $X_0$ is an $r$-curve, the hypothesis $r \neq r_1$ implies that $f_{r_1}$ is nonconstant on $X_0$, and so in particular, there exists a point $y \in \tilde{X}_0$ with $\tilde{f}_{r_1}(y) = 0$. From part (1) we see that $y$ cannot be an ideal point of $\tilde{X}_0$, and hence $y \in X_0^\nu$. By Proposition 4.8 we see that $Z_y(\tilde{f}_r) \geq Z_y(\tilde{f}_{r_1}) > 0$. Thus $\tilde{f}_r \equiv \tilde{f}_r(y) = 0$. □

*Proof of Addendum* 5.8. Suppose that $r$ is a cyclic filling slope. As $X_0$ is an $r$-curve, Proposition 5.4 (4) shows that $r$ is a boundary slope. We claim that rank $H_1(M) = 1$. To see this observe that as $M$ is small and has a connected boundary, its second Betti number must be zero. Further, since $\partial M$ is a torus, $M$ has Euler characteristic equal to 0. Therefore $0 = \mathrm{rank} H_0(M) - \mathrm{rank} H_1(M) = 1 - \mathrm{rank} H_1(M)$. Hence the hypotheses of Lemma 2.1 are satisfied by $M$. This result provides four possibilities for $M(r)$. The only one which is consistent with the hypotheses that $M$ is small and $r$ is a cyclic filling slope is the fourth one. It implies that $M(r) \cong S^1 \times S^2$ and further since $M$ is small, $r$ is a nonstrict boundary slope. Choose an ideal point $x$ on $\tilde{X}_0$ for which $\phi_x \not\equiv 0$. Then Proposition 4.7 (2) shows that $X_0$ must be index 2 virtually abelian. □

The next lemma will be used in Example 5.10 below.

LEMMA 5.9. *Let $X_0$ be a curve of $\mathrm{PSL}_2(\mathbf{C})$-characters of $\pi_1(M)$ which contains the character of an irreducible representation and let $R_0$ be the unique algebraic component of $\bar{t}^{-1}(X_0)$ for which $\bar{t}(R_0) = X_0$. Suppose that there are slopes $r$ and $r'$ on $\partial M$ such that $f_r$ is identically zero on $X_0$ while $f_{r'}$ is not.*

*Then $\rho(e(\alpha(r))) = \pm I$ for each representation $\rho \in R_0$. In particular, $X_0$ comes from a curve of $\mathrm{PSL}_2(\mathbf{C})$-characters of $\pi_1(M(r))$.*

*Proof.* By assumption $R_0$ contains an irreducible representation. Hence if we let $V$ be the set of irreducible representations in $R_0$, then $V$ is open and dense in $R_0$, and is also invariant under conjugation. It suffices to show that $\rho(e(\alpha(r))) = \pm I$ for each $\rho \in V$.

The hypothesis that $f_r \equiv 0$ on $X_0$ implies that for any $\rho \in V$, $\rho(e(\alpha(r)))$ is either $\pm I$ or a parabolic element of $\mathrm{PSL}_2(\mathbf{C})$. If the latter eventuality occurs for some $\rho_0 \in V$, we may apply [7, Prop. 1.5.4] to produce an open neighbourhood $U$ of $\rho_0$ in $V$ such that $\rho(e(\alpha(r)))$ is parabolic for each $\rho \in U$. Observe that $U$ is an open subset of the 4-dimensional variety $R_0$ (Lemma 4.1) and so contains an open 4-dimensional ball.

Since $e(\alpha(r))$ and $e(\alpha(r'))$ are commuting elements of $\pi_1(\partial M)$, it follows that for each $\rho \in U$, $\rho(e(\alpha(r'))$ is either $\pm I$ or a parabolic element of $\mathrm{PSL}_2(\mathbf{C})$. Therefore $f_{r'}$ vanishes on the set $\bar{t}(U) \subset X_0$. If $\bar{t}(U)$ consists of a finite number of characters, $\{\chi_{\rho_1}, \chi_{\rho_2}, \ldots, \chi_{\rho_n}\}$, then the irreducibility of these characters implies that $U$ is a subset of the $\mathrm{PSL}_2(\mathbf{C})$-orbits of $\rho_1, \rho_2, \ldots, \rho_n$. But this is impossible as $U$ is 4-dimensional. Hence $\bar{t}(U)$ is an infinite subset of $X_0$ on which $f_{r'}$ vanishes. It follows that $f_{r'}|X_0 \equiv 0$, which contradicts our hypotheses. Therefore we must have $\rho(e(\alpha(r))) = \pm I \in \mathrm{PSL}_2(\mathbf{C})$ for each $\rho \in R_0$ as claimed. $\square$

*Example* 5.10. Call a knot $K$ in a 3-manifold *small* if its exterior is small. Let $M$ be the exterior of a nontrivial small knot in $S^3$. Any irreducible component $X_0$ of $\bar{X}(M)$ is a curve ([4, Prop. 2.4]) and Proposition 5.7 shows that $\|\cdot\|_{X_0}$ is either a norm or a nonzero indefinite seminorm. We claim that if $X_0$ is an $r$-curve for some slope $r$ on $\partial M$, then $X_0 \subset \bar{X}(M(r))$. As a first case suppose that $X_0$ contains the character of an irreducible representation. By Addendum 5.8, $r$ cannot be the meridional slope on $\partial M$. Thus we may apply Proposition 5.7 (2) to see that $f_r|X_0 \equiv 0$. Since $\|\cdot\|_{X_0}$ is nonzero there is some slope $r'$ on $\partial M$ such that $f_{r'}|X_0$ is nonconstant. Lemma 5.9 now shows that $X_0$ comes from $\bar{X}(M(r))$. Suppose then that $X_0$ consists entirely of reducible characters. Then each point of $X_0$ is also the character of a diagonal representation. Hence $X_0$ lies in the image of $\bar{X}(H_1(M)) \to \bar{X}(M)$. As any representative of the longitudinal slope $\lambda$ of $K$ is homologically trivial in $M$, $f_\lambda|X_0 \equiv 0$. Thus $r = \lambda$ and $X_0 \subset \bar{X}(H_1(M(\lambda))) \subset \bar{X}(M(\lambda))$. $\square$

## 6. Culler-Shalen seminorms and finite or cyclic filling slopes

In this section we examine some general relations between finite or cyclic filling slopes and other exceptional slopes.

Suppose that $M$ is a compact, connected, orientable, irreducible 3-manifold whose boundary is a torus and let $X_0$ be a curve of $\text{PSL}_2(\mathbf{C})$-characters of $\pi_1(M)$. When $M$ is simple and non-Seifert and $X_0$ is the canonical curve in its $\text{PSL}_2(\mathbf{C})$-character variety, then $\|\cdot\|_{X_0}$ is the Culler-Shalen norm on $H_1(\partial M; \mathbf{R})$. In this case it is shown in [6, Cor. 1.1.4] and [3, Th. 2.3] that if $r$ is either a cyclic filling slope or a finite filling slope, but is not a strict boundary slope, then $\|\alpha(r)\|_{X_0}$ is bounded very close to $s$. In Theorem 6.2 below we show that this holds for an arbitrary curve $X_0$ of $\text{PSL}_2(\mathbf{C})$-characters.

Call a representation of a group $\Gamma$ to $\text{PSL}_2(\mathbf{C})$ *dihedral* if its image is a dihedral group. The following lemma is the key to estimating the seminorm of a finite or cyclic filling slope.

LEMMA 6.1. *Let $X_0 \subset \bar{X}(M)$ be a curve containing the character of an irreducible representation, and suppose that $r$ is a finite filling slope on $\partial M$. If $Z_x(\tilde{f}_r) > Z_x(\tilde{f}_\beta)$ for some element $\beta \in L$, and some ordinary point $x \in X_0^\nu$, then*

(1) $\nu(x)$ *is the character of an irreducible representation;*

(2) $\nu(x)$ *is a smooth point of $X_0$;*

(3) $Z_x(\tilde{f}_\beta) = 0$ *and*

$$Z_x(\tilde{f}_r) = \begin{cases} 2 & \text{if } \nu(x) \text{ is not the character of a dihedral representation} \\ 1 & \text{if } \nu(x) \text{ is the character of a dihedral representation.} \end{cases}$$

*Proof.* We shall follow the method of Section 4 of [3]. By Lemma 4.1, there is a unique algebraic component $R_0$ of $\bar{t}^{-1}(X_0)$ for which $\bar{t}(R_0) = X_0$. Moreover, $R_0$ has dimension 4. According to Proposition 4.8 there is a representation $\rho \in R_0$ which satisfies $\nu(x) = \chi_\rho$, $\rho(r) = \pm I$, and $\rho(\pi_1(M))$ is noncyclic. Since $\rho$ factors through the finite group $\pi_1(M(r))$ it then follows that it is irreducible. Arguing as in Lemmas 4.2, 4.3 and 4.4 of [3], we can calculate that the dimension of the Zariski tangent space to $R_0$ at $\nu(x)$ is 4. Thus $\rho$ is a smooth point of $R_0$. If $\rho$ is not a dihedral representation, then $\text{PSL}_2(\mathbf{C})$ acts freely on its orbit and we may use the method of [3, Lemma 4.5] to argue that $\nu(x)$ is a smooth point of $X_0$ and the Zariski tangent space of $X_0$ at $\nu(x)$ may be identified with $H^1(M; \text{Ad}\rho) \cong \mathbf{C}$. We may now proceed as in Section 4 of [3] to see that $Z_x(\tilde{f}_\beta) = 0$ and $Z_x(\tilde{f}_r) = 2$.

Assume then that $\rho$ is dihedral. In this case $\text{PSL}_2(\mathbf{C})$ does not act freely on the orbit of $\rho$ but has $\mathbf{Z}/2$ isotropy. Nevertheless, since $\rho$ is irreducible, hence stable [26, p. 53], and the orbit of $\rho$ is codimension 1 in $R_0$ and consists of smooth points, the $\text{PSL}_2(\mathbf{C})$-action on $R_0$ admits an analytic 2-disk slice $D$ at $\rho$ [26, Th. 1.2]. It follows that $\chi_\rho$ has a neighbourhood in $X_0$ which is analytically equivalent to the quotient of $D$ by the action of the $\mathbf{Z}/2$ isotropy group of $\rho$. This action is linear and so $\nu(x) = \chi_\rho$ is a smooth point of $X_0$ and $T_x X_0$ may be identified with $H^1(\pi; \mathrm{Ad} \circ \rho)/\{\pm 1\}$. This does not affect the calculation $Z_x(\tilde{f}_\beta) = 0$ as performed in Section 4 of [3], but it does affect the calculation of $Z_x(\tilde{f}_r)$. Indeed owing to the two-to-one branching of $\bar{t}|D$ we obtain $Z_x(\tilde{f}_r) = 1$ in this case. This completes the proof of Lemma 6.1. □

For each curve $X_0 \subset \bar{X}(M)$ define $L_{X_0}$ to be the subgroup of $L$ consisting of the elements whose $X_0$-seminorm is zero. Evidently

$$L_{X_0} = \begin{cases} 0 & \text{if } \|\cdot\|_{X_0} \text{ is a norm} \\ \{n\alpha(r) \mid n \in \mathbf{Z}\} & \text{if } X_0 \text{ is an } r\text{-curve} \\ L & \text{if } \|\cdot\|_{X_0} \equiv 0. \end{cases}$$

The proof of the following is essentially that of [3, Ths. 2.1 and 2.3]. One uses Proposition 4.8 and Lemma 6.1 to estimate the value of $\|\cdot\|_{X_0}$. Recall from Section 1 the terminology "C-, D-, T-, O-, I- and Q-type" finite filling slopes.

THEOREM 6.2. *Let $M$ be a compact, connected, orientable 3-manifold whose boundary is a torus and suppose that $X_0 \subset \bar{X}(M)$ is a curve which contains the character of an irreducible representation. Fix a slope $r_1$ on $\partial M$ such that $Z_x(\tilde{f}_{r_1}) \leq Z_x(\tilde{f}_\beta)$ for each ideal point $x$ of $\tilde{X}_0$ and for each $\beta \in L \setminus L_{X_0}$.*

(1) *If $r_1$ is a C-type filling slope on $\partial M$, then $\|\alpha(r_1)\|_{X_0} = s(X_0)$.*

(2) *Let $r_1$ be a D-type or a Q-type filling slope.*

  (i) *Then $\|\alpha(r_1)\|_{X_0} \leq s(X_0) + n$, where $n$ is the number of characters in $X_0$ of dihedral representations $\rho$ of $\pi_1(M)$ for which $\rho(r_1) = \pm I$.*

  (ii) *If $X_0$ is index 2 virtually irreducible, then there is an index 2 sublattice $\tilde{L} \subset L$ containing $\alpha(r_1)$ such that $\|\alpha(r_1)\|_{X_0} \leq \|\beta\|_{X_0}$ for any element $\beta \in \tilde{L} \setminus L_{X_0}$. In particular, $\|\alpha(r_1)\|_{X_0} \leq \min(s(X_0) + n, 2s(X_0))$.*

(3) *Let $r_1$ be a T-type filling slope.*

  (i) *Then $\|\alpha(r_1)\|_{X_0} \leq s(X_0) + 2$.*

  (ii) *If $X_0$ is index 2 and 3 virtually irreducible, then there is a sublattice $\tilde{L} \subset L$ of index $q = 2$ or $3$ such that $\alpha(r_1) \in \tilde{L}$ and $\|\alpha(r_1)\|_{X_0} \leq \|\beta\|_{X_0}$ for any element $\beta \in \tilde{L} \setminus L_{X_0}$.*

(4) *Let $r_1$ be an O-type filling slope.*

  (i) *Then $\|\alpha(r_1)\|_{X_0} \leq s(X_0) + 3$.*

(ii) *If we suppose that $\|\alpha(r_1)\|_{X_0} > s + 1$ and that $X_0$ is index 2, 3, and 4 virtually irreducible, then there is a sublattice $\tilde{L} \subset L$ of index $q = 2$, 3 or 4 such that $\alpha(r_1) \in \tilde{L}$ and $\|\alpha(r_1)\|_{X_0} \leq \|\beta\|_{X_0}$ for any element $\beta \in \tilde{L} \setminus L_{X_0}$.*

(5) *Let $r_1$ be an I-type filling slope.*

   (i) *Then $\|\alpha(r_1)\|_{X_0} \leq s(X_0) + 4$.*

   (ii) *If $X_0$ is index 2, 3, and 5 virtually irreducible, then there is a sublattice $\tilde{L} \subset L$ of index $q = 2$, 3, or 5 such that $\alpha(r_1) \in \tilde{L}$ and $\|\alpha(r_1)\|_{X_0} \leq \|\beta\|_{X_0}$ for any element $\beta \in \tilde{L} \setminus L_{X_0}$.*

*Remark* 6.3. There are various instances when the hypothesis on the multiplicities $Z_x(\tilde{f}_{r_1})$ for ideal points $x$ of $X_0$ given in the statement of Theorem 6.2, and also that of Corollary 6.5 below, holds. For instance it will hold if any one of the following occurs.

(i) The slope $r_1$ is not a boundary slope (Proposition 4.9).

(ii) $X_0$ is not index 2 virtually abelian and the slope $r_1$ is not a strict boundary slope (Proposition 4.9).

(iii) $X_0$ is an $r$-curve for some slope $r$ such that $\Delta(r, r_1) > 1$ and there is no closed essential surface in $M$ which remains incompressesible in $M(r)$ (Proposition 4.12).

(iv) $M$ is a small manifold, $X_0 \subset \bar{X}(M(r))$ and $r_1 \neq r$ (Proposition 5.7).

*Example* 6.4. Let $M = I(K)$, the twisted $I$-bundle over the Klein bottle. Then $M$ is small and there is a surjection $\pi_1(M) \to \mathbf{Z}/2 * \mathbf{Z}/2$ whose kernel is normally generated by the homotopy class of a loop on $\partial M$ which represents the fiber of a Seifert fibration on $I(K)$ whose base orbifold is the 2-disk with exactly two cone points of index 2 each. Let $r$ be the slope on $\partial M$ given by this fiber. According to Example 3.2 there is a curve $X_0 \subset \bar{X}(M(r)) \subset \bar{X}(M)$ which is necessarily an $r$-curve by Proposition 5.7. Now $I(K)$ admits $D$-type filling slopes whose distance to $r$ is arbitrarily large, and therefore by Theorem 6.2 (2)(ii), there must be an index 2 subgroup of $\pi_1(I(K))$ such that each $\rho \in \bar{R}(I(K))$ for which $\chi_\rho \in X_0$ restricts to a reducible representation. This is easily verified. Hence the condition that $X_0$ be index 2 virtually irreducible is necessary in Theorem 6.2 (2) (ii).

Suppose now that either $M$ is simple and non-Seifert or $M$ contains an essential torus but is not a cable on $I(K)$. It follows from the method of [3] that if there is a curve $X_0$ of $PSL_2(\mathbf{C})$-characters for which $\|\cdot\|_{X_0}$ is a norm, then the inequalities announced in Theorem 6.2 determine the bound $\Delta(r_1, r_2) \leq 5$ for any two finite or cyclic filling slopes $r_1, r_2$ on $\partial M$. We may also obtain information on finite or cyclic filling slopes when $\|\cdot\|_{X_0}$ is a nontrivial indefinite seminorm.

COROLLARY 6.5. *Let $M$ be a compact, connected, orientable 3-manifold whose boundary is a torus and suppose that there is an $r$-curve $X_0$ of $\mathrm{PSL}_2(\mathbf{C})$-characters of $\pi_1(M)$ which contains the character of an irreducible representation. Fix a slope $r_1$ on $\partial M$ such that $Z_x(\tilde{f}_{r_1}) \leq Z_x(\tilde{f}_\beta)$ for each ideal point $x$ of $\tilde{X}_0$ and for each $\beta \in L \setminus L_{X_0}$.*

(1) *If $r_1$ is a cyclic filling slope, then $\Delta(r, r_1) \leq 1$.*

(2) *If $r_1$ is a D-type or a Q-type filling slope, then $\Delta(r, r_1) \leq 1 + (n/s(X_0)) \leq 1+n$ where $n$ is the number of characters in $X_0$ of dihedral representations $\rho$ of $\pi_1(M)$ for which $\rho(r_1) = \pm I$. Further, if $X_0$ is index 2 virtually irreducible, then $\Delta(r, r_1) \leq 2$.*

(3) *If $r_1$ is a T-type filling slope, then $\Delta(r, r_1) \leq 1 + (2/s(X_0)) \leq 3$.*

(4) *If $r_1$ is an O-type filling slope, then $\Delta(r, r_1) \leq 1 + (3/s(X_0)) \leq 4$.*

(5) *If $r_1$ is an I-type filling slope, then $\Delta(r, r_1) \leq 1 + (4/s(X_0)) \leq 5$.*

*Proof.* Since $\|\alpha(r_1)\|_{X_0} = \Delta(r_1, r) s(X_0)$ (Proposition 5.4 (3)) and $s(X_0)$ is a positive integer, the corollary follows from Theorem 6.2. □

See Remark 6.3 for several conditions which guarantee that the hypothesis on $Z_x(\tilde{f}_{r_1})$ given in the statement of Corollary 6.5 holds.

The inequalities of Corollary 6.5 are sharp, as is illustrated by the following example.

*Example* 6.6. Let $M$ be the exterior of a right-handed trefoil knot and take $r$ to be the slope of the fiber of the Seifert fibering of $M$. Then $r$ corresponds to the rational number 6 under the usual identification of the slopes on the boundary of the exterior of a knot in $S^3$ with $\mathbf{Q} \cup \{\frac{1}{0}\}$ [37, Ch. 9]. Now the class $e(\alpha(r))$ is a central element of $\pi_1(M)$ and so any irreducible representation of $\pi_1(M)$ in $\mathrm{PSL}_2(\mathbf{C})$ must factor through the quotient group

$$\pi_1(M(r))/ << e(\alpha(r)) >> \; \cong \; \mathbf{Z}/2 * \mathbf{Z}/3.$$

Hence by Example 3.2, there is exactly one positive dimensional component of $\bar{X}(M)$ which contains the character of an irreducible representation, $X_0$ say, and it is isomorphic to a complex line. Note that $f_r \equiv 0$ on $X_0$ and since $M$ is small, Proposition 5.7 shows that $X_0$ is an $r$-curve. Now it follows from [32] that the finite or cyclic filling slopes on $\partial M$ coincide with the slopes $r$ which satisfy $1 \leq \Delta(r, r_0) \leq 5$. In Figure 6.1 we have depicted the lines in $H_1(\partial M); \mathbf{R})$ which contain the primitive classes in $L$ which correspond to these slopes. No finite or cyclic filling slope on $\partial M$ is a boundary slope and so by Remark 6.3 (i) the inequalities of Corollary 6.5 hold for this manifold. On the other hand, direct calculation shows that these inequalities are sharp. Notice that we must have $s(X_0) = 1$.

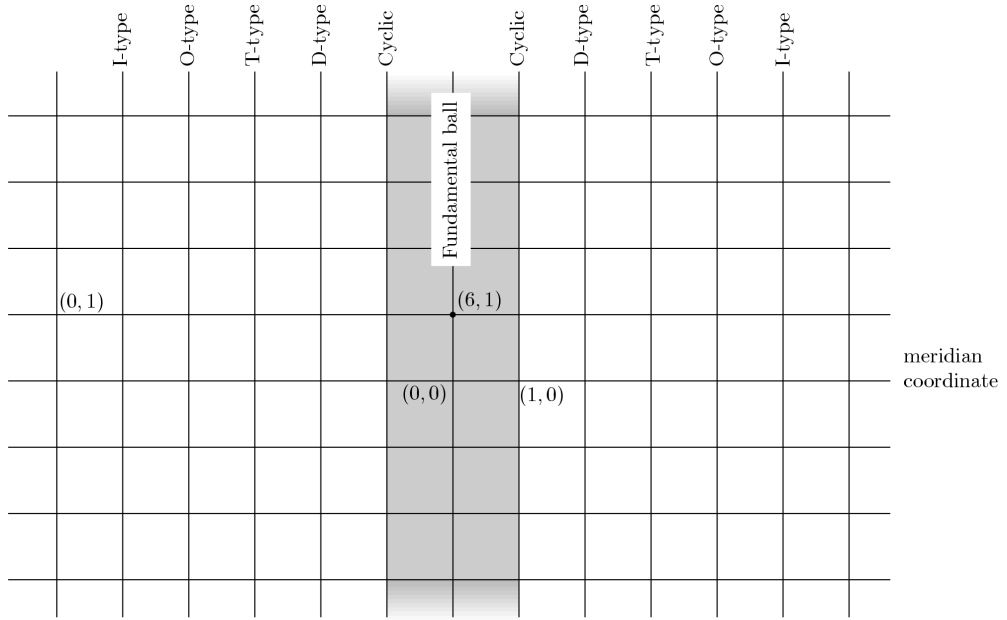

Figure 6.1: Dehn surgery on the trefoil knot.

We close this section with another corollary of Theorem 6.2.

COROLLARY 6.7. *Let $K$ be a knot in $S^3$ with exterior $M$ whose meridional slope is not a boundary slope (for instance $K$ could be a small knot [6, Th. 2.0.3]). Then any slope $r$ on $\partial M$ for which there exists some $r$-curve in $\bar{X}(M)$ is an integral slope.*

*Proof.* Let $X_0 \subset \bar{X}(M)$ be an $r$-curve. If $X_0$ consists entirely of characters of reducible representations, we use the fact that the character of a reducible representation is also the character of a diagonal one, and therefore is abelian. Hence $X_0 \subset \bar{X}(M(\lambda))$ where $\lambda$ is the longitudinal slope on $\partial M$. Then $r = \lambda$ is an integral slope. On the other hand, if $X_0$ contains the character of an irreducible representation, then we apply Remark 6.3 (i) and Corollary 6.5 (1) with $r_1$ the meridional slope to see that $r$ is necessarily an integral slope. □

## 7. Proof of Theorem 1.2

Let $M$ be a compact, connected, irreducible 3-manifold, with $\partial M$ a torus, which is neither a simple Seifert-fibered manifold nor a cable on $I(K)$. Fix slopes $r_1$ and $r_2$ on $\partial M$ and suppose that $M(r_1)$ is a reducible manifold while $M(r_2)$ has either a cyclic or a finite fundamental group. To prove Theorem 1.2 we must show that

(7.1) $\Delta(r_1, r_2) \leq 1$ if $r_2$ is a cyclic filling slope;

(7.2) $\Delta(r_1, r_2) \leq 5$ if $r_2$ is a finite filling slope as long as either $M(r_1) \neq \mathbf{R}P^3 \# \mathbf{R}P^3$ or $\pi_1(M(r_2))$ is not a $D$-type group or a $Q$-type group;

(7.3) if $M$ contains an essential torus then $\Delta(r_1, r_2) \leq 1$ as long as $M$ is neither

- a cable on $I(K)$, the twisted $I$-bundle over the Klein bottle, nor
- a cable on a simple, non-Seifert manifold $M_1$ which has filling slopes $r'_1, r'_2$ such that $M_1(r'_1) \cong \mathbf{R}P^3 \# \mathbf{R}P^3$, $M_1(r'_2)$ has a $D$-type or $Q$-type fundamental group, and $\Delta(r'_1, r'_2) \geq 8$.

*Proof of* (7.1) *and* (7.2) *when* rank $H_1(M) = 1$. According to Corollary 2.2, the reducibility of $M(r_1)$ implies that one of the following three possibilities occurs.

(7.4) $M(r_1) = L(p, s) \# L(q, t)$ is a connected sum of two lens spaces with $1 < p, q < \infty$; or

(7.5) $M$ contains a closed essential surface which remains essential in $M(r')$ whenever $\Delta(r_1, r') > 1$; or

(7.6) $M(r_1) = S^1 \times S^2$.

Suppose first of all that possibility (7.5) occurs, say $S$ is a closed surface in $M$ which remains essential in $M(r')$ as long as $\Delta(r_1, r') > 1$. As we have assumed that $\pi_1(M(r_2))$ is either cyclic or finite, $S$ must compress in $M(r_2)$ and thus $\Delta(r_1, r_2) \leq 1$. Hence (7.1) and (7.2) hold.

Next suppose that possibility (7.6) holds. Then $r_1$ is an infinite cyclic filling slope, and so when $r_2$ is a cyclic filling slope, the inequality $\Delta(r_1, r_2) \leq 1$ follows from the cyclic surgery theorem [6]. If $r_2$ is a finite filling slope, the inequality $\Delta(r_1, r_2) \leq 2$ follows from Lemma 2.4.

We may therefore suppose that $M(r_1)$ satisfies possibility (7.4), say

$$M(r_1) = L(p, s) \# L(q, t) \text{ where } 1 < p, q < \infty.$$

Note then that the following condition is satisfied.

(7.7) Any closed, essential surface in $M$ compresses in both $M(r_1)$ and $M(r_2)$.

From Example 3.2, we see that $\bar{X}(M(r_1))$ contains $[\frac{p}{2}][\frac{q}{2}] \geq 1$ curves, each one containing the character of an irreducible representation. Fix one, $X_0$ say, and note that by (7.7) and Proposition 5.6, either $\Delta(r_1, r_2) = 1$ or $\Delta(r_1, r_2) > 1$ and $X_0$ is an $r_1$-curve. Without loss of generality we may assume that the latter occurs. We claim that $r_2$ cannot be a boundary slope. Otherwise we could apply Lemma 2.1 with $r = r_2$. Consideration of our hypotheses and condition (7.7) shows that Lemma 2.1 (4) must hold. But this is impossible as it would imply that $M(r_2) \cong S^1 \times S^2$ is reducible, and so $\Delta(r_1, r_2) \leq 1$ by the main result of [18]. Thus $r_2$ is not a boundary slope. From Remark 6.3 (i) and Corollary 6.5 (1) we now see that $r_2$ cannot be a cyclic filling slope, so that (7.1) holds by default. On the other hand Corollary 6.5 also shows that (7.2) will hold as long as one of the following two conditions hold:

- $r_2$ is not a $D$-type or $Q$-type finite filling slope, or
- $X_0$ is index 2 virtually irreducible.

Now we can choose $X_0 \subset \bar{X}(M(r_1))$ so that it contains the character of a representation $\rho$ with image $< x, y \mid x^p = y^q = (xy)^7 = 1 > \cong \Delta(p, q, 7) \subset \mathrm{PSL}_2(\mathbf{C})$, the $(p, q, 7)$ triangle group. If we assume that $M(r_1) \neq \mathbf{R}P^3 \# \mathbf{R}P^3$, then $pq \geq 6$ and so $1/p + 1/q + 1/7 < 1$. Thus the image of $\rho$ is isomorphic to the orbifold fundamental group of a hyperbolic 2-orbifold. In particular, any finite index subgroup of its image is also such a group. It follows that the restriction of $\rho$ to any finite index subgroup is also irreducible. Thus the curve $X_0$ will be index 2 virtually irreducible. We conclude that the inequality $\Delta(r_1, r_2) \leq 5$ holds when either $M(r_1) \neq \mathbf{R}P^3 \# \mathbf{R}P^3$ or $r_2$ is neither a $D$-type nor a $Q$-type finite filling slope. This completes the proof of (7.2) when the first Betti number of $M$ is 1.

*Proof of (7.1) and (7.2) when* rank $H_1(M) > 1$. If rank $H_1(M) > 1$, then rank $H_1(M(r_2))$ is strictly positive, and so the only possibility is for $r_2$ to be an infinite cyclic filling slope. Any nonseparating, closed, orientable surface in $M(r_2)$ can be compressed to produce a nonseparating 2-sphere. Thus $M(r_2) = S^1 \times S^2 \# W$ where $W$ is a homotopy 3-sphere and hence $r_2$ is also a reducible filling slope. The inequality $\Delta(r_1, r_2) \leq 1$ now follows from [18]. □

*Proof of (7.3).* If $\pi_1(M(r_2)) \cong \mathbf{Z}$, then $M(r_2) \cong (S^1 \times S^2) \# W$ for some homotopy 3-sphere $W$, and so in particular $M(r_2)$ is reducible. Then the inequality $\Delta(r_1, r_2) \leq 1$ follows from [18].

Assume now that $M(r_2)$ has a finite fundamental group. It follows that each closed, orientable surface in $M$ must separate. There is also an upper bound, depending only on $M$, on the number of components of a family of disjoint essential tori in $M$, no two of which are parallel [24, Th. III.20]. Thus we can find an essential torus $T \subset \mathrm{int}\, M$ which splits $M$ into two pieces $N$ and $M_1$ where $N \cap M_1 = T$, $\partial N = \partial M \cup T$, $\partial M_1 = T$, and $M_1$ is simple.

Note that the reducibility of $M(r_1)$ implies that $N(\partial M; r_1)$ is either reducible or $\partial$-reducible. Further, $N(\partial M; r_2)$ is also $\partial$-reducible since $M(r_2)$ has a finite fundamental group.

Assume that $\Delta(r_1, r_2) > 1$. Given our hypotheses, Lemma 2.3 (1) implies that $N(\partial M; r_1)$ is irreducible. Thus $N(\partial M; r_1)$ is $\partial$-reducible and so $T$ compresses in both $N(\partial M; r_1)$ and $N(\partial M; r_2)$. It then follows from [6, Th. 2.0.1] that $N$ is a cable space, say of type $(m, n)$ where $n \geq 2$. Now [12, Lemma 7.2] shows that $N(\partial M; r_1) \cong N(\partial M; r_2) \cong S^1 \times D^2$ and so there are slopes $r'_1$ and $r'_2$ on $\partial M_1$ such that $M(r_1) = M_1(r'_1)$ and $M(r_2) = M_1(r'_2)$. Further, [12, Lemma 3.3] implies that $\Delta(r'_1, r'_2) = n^2 \Delta(r_1, r_2) \geq 8$. Hence if $M_1$ is not a Seifert manifold, (7.2) shows that $M_1(r'_1) \cong \mathbf{R}P^3 \# \mathbf{R}P^3$ and $\pi_1(M_1(r'_2))$ is a $D$-type or a $Q$-type group. If $M_1$ is a Seifert manifold, then since $M_1(r'_1) \cong M(r_1)$ is reducible and $M_1(r'_2) \cong M(r_2)$ has a finite fundamental group, Lemma 2.6 shows that $M_1(r'_1)$ is a Seifert manifold, and is therefore either $\mathbf{R}P^3 \# \mathbf{R}P^3$ or $S^1 \times S^2$. If the former case arises, the Seifert manifold $M_1$ must be $I(K)$, and so $M$ is a cable on $I(K)$. On the other hand, if the latter arises, Lemma 2.4 implies that $M$ is a cable on $I(K)$. This completes the proof of (7.3) and therefore that of Theorem 1.2. $\square$

The following example shows that Theorem 1.2 (1) is sharp when $M$ is simple and non-Seifert.

*Example* 7.8. Let $L$ be the Whitehead link. Identify the slopes on either of the components of the boundary of the exterior $Y$ of $L$ with $\mathbf{Q} \cup \{\frac{1}{0}\}$ through the use of the standard meridian-longitude coordinates. Let $M$ be the manifold obtained by Dehn filling $Y$ along one of the components of $\partial Y$ with slope 6. Now Jeff Weeks' SNAPPEA programme (available via anonymous ftp at ftp://geom.umn.edu/pub/software/snappea/) suggests that $M$ is simple and non-Seifert, and in fact, this can be verified directly. It is evident that $M(1/0)$ is the lens space $L(6, 1)$ and so the slope $\frac{1}{0}$ is a cyclic filling slope on $\partial M$. On the other hand, a simple application of the Kirby calculus shows that $M(1)$ is homeomorphic to the result of surgering the right-handed trefoil knot along the slope $\frac{6}{1}$. This latter manifold is homeomorphic to $L(3, 1) \# L(2, 1)$ ([37, 9.H.9]) and so in particular is reducible. Thus $M$ provides an example of a simple, non-Seifert manifold which admits a cyclic filling slope and a reducible filling slope of distance 1 from each other. We also note that $M(2)$ has fundamental group the finite group $D_{16} \times \mathbf{Z}/3$, where $D_{16}$ is the binary dihedral group $\{x, y; x^2 = (xy)^2 = y^4\}$ (of order 16).

*Remarks* 7.9. (1) In the case where rank $H_1(M) = 1$, our proof of parts (1) and (2) of Theorem 1.2 follows the same general outline as the original proof of Corollary 1.3 given in [17]. Both use [6, Th. 2.0.3] to deduce that the reducible manifold $M(r_1)$ is either $S^1 \times S^2$ or a connected sum of two nontrivial lens spaces. In [17] $M$ is the exterior of a knot $K$ in $S^3$ and work of D. Gabai on foliations [10] is applied to see that $M(r_1) \not\cong S^1 \times S^2$. Then the nontriviality of generalized triangle groups is invoked to conclude $r_1$ has distance 1 from the meridional slope of $K$. In our proof of Theorem 1.2 the seminorm method is used to prove the theorem when either of the possibilities for $M(r_1)$ arises. This is explicitly the case when $M(r_1)$ is a connected sum of two nontrivial lens spaces, where we use indefinite seminorms to achieve the result. When $M(r_1) \cong S^1 \times S^2$ our proof depends on the inequalities stated in Lemma 2.4. These inequalities are consequences of results in [3], where they are ultimately derived through the use of the Culler-Shalen norm determined by the canonical curve in the $\mathrm{PSL}_2(\mathbf{C})$-character variety of $M$ when $M$ is simple and non-Seifert.

(2) Examples of manifolds $M$ may be constructed which are cables on $I(K)$, and which admit a reducible filling slope $r_1$ and a finite filling slope $r_2$ such that $\Delta(r_1, r_2)$ is arbitrarily large.

## 8. Proof of Theorem 1.5

Let $M$ be a compact, connected, orientable, simple, non-Seifert 3-manifold with $\partial M$ a torus. Suppose that $r_1$ is a slope on $\partial M$ such that $M(r_1)$ has the fundamental group of a Seifert-fibered space which admits no Seifert fibration having base orbifold the 2-sphere with exactly three cone points. Let $r_2$ be a finite or cyclic filling slope on $\partial M$. In this section we shall examine the dimension of the $\mathrm{PSL}_2(\mathbf{C})$-character variety of a closed, connected, Seifert-fibered manifold. Theorem 1.5 will then follow when we combine this information with our results on Culler-Shalen seminorms to prove

(8.1) $\Delta(r_1, r_2) \leq 1$ if $M(r_2)$ has a cyclic fundamental group;

(8.2) $\Delta(r_1, r_2) \leq 5$ if $M(r_2)$ has a finite fundamental group unless $M(r_1)$ is either $\mathbf{R}P^3 \# \mathbf{R}P^3$ or a union of two copies of $I(K)$, and $\pi_1(M(r_2))$ is a $D$-type group or a $Q$-type group.

Throughout this section $W$ will denote a connected, closed Seifert-fibered space with base orbifold $\mathcal{F}$. Suppose that $W$ has precisely $q$ exceptional fibers, the $j^{\mathrm{th}}$ one of type $(\alpha_j, \beta_j)$ where $\gcd(\alpha_j, \beta_j) = 1$ and $0 < \beta_j < \alpha_j$. If $F$ denotes the surface underlying $\mathcal{F}$, we shall say that $\mathcal{F}$ is of the form $F(\alpha_1, \ldots, \alpha_q)$.

Let $g$ denote the genus of $F$ when $F$ is orientable and the maximal number of disjoint cross-caps in $F$ otherwise. The following two claims are proven in [24, VI.9-10].

- If $F$ is orientable, then there is an integer $\gamma$ for which the fundamental group of $W$ admits a presentation of the form

$$\pi_1(W) = \{a_1, b_1, ..., a_g, b_g, x_1, ..., x_q, h \mid h \text{ central}, \ x_j^{\alpha_j} = h^{\beta_j}, j = 1, ..., q,$$
$$h^\gamma = [a_1, b_1]...[a_g, b_g]x_1...x_q\}.$$

- If $F$ is nonorientable, then there is an integer $\gamma$ for which the fundamental group of $W$ admits a presentation of the form

$$\pi_1(W) = \{a_1, ..., a_g, x_1, \ldots, x_q, h \mid a_i h a_i^{-1} = h^{-1}, i = 1, \ldots, g,$$
$$x_j^{\alpha_j} = h^{\beta_j}, x_j h x_j^{-1} = h, j = 1, \ldots, q, h^\gamma = a_1^2...a_g^2 x_1...x_q\}.$$

The element $h \in \pi_1(W)$ which occurs in these presentations is represented by any regular fiber of the Seifert structure. It generates a normal cyclic subgroup $K$ of $\pi_1(W)$. Note that $K$ is central if $F$ is orientable. There is an exact sequence

(8.3) $\quad\quad\quad\quad\quad 1 \to K \to \pi_1(W) \to \pi_1^{\text{orb}}(\mathcal{F}) \to 1$

where $\pi_1^{\text{orb}}(\mathcal{F})$ is the orbifold fundamental group of $\mathcal{F}$ [42, Ch. 13].

Let $\chi(F)$ be the Euler characteristic of $F$ and recall that the orbifold Euler characteristic ([42, Ch. 13]) of the orbifold $\mathcal{F}$ is the rational number given by

$$\chi^{\text{orb}}(\mathcal{F}) = \chi(F) - \Sigma_{i=1}^q (1 - \frac{1}{\alpha_i})$$
$$= \begin{cases} 2 - 2g - \Sigma_{i=1}^q (1 - \frac{1}{\alpha_i}) & \text{if } F \text{ is orientable} \\ 2 - g - \Sigma_{i=1}^q (1 - \frac{1}{\alpha_i}) & \text{if } F \text{ is nonorientable.} \end{cases}$$

The orbifold $\mathcal{F}$ is called a *hyperbolic*, respectively *parabolic*, orbifold if it admits a hyperbolic, respectively parabolic, structure, and this condition is shown to be equivalent to the condition $\chi^{\text{orb}}(\mathcal{F}) < 0$, respectively $\chi^{\text{orb}}(\mathcal{F}) = 0$, in [42, Ch. 13]. A straightforward calculation proves that if $\mathcal{F}$ is parabolic, then it has one of the following forms.
  - a torus or a Klein bottle without cone points;
  - $S^2(2,2,2,2)$ or $S^2(p,q,r)$ where $(p,q,r)$ is one of the triples $(2,3,6)$, $(2,4,4)$ or $(3,3,3)$;
  - $\mathbf{R}P^2(2,2)$.

Thus a 2-dimensional, parabolic orbifold is finitely covered, in the sense of orbifolds, by a torus. It follows that

(8.4) a closed Seifert manifold whose orbifold is parabolic is finitely covered by an $S^1$-bundle over the torus.

We are interested in determining precise conditions for the existence of a curve in $\bar{X}(W)$ which contains the character of an irreducible representation. From Theorem 4.3 we see that if there is such a curve, then $W$ is Haken. We now examine to what extent the Haken condition suffices to guarantee the existence of such a curve.

LEMMA 8.5. *Suppose that $W$ is a closed, connected Seifert fibered manifold whose base orbifold $\mathcal{F}$ is parabolic.*

(1) *Any closed, essential surface in $W$ is a torus.*

(2) *There is a subgroup $\pi_0$ of finite index in $\pi_1(W)$ whose image under any irreducible representation $\rho : \pi_1(W) \to \mathrm{PSL}_2(\mathbf{C})$ is abelian.*

(3) *Let $\pi_0$ be a finite index subgroup of $\pi_1(W)$ as described in part (2) of this lemma. If $X_0 \subset \bar{X}(W)$ is a curve, then $\rho|\pi_0$ is reducible for each $\rho \in \bar{R}(W)$ such that $\chi_\rho \in X_0$.*

*Proof.* To prove part (1), let $\tilde{W} \to W$ be a finite cover where $\tilde{W}$ admits a locally trivial $S^1$-fiber bundle structure over the torus $\tilde{W} \to S^1 \times S^1$ (see (8.4)). Any closed, essential surface $S$ in $W$ is finitely covered by a closed, essential surface $\tilde{S}$ in $\tilde{W}$, so it suffices to prove that part (1) of the lemma holds for the latter manifold. But this is easily seen as any closed essential surface $\tilde{S} \subset \tilde{W}$ is isotopic to either a vertical surface or a horizontal one. Since the base orbifold of $\tilde{W}$ is a torus without cone points, any vertical surface is a torus. On the other hand, the restriction of the projection $\tilde{W} \to S^1 \times S^1$ to any horizontal surface in $\tilde{W}$ is a finite cover. In either eventuality we conclude that $\tilde{S}$ is a torus. This proves part (1).

To prove part (2), let $\rho : \pi_1(W) \to \mathrm{PSL}_2(\mathbf{C})$ be an irreducible representation. Recall the element $h \in \pi_1(W)$ represented by a regular fiber of $W$. We noted above that $h$ generates a normal, cyclic subgroup $K$ of $\pi_1(W)$, and so $\rho(\pi_1(W))$ is a subgroup of the normalizer of $K$ in $\mathrm{PSL}_2(\mathbf{C})$. The normalizer of a subgroup of $\mathrm{PSL}_2(\mathbf{C})$ generated by an upper triangular, parabolic element is upper triangular. Hence our assumption that $\rho$ is irreducible implies that $\rho(h)$ is diagonalisable. The normalizer of a subgroup of $\mathrm{PSL}_2(\mathbf{C})$ generated by a nontrivial diagonal element lies in

$$N = \{\pm \begin{pmatrix} z & 0 \\ 0 & z^{-1} \end{pmatrix}, \pm \begin{pmatrix} 0 & w \\ -w^{-1} & 0 \end{pmatrix} \mid z, w \in \mathbf{C} \setminus \{0\}\}.$$

The diagonal matrices form an index 2 abelian subgroup of $N$ and therefore if $\rho(h)$ is a nontrivial diagonal element of $\mathrm{PSL}_2(\mathbf{C})$, there is an index 2 subgroup of $\pi_1(W)$ such that the restriction of $\rho$ to it is abelian. On the other hand, if $\rho(h) = \pm I$, consider $\pi_1(\tilde{W}) \subset \pi_1(W)$ where $\tilde{W}$ is an $S^1$-bundle over a torus $T$ which admits a finite covering map to $W$. As $\rho(h) = \pm I$, $\rho|\pi_1(\tilde{W})$ factors through the abelian group $\pi_1(T) \cong \mathbf{Z} \oplus \mathbf{Z}$ (see (8.3)). Choose $\pi_0$ to be any finite index subgroup of $\pi_1(W)$ which is contained in the intersection of $\pi_1(\tilde{W})$ with all the index 2 subgroups of $\pi_1(W)$. Then our analysis above shows that $\pi_0$ satisfies the conclusion of part (2) of the lemma.

Finally we prove part (3). If $X_0$ consists of reducible characters, then the result is obvious. On the other hand, if it contains the character of an irreducible representation, then the argument from the proof of Lemma 4.6 (2) shows that $\rho|\pi_0$ is reducible for each $\rho \in \bar{R}(W)$ such that $\chi_\rho \in X_0$. □

LEMMA 8.6. *Suppose that $W$ is a closed, connected Seifert-fibered manifold.*

(1) *If $W$ is neither $S^1 \times S^2$ nor $\mathbf{R}P^3 \# \mathbf{R}P^3$, then it is irreducible.*

(2) *If $W$ is Haken, then $\mathcal{F}$ is either hyperbolic or parabolic.*

(3) *If $\mathcal{F}$ is parabolic, then $W$ is either Haken or $H_1(W)$ is finite and $\mathcal{F}$ is of the form $S^2(p,q,r)$. If $W$ is Haken, it is either an $S^1$-bundle over the torus, a union of two copies of $I(K)$, or $\mathcal{F} = S^2(p,q,r)$ and $W$ is a torus bundle over $S^1$.*

(4) *If $\mathcal{F}$ is hyperbolic, then either $W$ is Haken or $H_1(W)$ is finite and $\mathcal{F}$ is of the form $S^2(p,q,r)$.*

*Proof.* The proof of this lemma is contained for the most part in [24, VI.7] and [23, VI.15]. The fact that $W$ is of the form described in part (3) when $W$ is Haken and $\mathcal{F}$ is parabolic, follows from the list of the closed parabolic 2-orbifolds given above as well as Lemma 8.5 (1) and [24, VI.34]. □

Now it can be shown that the $\mathrm{PSL}_2(\mathbf{C})$-character variety of a Seifert-fibered space is positive dimensional if and only if this space is either Haken or reducible. The next lemma determines exactly when we can find a curve of such characters which contains the character of an irreducible representation.

LEMMA 8.7. *Let $W$ be a connected, closed Seifert-fibered space with base orbifold $\mathcal{F}$.*

(1) *There is a curve in $\bar{X}(W)$ which contains the character of an irreducible representation if and only if*

   (i) $W = \mathbf{R}P^3 \# \mathbf{R}P^3$ *or*

   (ii) *$W$ is Haken but $\mathcal{F}$ is neither a torus without cone points nor a 2-sphere with three cone points.*

(2) *There is a curve in $\bar{X}(W)$ which is index $q$ virtually irreducible for each $q \geq 1$ if and only if $W$ is Haken and $\mathcal{F}$ is a hyperbolic orbifold other than a 2-sphere with three cone points.*

*Proof.* Assume that there is a curve $X_0$ in $\bar{X}(W)$ which contains the character of an irreducible representation. Our first goal is to show that either $W = \mathbf{R}P^3 \# \mathbf{R}P^3$ or $W$ is Haken such that $\mathcal{F}$ is neither a torus without cone points nor a 2-sphere with three cone points.

To that end, note that by Theorem 4.3, $W$ admits an essential surface. Let $S$ be such a surface, chosen to have minimal genus. If $S$ is a 2-sphere then $W$ is reducible and so is either $S^2 \times S^1$ or $\mathbf{R}P^3 \# \mathbf{R}P^3$ (Lemma 8.6 (1)). Clearly $W$ cannot be $S^2 \times S^1$, and so is $\mathbf{R}P^3 \# \mathbf{R}P^3$. If the genus of $S$ is positive, then $W$ is Haken. Now $X_0$ contains a Zariski open set of irreducible characters. Further, there are only finitely many characters in $X_0$ corresponding to representations whose image is $\mathbf{Z}/2 \oplus \mathbf{Z}/2$. Thus there is a Zariski open set of irreducible characters in $X_0$ corresponding to representations with a nonabelian image. Suppose that $F$ is orientable and let $\rho \in \bar{R}(W)$ be any nonabelian, irreducible representation whose character lies in $X_0$. The fiber class $h \in \pi_1(W)$ is central and so $\rho(h) = \pm I$. Thus $\rho$ factors through $\pi_1^{\mathrm{orb}}(\mathcal{F})$ (see (8.3)). It follows that $X_0 \subset \bar{X}(\pi_1^{\mathrm{orb}}(\mathcal{F}))$. But then $\mathcal{F}$ cannot be a torus without cone points, as its orbifold fundamental group supports no nonabelian representations. Further, $\mathcal{F}$ is not of the form $S^2(p,q,r)$, as $\pi_1^{\mathrm{orb}}(S^2(p,q,r))$ is the $(p,q,r)$ triangle group (see (8.3)), which has a finite $\mathrm{PSL}_2(\mathbf{C})$-character variety.

Suppose now that there is a curve $X_0$ in $\bar{X}(W)$ which is index $q$ virtually irreducible for each $q \geq 1$. From the previous paragraph, we see that $W$ is either Haken or $\mathbf{R}P^3 \# \mathbf{R}P^3$. The latter case cannot occur, for then $\pi_1(W) \cong \mathbf{Z}/2 * \mathbf{Z}/2$ and so has as an index-2 normal subgroup isomorphic to $\mathbf{Z}$. Hence by Lemma 4.6 (2) $X_0$ would be index-2 virtually reducible, contrary to our assumptions. Thus $W$ is Haken and so by Lemmas 8.6 (2) and 8.5 (3) we see that $\mathcal{F}$ is hyperbolic.

Next we shall construct curves in $\bar{X}(W)$ which contain irreducible characters.

The fundamental group of the manifold $\mathbf{R}P^3 \# \mathbf{R}P^3$ is $\mathbf{Z}/2 * \mathbf{Z}/2$ and so by Example 3.2 there is a curve in $\bar{X}(\mathbf{R}P^3 \# \mathbf{R}P^3)$ containing the character of an irreducible representation. Suppose then that the manifold $W$ is Haken. The orbifold $\mathcal{F}$ is necessarily parabolic or hyperbolic. If it is parabolic but not a torus without cone points or a 2-sphere with exactly three cone points, then Lemma 8.6 (3) shows that $\pi_1(W)$ surjects onto $\mathbf{Z}/2 * \mathbf{Z}/2$. Thus we may again apply Example 3.2 to find a curve in $\bar{X}(W)$ which contains the character of an irreducible representation.

Finally suppose that $W$ is Haken and $\mathcal{F}$ is hyperbolic, but is not a 2-sphere with three cone points. To complete the proof of the lemma it suffices to produce a representation $\rho \in \bar{R}(W)$ whose restriction to each finite index subgroup of $\pi_1(W)$ is irreducible, and whose character lies in a positive dimensional component of $\bar{X}(W)$.

To that end, let $\rho$ be a discrete, faithful, representation of $\pi_1^{\mathrm{orb}}(F)$ into $\mathrm{PSL}_2(\mathbf{R}) \subset \mathrm{PSL}_2(\mathbf{C})$. This is a holonomy representation of some hyperbolic structure on $\mathcal{F}$. Now there is a natural inclusion $\bar{R}(\pi_1^{\mathrm{orb}}(\mathcal{F})) \subset \bar{R}(W)$ (see (8.3)) and so we shall consider $\rho$ as an element of $\bar{R}(W)$. Notice that any finite orbifold cover of $\mathcal{F}$ is closed and hyperbolic, and so the restriction of $\rho$ to any finite index subgroup of $\pi_1^{\text{orb}}(\mathcal{F})$ is irreducible. To see that $\chi_\rho$ lies in a positive dimensional component of $\bar{X}(W)$, consider the Teichmüller space $\mathcal{T}_\mathcal{F}$ of the hyperbolic orbifold $\mathcal{F}$. According to [42, Cor. 13.3.7], this is a Euclidean space of dimension

$$-3\chi(F) + 2q = \begin{cases} -6 + 6g + 2q & \text{if } F \text{ is orientable,} \\ -6 + 3g + 2q & \text{if } F \text{ is nonorientable.} \end{cases}$$

It is simple to show that the only hyperbolic orbifolds for which this quantity is zero are orbifolds of the form $S^2(p,q,r)$, and as we have purposefully excluded this possibility, $\mathcal{T}_\mathcal{F}$ is positive dimensional. Now there is a function $\mathcal{T}_\mathcal{F} \to \bar{X}(\pi_1^{\text{orb}}(\mathcal{F})) \subset \bar{X}(W)$ which sends a hyperbolic structure to the character of its holonomy representation, and this function is injective. Thus $\chi_\rho$ is contained in a positive dimensional component of $\bar{X}(W)$. □

*Proof of* (8.1) *and* (8.2). Recall the notation from the beginning of this section. Observe that if $M(r_1)$ is reducible, then (8.1) and (8.2) follow from Theorem 1.2. Hence, without loss of generality, we can assume that $M(r_1)$ is irreducible.

Suppose, first of all, that $M(r_1)$ has the fundamental group of a Seifert-fibered space admitting a fibering over the 2-sphere with at most two cone points. Then $\pi_1(M(r_1))$ is a cyclic group and so in this case, (8.1) and (8.2) are consequences of the cyclic surgery theorem [6] and the finite surgery theorem [3, Th. 1.1]. Hence we may assume below that $M(r_1)$ has the fundamental group of a Seifert-fibered space which admits no fibering over the 2-sphere with *at most* three cone points. Then in particular, $\pi_1(M(r_1))$ is infinite. Peter Scott has proven [39] that a closed, irreducible 3-manifold whose fundamental group is infinite and isomorphic to the fundamental group of a Seifert-fibered manifold, is in fact Seifert fibered. Thus $M(r_1)$ admits a Seifert fibering. Let $\mathcal{F}$ be its base orbifold. A simple calculation shows that a Seifert manifold whose base orbifold has a positive orbifold Euler characteristic is either reducible or has a finite fundamental group. Thus $\chi^{\text{orb}}(\mathcal{F}) \leq 0$, and hence $\mathcal{F}$ is either parabolic or hyperbolic.

If $\mathcal{F}$ is parabolic, then given our hypothesized constraints on $M(r_1)$, we may apply Lemma 8.6 (3) to see that $M(r_1)$ is either an $S^1$-bundle over the torus or the union of two copies of $I(K)$. Assume first of all that the former case arises and that the $S^1$-bundle has Euler number $b \in \mathbf{Z}$. Then $H_1(M(r_1)) \cong \mathbf{Z} \oplus \mathbf{Z} \oplus \mathbf{Z}/b$, and so rank $H_1(M) \geq 2$. It follows that there can never be a finite filling of $M$, and so $r_2$ is an infinite cyclic filling slope. But for this to occur we must have $H_1(M) \cong \mathbf{Z} \oplus \mathbf{Z}$ and therefore the natural homomorphism $H_1(M) \to H_1(M(r_1))$ is an isomorphism. Then if $\alpha(r_1), \alpha(r_2)$ are primitive classes associated to the slopes $r_1, r_2$, the image of $\alpha(r_1)$ in $H_1(M)$ must be trivial. Choose $\alpha_2 \in H_1(\partial M)$ which forms a basis of $H_1(\partial M)$ with $\alpha(r_1)$. Write $\alpha(r_2) = m\alpha(r_1) + n\alpha_2$ and notice that $|n| = \Delta(r_1, r_2)$. Then $\mathbf{Z} \cong H_1(M(r_2)) \cong H_1(M)/ < n\alpha_2 > \cong (\mathbf{Z} \oplus \mathbf{Z})/ < n\xi >$ for some $\xi \in \mathbf{Z} \oplus \mathbf{Z}$. Clearly then $\Delta(r_1, r_2) = |n| = 1$. Thus both (8.1) and (8.2) hold.

Next assume that $M(r_1)$ is the union of two copies of $I(K)$. Then Lemma 8.7 (1) implies that there is a curve $X_0 \subset \bar{X}(M(r_1)) \subset \bar{X}(M)$ which contains the character of an irreducible representation. If there is an ideal point $x \in \tilde{X}_0$ such that $\tilde{f}_\alpha(x) \in \mathbf{C}$ for all $\alpha \in H_1(\partial M)$, then according to Proposition 5.6, there is an essential, closed surface $S \subset M$ which cannot compress in both $M(r_1)$ and $M(r)$ for any slope $r$ on $\partial M$ for which $\Delta(r_1, r) > 1$. Now under the hypotheses of either (8.1) or (8.2), $S$ compresses in $M(r_2)$. Further, since $M$ is simple, the genus of $S$ is at least 2. But then Lemma 8.5 (1) implies that $S$ compresses in $M(r_1)$ and so $\Delta(r_1, r_2) \leq 1$. Thus both (8.1) and (8.2) hold. Finally assume that for each ideal point $x$ of $\tilde{X}_0$, $f_{r'}$ has a pole at $X$ for each slope $r' \neq r_1$. In particular $Z_x(f_{r_2}) = 0$. An application of Corollary 6.5 (where $r_1, r_2$ of our current situation correspond respectively to $r, r_1$ of Corollary 6.5) completes the proof of (8.1) and (8.2) when $\mathcal{F}$ is parabolic.

Assume now that $\mathcal{F}$ is hyperbolic. According to Lemma 8.6 (4), $M(r_1)$ is Haken, and since we have assumed that $M(r_1)$ does not admit a Seifert structure whose base orbifold is of the form $S^2(p, q, r)$, we can apply Lemma 8.7 (2) to find a curve $X_0 \subset \bar{X}(\pi_1^{\mathrm{orb}}(\mathcal{F})) \subset \bar{X}(M(r_1)) \subset \bar{X}(M)$ which is index $q$ virtually irreducible for each $q \geq 1$.

First assume that there is an ideal point $x \in \tilde{X}_0$ such that $\tilde{f}_\alpha(x) \in \mathbf{C}$ for all $\alpha \in H_1(\partial M)$. According to Proposition 5.6, there is an essential, closed surface $S$ in $M$, associated to the ideal point $x$, which cannot compress in both $M(r_1)$ and $M(r)$ for any slope $r$ on $\partial M$ for which $\Delta(r_1, r) > 1$.

*Claim.* $S$ compresses in $M(r_1)$.

*Proof of the claim.* If we assume that $S$ remains incompressible in $M(r_1)$, then with respect to the given Seifert structure on $M(r_1)$, $S$ may be isotoped to be either horizontal or vertical ([24, VI.34]). In the latter eventuality $S$ would necessarily be a torus, which is ruled out by the fact that $M$ is simple. Thus $S$ must be horizontal. But then by Lemma 2.7, $\pi_1(S)$ is a normal subgroup of $\pi_1(M(r_1))$ with $\pi_1(M(r_1))/\pi_1(S) \cong \mathbf{Z}$ or $\mathbf{Z}/2 * \mathbf{Z}/2$. Set $\Gamma = \pi_1(M(r_1))$ and let $R_0 \subset \bar{R}(\Gamma)$ be the unique subvariety for which $\bar{t}(R_0) = X_0$ (Lemma 4.1). Recall from Section 4 that there exist

- a central $\mathbf{Z}/2$-extension $\phi : \hat{\Gamma} \to \Gamma$,

- a subvariety $S_0 \subset R(\hat{\Gamma})$ and a regular mapping $\phi_* : S_0 \to R_0$ where $\phi_*(\hat{\rho})$ is the unique homomorphism $\rho \in R_0$ for which $\hat{\rho} = \phi_*(\hat{\rho}) \circ \phi$, and

- tautological representations $\hat{P} : \hat{\Gamma} \to \mathrm{SL}_2(F)$ and $P : \Gamma \to \mathrm{PSL}_2(F)$ where $F$ is the function field of $S_0$ and such that $\hat{P}(\hat{\gamma}) = P(\phi(\hat{\gamma}))$ for each $\hat{\gamma} \in \hat{\Gamma}$.

Now by construction, $\pi_1(S)$ lies in an edge stabiliser of the action of $\pi_1(M(r_1))$ on some tree associated to the ideal point $x$. Hence Proposition 4.4 implies that either $P(\pi_1(S)) = \{\pm I\}$ or there is an index-2 subgroup of $\Gamma$ such that $P|\Gamma_0$ is diagonalisable. Note that the latter case cannot arise because if it did, Lemma 4.6 would imply that $X_0$ is index-2 virtually reducible, contrary to our choice of $X_0$. Hence $P(\pi_1(S)) = \{\pm I\}$. But then $P$ factors through $\pi_1(M(r_1))/\pi_1(S) \cong \mathbf{Z}$ or $\mathbf{Z}/2 * \mathbf{Z}/2$. It follows that there is an index-2 subgroup $\Gamma_0$ of $\Gamma$ such that $P|\Gamma_0$ is reducible. Lemma 4.6 now implies that $X_0$ is index-2 virtually reducible, contrary to our hypothesis. Thus $S$ must compress in $M(r_1)$. □

Observe that our hypotheses imply that $S$ compresses in $M(r_2)$. Thus from the choice of $S$ we have $\Delta(r_1, r_2) \leq 1$, so that both (8.1) and (8.2) hold.

As a final case for consideration, suppose that for any ideal point $x \in \tilde{X}_0$, there is some $\alpha \in H_1(\partial M)$ for which $\tilde{f}_\alpha$ has a pole at $x$. Such an $f_\alpha$ is nonconstant on $X_0$ and so as $f_{r_1}|X_0 \equiv 0$, it follows that $X_0$ is an $r_1$-curve. Fix an ideal point $x$ of $\tilde{X}_0$. Now $\tilde{f}_{r_1}(x) = 0 \in \mathbf{C}$ and so according to Proposition 4.7, $\tilde{f}_r$ has a pole at $x$ for each slope $r \neq r_1$ on $\partial M$. If $\beta \in H_1(\partial M) \setminus \{n\alpha(r_1) \mid n \in \mathbf{Z}\}$, then there is a slope $r \neq r_1$ such that $\beta = n\alpha(r)$ for some nonzero $n \in \mathbf{Z}$. Now Proposition 4.7 also tells us that $\tilde{f}_\beta$ has a pole of degree $|n|\Pi_x(\tilde{f}_r) > 0$ at $\tilde{x}$, so in particular $Z_x(\tilde{f}_\beta) = 0 < Z_x(\tilde{f}_{r_1})$. The hypotheses of Corollary 6.5 are therefore satisfied (where $r_1, r_2$ of our current situation correspond respectively to $r, r_1$ of Corollary 6.5) and as $X_0$ is index-2 virtually irreducible, (8.1) and (8.2) follow from this corollary. This completes the proof of Theorem 1.5. □

## 9. Proof of Theorem 1.6

Let $M$ be a non-Seifert, compact, connected, orientable, irreducible 3-manifold with $\partial M$ a torus and suppose that $M$ contains an essential torus. Consider two slopes $r_1$ and $r_2$ on $\partial M$ such that $M(r_1)$ has the fundamental group of a Seifert-fibered space and $M(r_2)$ has a finite or cyclic fundamental group. To prove Theorem 1.6 we must show that

(9.1) If $\Delta(r_1, r_2) > 1$ then $M$ is a cable on a manifold $M_1$ which is either simple or Seifert fibered. Furthermore, $M_1$ admits a finite or cyclic filling according to whether $r_2$ is a finite or a cyclic filling slope.

(9.2) If $\Delta(r_1, r_2) > 1$ where $r_2$ is a cyclic filling slope and $M(r_1)$ has the fundamental group of a Seifert-fibered space which admits no Seifert fibration having base orbifold the 2-sphere with exactly three cone points, then $M$ is a cable on a Seifert manifold admitting a cyclic filling.

We may suppose below that $M(r_1)$ is irreducible, as otherwise the result follows from Theorem 1.2. We may also suppose that $M(r_1)$ has an infinite fundamental group by [3, Th. 1.2]. Thus $M(r_1)$ can be taken to be a Seifert-fibered manifold [39]. It will be assumed below that $M(r_1)$ is endowed with a fixed Seifert structure.

The verification that (9.1) and (9.2) hold will be accomplished by examining how the essential tori in $M$ behave in the Seifert manifold $M(r_1)$. If $T$ is such a torus, it will either compress in $M(r_1)$ or be isotopic to an essential horizontal torus or an essential vertical torus. Most of the work which follows involves an analysis of the latter two cases.

For the remainder of this section $N$ will denote the piece of the torus decomposition of $M$ which contains $\partial M$. Note that $N$ has at least two boundary components and is either Seifert-fibered or simple and non-Seifert. Further, it is possible that some component of $\partial N \setminus \partial M$ corresponds to a nonseparating torus in $\text{int}(M)$. Hence it may occur that distinct boundary components of $N$ are identified under the natural map $N \to M$. The image of $N$ in $M$ will be denoted by $N_M$.

LEMMA 9.3. (1) *If $T$ is a separating, essential torus in the interior of $M$, say $M = P \cup_T M_1$ where $\partial P = \partial M \cup T$ and $\partial M_1 = T$, then $M(r_2) \cong M_1(r_2')\#V$ for some slope $r_2'$ on $T$ and for some closed 3-manifold $V$.*

(2) $N(\partial M; r_2) \cong (S^1 \times D^2)\#W$ *for some 3-manifold $W$.*

(3) *Suppose that $N$ is Seifert fibered and that $\phi_N$ is the slope on $\partial M$ corresponding to the fiber of $N$. If $r_2 \neq \phi_N$, then $N$ is a cable space and $N(\partial M; r_2) \cong S^1 \times D^2$.*

*Proof.* Let $T$ be a torus as described in the hypotheses of part (1). Clearly $T$ compresses in $M(r_2)$ and therefore it must compress in $P(\partial M; r_2)$. Thus $P(\partial M; r_2) \cong (S^1 \times D^2)\#V$ for some closed 3-manifold $V$. Let $r_2'$ be the slope on $T$ corresponding to the meridian of the $S^1 \times D^2$ factor of $P(\partial M; r_2)$. Then $M(r_2) \cong M_1(r_2')\#V$ and hence (1) holds.

To prove (2), just observe that as the image of each component of $\partial(N(\partial M; r_2))$ compresses in $M(r_2)$, at least one of them compresses in $N(\partial M; r_2)$. Hence $N(\partial M; r_2) \cong (S^1 \times D^2)\#W$ for some 3-manifold $W$.

Now suppose that $N$ is Seifert fibered and that $r_2 \neq \phi_N$. Then the Seifert structure on $N$ extends over $N(\partial M; r_2)$ and so as $\partial(N(\partial M; r_2))$ must compress in $N(\partial M; r_2)$, $N(\partial M; r_2)$ is a solid torus. This implies that $N$ is a cable space. □

Write $\partial N = \partial M \cup T_1 \cup \ldots \cup T_m$ where $T_1, T_2, \ldots, T_m$ are tori. By Lemma 9.3 (2) we have $N(\partial M; r_2) \cong (S^1 \times D^2) \# W$, and without loss of generality we may suppose that $T_1$ is the boundary of the $(S^1 \times D^2)$ factor. Let $s_1$ be the slope on $T_1$ which corresponds to the meridian of this solid torus. Observe that $\partial W = T_2 \cup \ldots \cup T_m$.

LEMMA 9.4. (1) *If the image of some component $T_i$ of $\partial N \setminus \partial M$ compresses in $M(r_1)$, then either $\Delta(r_1, r_2) \leq 1$ or $M$ is cabled.*

(2) *Suppose that the image in $M(r_1)$ of $T_1 \cup T_2 \cup \ldots \cup T_m$ forms a set of essential tori. Then*

   (i) *this image is isotopic in $M(r_1)$ to either a collection of vertical tori or a collection of horizontal tori;*

   (ii) *$N(\partial M; r_1)$ and $\overline{M \setminus N_M}$ are Seifert-fibered.*

(3) *If the image in $M(r_1)$ of $T_1 \cup T_2 \cup \ldots \cup T_m$ is isotopic to a set of horizontal tori, then*

   (i) *$N(\partial M; r_1) \cong I(K)$ or $N(\partial M; r_1) \cong S^1 \times S^1 \times I$;*

   (ii) *$\overline{M \setminus N_M}$ is a possibly empty union of copies of $I(K)$.*

(4) *Suppose that $N$ is Seifert-fibered and $\phi_N$ is the slope on $\partial M$ corresponding to the fiber of this structure. If the image in $M(r_1)$ of $T_1 \cup T_2 \cup \ldots \cup T_m$ is a set of essential tori, then either*

   (i) *$\Delta(r_1, r_2) \leq 1$, or*

   (ii) *$N$ is a cable space, $N(\partial M; r_1) \cong I(K)$, $\Delta(r_1, \phi_N) = 2$, and $M = N \cup M_1$ where $M_1$ is Seifert-fibered and admits a finite or cyclic filling slope according to whether $r_2$ is a finite or a cyclic filling slope.*

*Proof.* The image of $T_i$ in $M(r_2)$ is compressible so part (1) follows from [6, Th. 2.0.1].

According to [24, VI.34], the hypotheses of part (2) imply that the image in $M(r_1)$ of each component of $\partial N \setminus \partial M$ is isotopic to either a vertical surface or a horizontal surface. Since the image of $\partial N \setminus \partial M$ in $M(r_1)$ is a disjoint family of tori, the first part of (2) now follows from an application of Lemma 2.8 (1).

Suppose now that the image $T_i'$ of some component $T_i$ of $\partial N \setminus \partial M$ is isotopic to a horizontal surface in $M(r_1)$. According to [24, VI.34], either $M(r_1) \cong I(K) \cup_{T_i'} I(K)$ or $M(r_1)$ is the total space of a locally trivial torus bundle over the circle for which $T_i'$ is a fiber. In either case, any essential torus in $M(r_1)$ which is disjoint from $T'_i$ is parallel to it (Lemma 2.8 (2)). Thus if $m > 1$ we must have $m = 2$ and $N(\partial M; r_1) \cong S^1 \times S^1 \times I$. In this case, $\overline{M \setminus N_M} = \overline{M(r_1) \setminus N_M(\partial M; r_1)}$ is either empty or $I(K) \cup I(K)$, and so is Seifert fibered. On the other hand if $m = 1$, then $T'_1$ separates $M(r_1)$ and so we necessarily have $M = N \cup_{T'_1} I(K)$ and $N(\partial M; r_1) \cong I(K)$. This proves part (3) of the lemma and the second part of (2) when the tori from $\partial N \setminus \partial M$ are isotopic to horizontal surfaces. To complete the proof of the second part of (2) simply observe that if the image in $M(r_1)$ of $T_1 \cup T_2 \cup \ldots \cup T_m$ is isotopic to a set of essential, vertical tori, then the Seifert structure on $M(r_1)$ restricts to one on $\overline{M \setminus N_M}$ and pulls back to one on $N(\partial M; r_1)$ (recall that in general, $N \not\subset M$).

Finally we prove (4). Assume that $N$ is Seifert-fibered and that the image in $M(r_1)$ of $T_1 \cup T_2 \cup \ldots \cup T_m$ is a set of essential tori. First we claim that either $N(\partial M; r_1) \cong I(K)$ or $N(\partial M; r_1) \cong S^1 \times S^1 \times I$. By (3), this holds if the tori from $\partial N \setminus \partial M$ are horizontal in $M(r_1)$, so suppose that they are vertical. As we noted in the previous paragraph, this condition implies that the Seifert structure on $M(r_1)$ restricts to one on $\overline{M \setminus N_M}$ and pulls back to one on $N(\partial M; r_1)$. We can construct another Seifert structure on $N(\partial M; r_1)$ as follows. The slope $r_1$ cannot be $\phi_N$, the slope on $\partial M$ corresponding to the fiber of the given Seifert structure on $N$, because otherwise the image of each component of $\partial N \setminus \partial M$ would compress in $M(r_1)$. Hence the Seifert structure on $N$ extends over $N(\partial M; r_1)$.

Observe that these two Seifert fiberings on $N(\partial M; r_1)$ cannot have isotopic fibers on $\partial(N(\partial M; r_1)) \neq \emptyset$, for if this were to occur, a Seifert structure on $M$ could be constructed, contrary to our hypotheses. Thus $N(\partial M; r_1)$ admits two Seifert structures with nonisotopic fibers on at least one component of its boundary. The only possibility is for $N(\partial M; r_1)$ to be either $S^1 \times D^2, I(K)$, or $S^1 \times S^1 \times I$ [24, Lemma VI.20]. We may exclude the first possibility by the assumed $\partial$-irreducibility of $N(\partial M; r_1)$.

Now we return to the proof of (4). Observe that if $r_2 \neq \phi_N$, then Lemma 9.3 (3) shows that $N$ is a cable space and $N(\partial M; r_2) \cong S^1 \times D^2$. In particular, $M = N \cup_{T_1} M_1$ is a cable on a Seifert manifold $M_1$ (by part (2)) and $M_1(s_1) \cong M(r_2)$ where $s_1$ is the slope on $\partial M_1$ corresponding to the meridian of $N(\partial M; r_2)$. Since $N$ is a cable space, consideration of the constraints on $N(\partial M; r_1)$ determined above now shows that $N(\partial M; r_1) \cong I(K)$. But for this to occur $N$ must be an $(m, 2)$-cable space for some odd integer $m$, and $\Delta(r_1, \phi_N) = 2$.

The final case to consider is when $r_2 = \phi_N$. Then $r_1 \neq \phi_N$ and so the Seifert structure on $N$ extends over $N(\partial M; r_1)$. In particular the latter structure has a fiber of multiplicity $\Delta(r_1, \phi_N) = \Delta(r_1, r_2)$. Consideration of the multiplicities of the exceptional fibers of the possible Seifert fibrations on $I(K)$ and $S^1 \times S^1 \times I$ [24, VI.17] shows that $\Delta(r_1, r_2) = \Delta(r_1, \phi_N) \leq 2$ and if it equals 2, then $N$ is a cable space and $N(\partial M; r_1) = I(K)$. Since we have assumed that the image of $\partial N \setminus \partial M$ is essential in $M(r_1)$, it now follows from part (2) of this lemma that $M = N \cup M_1$ where $M_1$ is Seifert-fibered. Now as $r_2 = \phi_N$, $N(\partial M; r_2) \cong (S^1 \times D^2) \# \mathbf{R}P^2$, and hence $M(r_2) \cong \mathbf{R}P^2 \# M_1(r_2')$. Our hypotheses on $\pi_1(M(r_2))$ now imply that $M_1(r_2')$ is simply connected, and as $M_1$ is Seifert-fibered, it is in fact $S^3$. This completes the proof of (4). □

LEMMA 9.5. *Suppose that $N$ is simple and the image in $M(r_1)$ of $T_1 \cup T_2 \cup \ldots \cup T_m$ is a set of essential tori. Then each closed, essential surface in $N$ compresses in $N(\partial M; r_1)$.*

*Proof.* Suppose that $S$ is a closed, essential surface in $N$. As $N$ is simple, the genus of $S$ is at least 2. Hence if $S$ remains essential in the Seifert manifold $M(r_1)$, it must be isotopic to a horizontal surface [24, VI.34]. But the image in $M(r_1)$ of the torus $T_1$, call it $T_1'$, is essential and disjoint from $S$, and therefore since $S$ is isotopic to a horizontal surface in $M(r_1)$, Lemma 2.8 (2) implies that $T_1'$ is isotopic to $S$ in $M(r_1)$. Thus the genus of $S$ is 1, contrary to the fact that $N$ is simple. Thus $S$ must compress in $M(r_1)$. Finally note that since the image of $\partial N \setminus \partial M$ is essential in $M(r_1)$, $S$ must compress in $N(\partial M; r_1)$. □

LEMMA 9.6. *Suppose that $N$ is simple, non-Seifert and that each component of $\partial N \setminus \partial M$ corresponds to an essential torus in $M(r_1)$. Then $\Delta(r_1, r_2) \leq 1$.*

*Proof.* The idea behind the proof is to use Theorem 1.5 to show that there are slopes on the components of $\partial N \setminus \partial M$ such that by filling $N$ along these slopes we obtain a simple, non-Seifert manifold $Q$ for which $Q(r_1)$ has a cyclic fundamental group while $Q(r_2)$ is either reducible or has a cyclic fundamental group. In either event, Corollary 1.4 will now imply that $\Delta(r_1, r_2) \leq 1$.

We shall appeal to the following two results in the proof.

(9.7) [43, Th. 2.6] There are only finitely many slopes on a toral boundary component of a compact, orientable, simple, non-Seifert 3-manifold which yield either a nonsimple manifold or a Seifert manifold.

(9.8) [21] There are only finitely many slopes on a toral boundary component $T$ of a compact, irreducible 3-manifold, corresponding to the boundary of an incompressible, boundary-incompressible surface whose boundary lies entirely in $T$.

Recall from Lemma 9.3 (2) that $N(\partial M; r_2) \cong (S^1 \times D^2) \# W$ for some 3-manifold $W$ whose boundary is $T_2 \cup T_3 \cup \ldots \cup T_m$. It was shown in Lemma 9.4 (2) that $N(\partial M; r_1)$ and $\overline{M \setminus N_M}$ admit Seifert fibrations. Fix such structures and let $\phi_i$ be the slope on $T_i$ corresponding to the fiber of the structure on $N(\partial M; r_1)$. If we assume that the image in $M(r_1)$ of the components of $\partial N \setminus \partial M$ is isotopic to a collection of vertical tori, then the Seifert structure on $M(r_1)$ pulls back to one on $N(\partial M; r_1)$ and restricts to one on $\overline{M \setminus N_M}$. We shall assume in this case that these are the Seifert structures which have been chosen.

Fix a collection of positive integers $k_1, k_2, \ldots, k_m \geq 1$ and select slopes $s'_m, s'_{m-1}, \ldots, s'_1$ on $T_m, T_{m-1}, \ldots, T_1$ inductively as follows. Set $N_{m+1} = N$ and let $1 \leq j \leq m$. Suppose that a set of slopes $s'_m, s'_{m-1}, \ldots, s'_{j+1}$ have been chosen (the empty set when $j = m$) so that

- $\Delta(s'_i, \phi_i) = k_i$, $j+1 \leq i \leq m$, and

- $N_j$, the $(s'_{j+1}, s'_{j+2}, \ldots, s'_m)$-filling of $N$, is simple, non-Seifert and any closed, essential surface in $N_j$ is isotopic to one in $N_{j+1}$.

According to (9.7) and (9.8), there are only finitely many slopes $r$ on $T_j$ for which $N_{j+1}(T_j, r)$ is either not simple or Seifert-fibered, or for which $N_{j+1}(T_j, r)$ contains a closed, essential surface which is not isotopic into $N_{j+1}$. Select $s'_j$ to be any slope on $T_j$ not amongst this finite number of exceptions and which also satisfies $\Delta(s'_j, \phi_j) = k_j$. Now proceed inductively to construct $s'_m, s'_{m-1}, \ldots, s'_1$.

In what follows, we shall reserve the notation $s'_i$ for any slope on $T_i$ which has been chosen by such an inductive process.

Denote by $Q$ the $(s'_1, s'_2, \ldots, s'_m)$ filling of $N$. Then by construction $Q$ is simple and non-Seifert and any closed, essential surface in $Q$ is isotopic into $N$. Moreover, as $k_1, k_2, \ldots, k_m \geq 1$, the Seifert structure on $N(\partial M; r_1)$ extends over $Q(r_1)$ in such a way that the core of the solid torus attached to $T_i$ has multiplicity $k_i$.

Recall from Lemma 9.3 (2) that $N(\partial M; r_2) \cong (S^1 \times D^2) \# W$ and that $s_1$ is the slope of the meridian of the $S^1 \times D^2$ factor. If $W(s'_2, \ldots, s'_m)$ denotes the $(s'_2, \ldots, s'_m)$ filling of $W$, then

$$Q(r_2) \cong L_{\Delta(s_1, s'_1)} \# W(s'_2, \ldots, s'_m)$$

where $L_{\Delta(s_1, s'_1)}$ is a lens space whose first homology group has order $\Delta(s_1, s'_1)$. We remark that the value of $\Delta(s_1, s'_1)$ is not independent of $k_1$; for instance if $s_1 = \phi_1$ they are equal. Nevertheless, we can always find a value of $k_1$ for which there is a slope $s'_1$ on $T_1$ satisfying $\Delta(s'_1, s_1) > 1$ say, or $\Delta(s'_1, s_1) = 1$.

*Claim* 1. If $N$ contains a closed, essential, nonseparating surface, then $\Delta(r_1, r_2) \leq 1$.

*Proof of Claim* 1. Amongst all closed, essential, nonseparating surfaces in $N$, we may choose one, $S$ say, whose genus $g$ is minimal. Since $N$ is simple, $g \geq 2$. Observe that $S$ must remain incompressible in $Q$, for otherwise we could perform a sequence of compressions on $S$ in $Q$ to produce a closed, essential, nonseparating surface $S_0$ whose genus would be strictly less than $g$. But by the choices made in the construction of $Q$ from $N$, $S_0$ is isotopic to a surface in $N$, and so its genus is at least $g$. Thus $S$ must remain essential in $Q$. A similar argument shows that the genus of any closed, nonseparating surface in $Q$ is greater than or equal to $g$, and so as $Q$ is simple, $S$ is Thurston norm minimizing there.

Next we note that since $Q$ contains a closed, nonseparating surface, rank $H_2(Q) > 0$. As the boundary of $Q$ is a torus, its Euler characteristic is zero, and hence rank $H_1(Q) = $ rank $H_2(Q) + $ rank $H_0(Q) > 1$. We may now apply [10, p. 462, Cor.] to see that there is at most one slope $r$ on $\partial Q$ such that $S$ compresses in $Q(r)$. Now according to Lemma 9.5, $S$ compresses in $N(\partial M; r_1)$, and hence also in $Q(r_1)$. Thus $S$ remains incompressible in $Q(r_2)$. We shall see that this is impossible if $\Delta(r_1, r_2) > 1$.

Let us suppose that $\Delta(r_1, r_2) > 1$. Consider the manifold $P$ obtained by cutting $N$ open along $S$. By Lemma 9.5, one of the two copies $S', S''$ of $S$ in $\partial P$, say $S'$, compresses in $P(\partial M; r_1)$. On the other hand, $T_1$ compresses in $N(\partial M; r_2)$ and so either it compresses in $N(\partial M; r_2) \setminus S$ or $S$ compresses in $N(\partial M; r_2)$. In either event $P(\partial M; r_2)$ is $\partial$-reducible. Since $\Delta(r_1, r_2) > 1$, [6, Th. 2.4.5] shows that there is a unique component of $\partial P \setminus \partial M$ which compresses in both $P(\partial M; r_1)$ and $P(\partial M; r_2)$. We have assumed that the image of $\partial N \setminus \partial M$ is essential in $M(r_1)$, and so it is incompressible in $P(\partial M; r_1)$. Thus $S'$ or $S''$ must compress in $P(\partial M; r_2)$, which implies that $S$ compresses in $N(\partial M; r_2) \subset Q(r_2)$. But this contradicts what was deduced in the previous paragraph. Thus we must have $\Delta(r_1, r_2) \leq 1$. □

We may therefore suppose for the remainder of the proof of Lemma 9.6 that $N$ does not contain a closed, essential, nonseparating surface.

*Claim* 2. Either $\Delta(r_1, r_2) \leq 1$ or $W(s'_2, \ldots, s'_m)$ is a possibly trivial lens space.

*Proof of Claim* 2. Without loss of generality we may take $W(s'_2, \ldots, s'_m) \neq S^3$. As observed above there is a value of $k_1$ for which we can choose a slope $s'_1$ satisfying $\Delta(s_1, s'_1) > 1$. In this case $Q(r_2) \cong L_{\Delta(s_1, s'_1)} \# W(s'_2, \ldots, s'_m)$ is a reducible manifold and therefore $r_2$ is a boundary slope on $\partial Q$. If rank $H_1(Q) > 1$, then rank $H_2(Q) = $ rank $H_1(Q) - $ rank $H_0(Q) > 0$. Thus $Q$ contains a closed, essential, nonseparating surface, and therefore our construction of $Q$ implies that the same is true for $N$, contrary to our hypotheses. Hence rank $H_1(Q) = 1$, and so we may apply Corollary 2.2 to see that either (i) $Q(r_2)$ is a connected sum of two nontrivial lens spaces, or (ii) $Q$ contains a closed essential surface $S$ such that $S$ is compressible in $Q(r_2)$ and is incompressible in $Q(r)$ whenever $\Delta(r, r_2) > 1$, or (iii) $Q(r_2) \cong S^1 \times S^2$. The third possibility is inconsistent with our choice of $s'_1$. If possibility (ii) occurs then by the construction of $Q$, the surface $S$ may be assumed to lie in $N$. But then according to Lemma 9.5, $S$ compresses in $N(\partial M; r_1)$ and hence in $Q(r_1)$. Thus $\Delta(r_1, r_2) \leq 1$ by (ii). We may therefore assume that possibility (i) arises, and so $W(s_2', \ldots, s_m')$ is a lens space.

*Case* 1. $m > 1$. Choose $k_2 > 1$. Next choose a value of $k_1$ for which there is a slope $s_1'$ satisfying $\Delta(s_1, s_1') = 1$ and the manifold $N(s_1', s_3', s_4', \ldots, s_m')$ is simple, non-Seifert. Then $Q(r_2) \cong W(s_2', \ldots, s_m')$ and so by Claim 2 we may assume that $r_2$ is a cyclic filling slope on $\partial Q$. On the other hand $r_1$ is a Seifert filling slope on $\partial Q$. According to Theorem 1.5, either $\Delta(r_1, r_2) \leq 1$ or $Q(r_1)$ admits a Seifert stucture whose base orbifold is of the form $S^2(p, q, r)$ where $2 \leq p \leq q \leq r$. Suppose that the latter case occurs. Then by [24, VI.17], either $Q(r_1)$ has a unique Seifert structure, or $(p, q, r) = (2, 2, n)$ for some $n \geq 2$ and $Q(r_1)$ admits exactly one other Seifert structure whose base orbifold is of the form $\mathbf{R}P^2(a)$ for some $a \geq 2$. Now as we noted above, the Seifert structure on $N(\partial M; r_1)$ extends over $Q(r_1)$ in such a way that there is an exceptional fiber of multiplicity $k_2 > 1$ corresponding to the core of the solid torus attached to $T_2$. Thus consideration of the possible Seifert structures on $Q(r_1)$ shows that $N'$, the $(r_1, s_1', s_3', s_4', \ldots, s_m')$ filling of $N$, is a Seifert-fibered manifold whose base orbifold is either a Moebius band without cone points or a 2-disk with exactly two cone points. In either event, $N'$ admits infinitely many distinct cyclic filling slopes on $T_2$. We may therefore find some new $k_2$ and $s_2'$ so that the new $Q = N(s_1', s_2', s_3', \ldots, s_m')$ is both simple and non-Seifert (recall that $s_1'$ was selected so that $N(s_1', s_3', s_4', \ldots, s_m')$ is simple, non-Seifert) and for which $r_1$ becomes a cyclic filling slope on $\partial Q$. Since $r_2$ is also a cyclic filling slope on $\partial Q$, the cyclic surgery theorem [6] implies that $\Delta(r_1, r_2) \leq 1$. Thus Lemma 9.6 holds when $m > 1$.

*Case* 2. $m = 1$. Now as $m = 1$, $\partial N = \partial M \cup T_1$ and so $T_1$ separates $M$ into two pieces $M = N \cup_{T_1} M_1$. According to Lemma 9.4 (2), $M_1$ is a Seifert-fibered space. We also note that we may suppose that $W$ is a possibly trivial lens space (see Claim 2 above). We claim that either $W \cong S^3$ or $M_1(s_1) \cong S^3$. To see this, recall that from the definition of the slope $s_1$ on $T_1$ we have a decomposition $M(r_2) = M_1(s_1) \# W$. Since $\pi_1(M(r_2))$ is finite or cyclic, one of $\pi_1(W)$ and $\pi_1(M_1(s_1))$ is trivial. But as $W$ is a possibly trivial lens space and $M_1$ admits a Seifert fibering, it follows that either $W \cong S^3$ or $M_1(s_1) \cong S^3$. We shall consider these two cases separately.

Suppose first of all that $W \cong S^3$. Then any filling of $N(\partial M; r_2) \cong (S^1 \times D^2) \# W \cong S^1 \times D^2$ is a cyclic filling. If we choose $k_1 > 1$, then $Q(r_2)$ has a cyclic fundamental group while $Q(r_1)$ is Seifert with an exceptional fiber of multiplicity $k_1 > 1$ corresponding to the core of solid torus attached to $T_1$. We may invoke Theorem 1.5 as in the case $m > 1$ to deduce that either $\Delta(r_1, r_2) \leq 1$ or $N(\partial M; r_1)$ admits a Seifert structure whose base orbifold is a Moebius band without cone points or a 2-disk with exactly two cone points. In either event we can now rechoose $k_1$ and $s_1'$ so that both $r_1$ and $r_2$ are cyclic filling slopes of the (new) simple, non-Seifert manifold $Q$. The cyclic surgery theorem now implies that $\Delta(r_1, r_2) \leq 1$.

Finally assume that $M_1(s_1) \cong S^3$. We claim that in this case $s_1 \neq \phi_1$. To see this note that $M_1$ is a torus knot complement and so in particular $M_1$ is not a union of copies of $I(K)$. It now follows from Lemma 9.4 (3) that the image of $T_1$ in $M(r_1)$ is isotopic to a vertical torus. Thus the Seifert structures on $N(\partial M; r_1)$ and $M_1$ are restrictions of the one on $M(r_1)$. It follows that $\phi_1$ is also the fiber of a Seifert structure on the torus knot exterior $M_1$. Hence $M_1(\phi_1)$ is a connected sum of two nontrivial lens spaces, and so in particular $S^3 \cong M_1(s_1) \not\cong M_1(\phi_1)$. Thus $s_1 \neq \phi_1$. But then we may find a $k_1$ and subsequently choose $s_1'$ so that $\Delta(s_1', \phi_1) = k_1 > 1$ and $\Delta(s_1', s_1) = 1$. It follows that the Seifert structure on $N(\partial M; r_1)$ extends over $Q(r_1)$ in such a way that there is an exceptional fiber of multiplicity $k_1 > 1$ corresponding to the core of the solid torus attached to $T_1$. On the other hand, for this choice of $s_1'$ we have $Q(r_2) = W$ and therefore by Claim 2, $r_2$ is a cyclic filling slope on $\partial Q$. As in our previous cases, either $\Delta(r_1, r_2) \leq 1$ or $N(\partial M; r_1)$ admits a Seifert structure whose base orbifold is a Moebius band without cone points or a 2-disk with exactly two cone points. Hence there are new choices of $k_1$ and $s_1'$ for which $r_1$ becomes a cyclic filling slope on the boundary of the new $Q$, while $Q(r_2) \cong L_{\Delta(s_1, s_1')} \# W$ is either reducible or has a cyclic fundamental group. In either case, Corollary 1.4 implies that $\Delta(r_1, r_2) \leq 1$. We have therefore completed the proof of Lemma 9.6. □

*Proof of* (9.1) *and* (9.2). The proof will proceed through the consideration of several cases.

*Case* 1. The manifold $M$ is not cabled. If some component of $\partial N \setminus \partial M$ compresses in $M(r_1)$ then by Lemma 9.4 (1) we see that $\Delta(r_1, r_2) \leq 1$. Thus we may assume that $\partial N \setminus \partial M$ is sent to a collection of essential tori in $M(r_1)$. If $N$ is simple and non-Seifert then Lemma 9.6 implies $\Delta(r_1, r_2) \leq 1$, while if $N$ is Seifert-fibered, Lemma 9.4 (4) gives the same result. Thus both (9.1) and (9.2) hold.

We shall assume for the remainder of the proofs of (9.1) and (9.2) that $M$ is cabled. Write $M = C_1 \cup_T M_1$ where $C_1$ is a cable space, $\partial C_1 = \partial M \cup T$, and $\partial M_1 = T$. According to Lemma 9.3 (1), there is a slope $r_2'$ on $\partial M_1$ such that $M(r_2) \cong M_1(r_2') \# V$ for some closed 3-manifold $V$. It follows that

$$(9.9) \qquad \pi_1(M_1(r_2')) \text{ is } \begin{cases} \text{finite if } \pi_1(M(r_2)) \text{ is finite} \\ \text{cyclic if } \pi_1(M(r_2)) \text{ is cyclic.} \end{cases}$$

*Case* 2. The manifold $M_1$ is a Seifert-fibered space. If $M_1$ is a Seifert-fibered manifold then (9.9) shows that both (9.1) and (9.2) hold.

*Case* 3. All the tori in $\partial N \setminus \partial M$ correspond to essential tori in $M(r_1)$.

As we have assumed that $M$ is cabled, $N$ must be Seifert-fibered. An application of Lemma 9.4 (4) now completes this case.

Appealing to Cases (2) and (3) shows that we may assume for the rest of the proof that the image $T' \subset M(r_1)$ of some torus in $T_i \subset \partial N \setminus \partial M$ is compressible. Furthermore, we may suppose that $\Delta(r_1, r_2) > 1$. Then as the image of $T_i$ in $M(r_2)$ also compresses, it follows from [6, Th. 2.0.1] that $T'$ and $\partial M$ cobound a cable space in $M$. Without loss of generality we shall suppose that $T' = T$ and $N = C_1$. According to [12, Lemma 7.2], our hypotheses imply that $C_1(\partial M; r_1) \cong C_1(\partial M; r_2) \cong S^1 \times D^2$ and there are slopes $r_1', r_2'$ on $\partial M_1$ such that $M_1(r_1') \cong M(r_1)$, $M_1(r_2') \cong M(r_2)$. Further, if we suppose that $C_1$ is a cable space of type $(m, n)$, $n \geq 2$, then $\Delta(r_1', r_2') = n^2 \Delta(r_1, r_2) \geq 8$ and $\Delta(r_1', \phi_1) = \Delta(r_1', \phi_1) = n$ where $\phi_1$ is the slope on $T$ corresponding to the Seifert structure on $C_1$ ([12, Lemma 3.3]).

*Case* 4. The manifold $M_1$ is simple, non-Seifert, $\Delta(r_1, r_2) > 1$, and the image of $\partial N \setminus \partial M$ compresses in $M(r_1)$.

In this case $M$ is cable on a simple manifold and so by (9.9), (9.1) holds. Assume now that $r_2$ is a cyclic filling slope and that $M(r_1)$ admits no Seifert fibering whose base orbifold is the 2-sphere with exactly three cone points. Now $M(r_1) \cong M_1(r_1')$ and $M(r_2) \cong M_1(r_2')$ and so by Theorem 1.5 (1), $\Delta(r_1', r_2') \leq 1$. But by construction $\Delta(r_1', r_2') \geq 8$. Hence under the hypotheses of (9.2), this case does not arise.

*Case* 5. The manifold $M_1$ is neither simple nor Seifert, $\Delta(r_1, r_2) > 1$, and the image of $\partial N \setminus \partial M$ compresses in $M(r_1)$.

As $M_1$ is not Seifert fibered and contains an essential torus, we can apply Case 1, replacing $r_1$ by $r_1'$ and $r_2$ by $r_2'$, to conclude that $M_1$ is cabled. Write $M_1 = C_2 \cup_{T''} M_2$ where $C_2$ is a cable space, $\partial C_2 = \partial M_1 \cup T''$, and $\partial M_2 = T''$. Notice that as $\Delta(r_1, r_2) > 1$ and $T''$ compresses in $M(r_2)$, it cannot compress in $M(r_1)$ since otherwise [6, Th. 2.0.1] would imply that $C_1 \cup_T C_2$ is a cable space, which is clearly false. Thus $T''$ remains essential in $M(r_1)$. But then $M_2$ must be a Seifert-fibered space, for the Seifert structure on $M(r_1)$ restricts to one on $M_2$ if $T''$ is isotopic to a vertical torus in $M(r_1)$, while $M_2 \cong I(K)$ if $T''$ is isotopic to a horizontal torus in $M(r_1)$ [24, VI.34]. Since $M_1$ is not Seifert-fibered, it follows that $C_2$ is the piece of the torus decomposition of $M_1$ which contains $T$.

Consider $\phi_1'$, the slope on $T$ corresponding to the fiber slope of the Seifert structure on $C_2$. As $T''$ compresses in $M_1(r_2') \cong M(r_2)$, we have $\Delta(r_2', \phi_1') \leq 1$ ([12, Lemma 7.2]). On the other hand since we have assumed that $\Delta(r'_1, r'_2) > 1$, we may apply Lemma 9.4 (4) to $M_1$ and $C_2$ to deduce that $\Delta(r'_1, \phi'_1) = 2$. Thus we have four slopes $r'_1, r'_2, \phi_1, \phi'_1$ on $T$ and they satisfy the equations $\Delta(r'_1, \phi_1) = n$, $\Delta(r'_2, \phi_1) = n$, $\Delta(r'_1, \phi'_1) = 2$, and $\Delta(r'_2, \phi'_1) \leq 1$. The argument given in the last paragraph of the proof of Lemma 2.4 may be used to deduce that $n\Delta(r_1, r_2)\Delta(\phi_1, \phi'_1)$ is either 1, 2, or 3. We have selected $n$ to be at least 2, and so $\Delta(r_1, r_2) = 1$, contrary to our assumptions. Hence this case does not arise under the hypotheses of either (9.1) or (9.2).

This completes the proofs of (9.1) and (9.2), and therefore the proof of Theorem 1.6. □

We shall finish this section by giving two constructions of cablings on a Seifert manifold $M_1$ which admit Seifert and finite filling slopes arbitrarily distanced from each other. The following notation will be common to both.

- $\phi_1$ is an *oriented* slope on $\partial M_1$ corresponding to the fiber of the Seifert structure on $M_1$ and $\beta \in H_1(\partial M_1)$ is chosen so that it forms a basis of $H_1(\partial M_1)$ along with the homology class $\alpha(\phi_1)$ carried by $\phi_1$.

- $C$ is a cable space of type $(m,n)$, $n \geq 2$, whose boundary consists of two tori $T_+, T_-$. We use $\phi_+$ and $\phi_-$ to denote oriented slopes on $T_+$ and $T_-$ corresponding to the fiber of the (unique) Seifert structure on $C$ and we use $\alpha(\phi_\pm)$ to denote the associated primitive homology classes in $H_1(T_\pm)$.

- Choose bases $\{\mu_\pm, \lambda_\pm\}$ for $H_1(T_\pm)$ according to the prescriptions from [16]:
$$\alpha(\phi_+) = mn\mu_+ + \lambda_+ \text{ and } \alpha(\phi_-) = m\mu_- + n\lambda_-.$$

It is shown in [12] that if $k \in \mathbf{Z}$ and $r$ is the slope on $T_+$ for which $\alpha(r) = \pm((kmn+1)\mu_+ + k\lambda_+)$ (i.e. any slope $r$ which satisfies $\Delta(r, \phi_+) = 1$), then

(9.10)
$$C(T_+; r) \cong S^1 \times D^2 \text{ where } \alpha(\{*\} \times \partial D^2) = \pm((kmn+1)\mu_- + kn^2\lambda_-) \in H_1(T_-).$$

*Example* 9.11. Take $M_1$ to be a Seifert-fibered manifold whose base orbifold is $D^2(2,3,5)$. Any filling of $M_1$ whose slope $r$ satisfies $\Delta(r, \phi_1) = 1$ is a finite filling, while any filling whose slope $r$ satisfies $\Delta(r, \phi_1) > 1$ is Seifert fibered in such a way as to admit an essential, vertical torus. Define $M$ by gluing $C$ and $M_1$ together along $T_-$ and $\partial M_1$ so that $\mu_-$ corresponds to $\beta$ and $\lambda_-$ corresponds to $\alpha(\phi_1)$. Now $M$ is not Seifert, for if it were, $T_-$ would necessarily be a vertical torus, and so the structure would restrict to one on both $C$ and $M_1$. By the uniqueness of the Seifert structure on $C$, we would then have $\alpha(\phi_-)$ identified with $\pm\alpha(\phi_1)$ under the gluing map, i.e. $\pm\alpha(\phi_1) = m\beta + n\alpha(\phi_1)$, which contradicts the fact that $n \geq 2$. Thus $M$ is not Seifert.

Fix $k \in \mathbf{Z}$ and choose $r_1$ and $r_2$ to be the oriented slopes satisfying $\alpha(r_1) = (kmn + 1)\mu_+ + k\lambda_+$ and $\alpha(r_2) = \mu_+$. Using (9.10) above we have $M(r_1) = M_1(r_1')$ where $\pm\alpha(r_1') = (kmn + 1)\beta + kn^2\alpha(\phi_1)$ and $M(r_2) = M_1(r_2')$ where $\pm\alpha(r_2') = \beta$. Now $\Delta(r_1', \phi_1) = |kmn + 1|$ and so as noted above, $M(r_1)$ will be Seifert with an essential, vertical torus as long as $|kmn + 1| > 1$. Also $\Delta(r_2', \phi_1) = 1$ implies that $M(r_2) = M_1(r_2')$ by (9.10) above, and so $r_2$ is a finite filling slope. Finally, $\Delta(r_1, r_2) = |k|$, which can be made arbitrarily large. □

*Example* 9.12. For this example let $C$ be a cable space of type $(1, 2)$, let $M_1$ have base orbifold $D^2(2, 2)$, and let $\phi_1$ correspond to this fibering of $M_1 \cong I(K)$. There is only one other slope $\phi_1'$ on $\partial M_1$ which corresponds to the fiber of some other Seifert fibering of $M_1$ [24, VI.17]. We remark that any filling of $M_1$ whose slope $r$ satisfies $\Delta(r, \phi_1) > 1$ is a finite, noncyclic filling. We shall suppose that $\beta \neq \pm\alpha(\phi_1')$. Fix an integer $k$ and define $M = C \cup_{T_- = \partial M_1} M_1$ glued together so that $\mu_-$ corresponds to $(2k+1)\beta - 2\alpha(\phi_1)$ and $\lambda_-$ corresponds to $-k\beta + \alpha(\phi_1)$. From our choice of $\beta$ we may argue as in the previous example that $M$ is not Seifert-fibered.

Let $r_1$ be the oriented slope on $\partial M$ for which $\alpha(r_1) = 4k\mu_+ + (2k + 1)\lambda_+$ and note that $\Delta(r_1, \phi_+) = 2$. Thus $C(\partial M; r_1)$ has base orbifold $D^2(2, 2)$, i.e. $C(\partial M; r_1) \cong I(K)$. It can be verified that the slope $k\mu_- + (2k + 1)\lambda_-$ is the slope on $C(\partial M; r_1)$ which yields $S^1 \times S^2$ upon filling, and therefore it is the slope of the other Seifert fibering of $C(\partial M; r_1)$. But $k\mu_- + (2k + 1)\lambda_-$ corresponds to $\alpha(\phi_1)$ under our gluing map, and so $M(r_1)$ is Seifert-fibered with base orbifold $\mathbf{R}P^2(2, 2)$. Now choose an oriented slope $r_2$ on $\partial M$ so that $\alpha(r_2) = \mu_+$ and note that $\Delta(r_2, \phi_+) = 1$ so that $C(\partial M; r_2) \cong S^1 \times D^2$ and has meridian curve whose homology class on $T_-$ is $\pm\mu_-$. By construction, the latter curve has distance $|2k+1|$ from $\phi_1$ and so $M(r_2)$ has a finite fundamental group, as long as $|2k + 1| \geq 2$. Finally note that $\Delta(r_1, r_2) = |2k + 1|$ can be made arbitrarily large. □

## 10. Proof of Theorem 1.8

Theorem 1.8 is a consequence of the following lemma and example.

LEMMA 10.1. *Let $K$ be a small hyperbolic knot in $S^3$. Suppose that $K$ has a boundary slope $r = m/n$ with $n > 1$. Then $M(r)$ is a Haken manifold with 0-dimensional $\mathrm{PSL}_2(\mathbf{C})$-character variety. If we further assume that $n > 2$, then $M(r)$ is also a hyperbolic manifold.*

*Remark.* Note that the inequality $\dim X(M(r)) \leq \dim \bar{X}(M(r)) = 0$ shows that $X(M(r))$ is also 0-dimensional.

*Proof of Lemma* 10.1. By the main result of [17], or Corollary 1.3, $M(r)$ is an irreducible manifold and so we may apply Lemma 2.1 to see that $M(r)$ is Haken. Suppose that $\dim \bar{X}(M(r)) > 0$ and choose a curve $X_0 \subset \bar{X}(M(r)) \subset \bar{X}(M)$. Now Proposition 5.7 (1) shows that $\|\cdot\|_{X_0}$ is nontrivial, while $f_r \equiv 0$ on $X_0$ by construction. Thus $X_0$ must be an $r$-curve. But this is impossible according to the conclusion of Corollary 6.7. Hence $\dim \bar{X}(M(r)) = 0$

Now consider the case where $n > 2$. Gordon and Luecke have shown that if $r'$ is a slope on the boundary of a hyperbolic knot exterior, then $M(r')$ is simple as long as $\Delta(r', \mu) > 2$ ([19], [20]). Thus in our case, $M(m/n)$ is a simple, Haken manifold. Further note that by Corollary 1.7, $M(m/n)$ cannot be Seifert-fibered. Therefore according to Thurston [43, Th. 2.5], $M(m/n)$ is a hyperbolic manifold. □

The following example gives an infinite family of knots in $S^3$ satisfying the conditions of Lemma 10.1.

*Example* 10.2. Let $K_m$ be the pretzel knot of type $(-1/3, 1/(2m+1), 5/18)$. Then for $|n| \geq 3$, the knot $K_{4n+6}$ satisfies the conditions of Lemma 10.1.

*Proof.* First of all we shall apply the results of [22] to show that for each $n$, $r = \frac{16n^2+22n+1}{n}$ is a boundary slope of $K_{4n+6}$ (cf. [22, Prop. 2.2]). We shall assume the reader is familiar with the notation and terminology of [22].

Let $(\gamma_1, \gamma_2, \gamma_3)$ be the edge-path system of $K_m$ where $\gamma_1$ is the constant path at the point $(m-6, 2m, 2-m)$, $\gamma_2$ is the edge path $< 1/(2m+1) >$ $\to (m-7) < 1/0 > + < 1/(2m+1) >$, and $\gamma_3$ is the edge path $< 5/18 >$ $\to < 2/7 > \to (m-9) < 1/3 > + 3 < 2/7 >$. If $S$ is the candidate surface determined by $(\gamma_1, \gamma_2, \gamma_3)$, then by [22, Prop. 2.1], $S$ is incompressible. We must calculate the boundary slope of $S$.

Now
$$\tau(S) = 2[0 + (m-7)/(m-6) - 1 - (m-9)/(m-6)] = -2 + 4/(m-6).$$

To compute $\tau(S_0)$ for a Seifert surface $S_0$, use the edge path system

$$< -1/3 > \to < -1/2 > \to < -1 > \to < 1/0 >,$$
$$< 1/(2m+1) > \to < 1/2m > \to < 1/(2m-1) >$$
$$\to \ldots \to < 1/2 > \to < 1 > \to < 1/0 >,$$

and
$$< 5/18 > \to < 2/7 > \to < 1/4 > \to < 0 > \to < 1/0 >,$$

so that $\tau(S_0) = 2[1 - (2m+1) - 0] = -4m$. Thus the boundary slope $r_m$ of $S$ is $r_m = \tau(S) - \tau(S_0) = \frac{4}{m-6} - 2 + 4m$. Taking $m = 4n+6$, we have $r_{4n+6} = \frac{16n^2+22n+1}{n}$. Finally, we may apply [36] to see that $K_{4n+6}$ is small and hyperbolic. □

*Proof of Theorem* 1.8. Let $M_n(r_{4n+6})$ denote the manifold obtained by Dehn surgery on $K_{4n+6}$ with slope $r_{4n+6} == \frac{16n^2+22n+1}{n}$. By Lemma 10.1, $M_n(r_{4n+6})$ will be a hyperbolic Haken manifold with 0-dimensional character variety as long as $n > 2$. Finally $H_1(M_n(r_{4n+6})) = \mathbf{Z}/16n^2 + 22n + 1$, and so the manifolds $M_n(r_{4n+6})$, $n \geq 3$ are mutually distinct. □

Question 10.3. Is there a hyperbolic knot in $S^3$ which is not small such that the $\mathrm{PSL}_2(\mathbf{C})$-character variety of its exterior only has 1-dimensional components?

UNIVERSITÉ DU QUÉBEC À MONTRÉAL, MONTRÉAL, CANADA
*E-mail address*: boyer@math.uqam.ca
*E-mail address*: zhang@math.uqam.ca